\documentclass[11pt,twoside]{article}
\newcommand{\ifims}[2]{#1} 
\newcommand{\ifAMS}[2]{#1}   
\newcommand{\ifau}[4]{#1}  
\newcommand{\ifbook}[2]{#1}   
\usepackage{ifthen}
\usepackage[ampersand]{easylist}
\usepackage{amsmath,amssymb,amsthm}
\usepackage{mathtools}
\usepackage{natbib}
\usepackage{epsfig,graphicx}
\usepackage{comment}
\usepackage{color}
\usepackage{srcltx}
\usepackage[mathscr]{eucal}
\usepackage[math]{easyeqn}
\usepackage{etoolbox}
\usepackage{hyperref}
\hypersetup{
            colorlinks,
            linkcolor=hookersgreen,
            linktoc=blue,
            citecolor=ultramarine,
            urlcolor=black,
            filecolor=black
            }
\usepackage{nameref}
\usepackage{multirow}

\ifims{
\textheight=22cm
\textwidth=14.8cm
\topmargin=0pt
\oddsidemargin=1.0cm
\evensidemargin=1.0cm
\linespread{1.3}

}{ 
}

\numberwithin{equation}{section}
\numberwithin{figure}{section}
\newcounter{example}[section]
\numberwithin{example}{section}
\newcounter{remark}[section]
\numberwithin{remark}{section}
\newtheorem{theorem}{Theorem}[section]

\newtheorem{lemma}[theorem]{Lemma}
\newtheorem{corollary}[theorem]{Corollary}

\newtheorem{exmp}[example]{Example}
\newtheorem{rmrk}[remark]{Remark}
\newenvironment{example}{\begin{exmp}\rm}{\end{exmp}}
\newenvironment{remark}{\begin{rmrk}\rm}{\end{rmrk}}

\ifbook{

  }
  {

  }

\renewcommand{\(}{$\,}
\renewcommand{\)}{\,$}

\def\nquad{\hspace{-1cm}}
\def\eqdef{\stackrel{\operatorname{def}}{=}}

\def\tow{\stackrel{w}{\longrightarrow}}

\def\block{\operatorname{block}}

\newcommand{\cc}[1]{\mathscr{#1}}
\newcommand{\bb}[1]{\boldsymbol{#1}}

\renewcommand{\bar}[1]{\overline{#1}}
\renewcommand{\hat}[1]{\widehat{#1}}
\renewcommand{\tilde}[1]{\widetilde{#1}}

\newcommand{\thankstitle}[1]{\ifthenelse{\equal{#1}{}}{}{\thanks{#1}}}
\newcommand{\thanksau}[1]{\ifthenelse{\equal{#1}{}}{}{\thanks{#1}}}

\newcommand{\aua}[6]
{\def\authora{#1}
\def\runauthora{#2}
\def\addressa{#3}
\def\emaila{#4}
\def\affiliationa{#5}
\def\thanksa{#6}}

\def\theauthors{
\ifau{ 
  \author{
    \authora
    \thanksau{\thanksa}
    \\[5.pt]
    \addressa \\
    \texttt{ \emaila}
  }
}
{  
  \author{
    \authora
    \thanksau{\thanksa}
    \\[5.pt]
    \addressa \\
    \texttt{ \emaila}
    \and
    \authorb
    \thanksau{\thanksb}
    \\[5.pt]
    \addressb \\
    \texttt{ \emailb}
  }
}
{   
  \author{
    \authora
    \thanksau{\thanksa}
    \\[5.pt]
    \addressa \\
    \texttt{ \emaila}
    \and
    \authorb
    \thanksau{\thanksb}
    \\[5.pt]
    \addressb \\
    \texttt{ \emailb}
    \and
    \authorc
    \thanksau{\thanksc}
    \\[5.pt]
    \addressc \\
    \texttt{ \emailc}
  }
} {   
  \author{
    \authora
    \thanksau{\thanksa}
    \\[5.pt]
    \addressa \\
    \texttt{ \emaila}
    \and
    \authorb
    \thanksau{\thanksb}
    \\[5.pt]
    \addressb \\
    \texttt{ \emailb}
    \and
    \authorc
    \thanksau{\thanksc}
    \\[5.pt]
    \addressc \\
    \texttt{ \emailc}
    \and
    \authord
    \thanksau{\thanksd}
    \\[5.pt]
    \addressd \\
    \texttt{ \emaild}
  }
}
}
\def\therunauthors{ \hfill \textsc{\small
 \ifau{\runauthora}
      {\runauthora and \runauthorb}
      {\runauthora, \runauthorb, and \runauthorc}
      {\runauthora, \runauthorb, \runauthorc, and \runauthord}
 }
 \hfill
 }

\renewcommand{\Gamma}{\varGamma}
\renewcommand{\Pi}{\varPi}
\renewcommand{\Sigma}{\varSigma}
\renewcommand{\Delta}{\varDelta}
\renewcommand{\Lambda}{\varLambda}
\renewcommand{\Psi}{\varPsi}
\renewcommand{\Phi}{\varPhi}
\renewcommand{\Theta}{\varTheta}
\renewcommand{\Omega}{\varOmega}
\renewcommand{\Xi}{\varXi}
\renewcommand{\Upsilon}{\varUpsilon}

\def\argmax{\operatornamewithlimits{argmax}}
\def\argmin{\operatornamewithlimits{argmin}}

\def\av{\bb{a}}

\def\fv{\bb{f}}
\def\gv{\bb{g}}

\def\uv{\bb{u}}

\def\wv{\bb{w}}
\def\xv{\bb{x}}

\def\Yv{\bb{Y}}

\def\alphav{\bb{\alpha}}
\def\betav{\bb{\beta}}

\def\deltav{\bb{\delta}}

\def\epsv{\bb{\varepsilon}}
\def\etav{\bb{\eta}}
\def\gammav{\bb{\gamma}}

\def\lambdav{\bb{\lambda}}

\def\psiv{\bb{\psi}}

\def\thetav{\bb{\theta}}
\def\upsilonv{\bb{\upsilon}}
\def\xiv{\bb{\xi}}

\usepackage{color}

\definecolor{blue(pigment)}{rgb}{0.2, 0.2, 0.6}
\definecolor{ultramarine}{rgb}{0.07, 0.04, 0.56}
\definecolor{darkspringgreen}{rgb}{0.09, 0.45, 0.27}
\definecolor{hookersgreen}{rgb}{0.0, 0.44, 0.0}
\definecolor{plum(traditional)}{rgb}{0.56, 0.27, 0.52}
\definecolor{purple(html/css)}{rgb}{0.5, 0.0, 0.5}
\definecolor{magenta(dye)}{rgb}{0.79, 0.08, 0.48}

\def\R{I\!\!R}
\def\E{I\!\!E}

\def\P{I\!\!P}

\def\kappa{\varkappa}

\def\Frobg{\Lambda}

\def\diag{\operatorname{diag}}

\def\Fr{\operatorname{Fr}}

\def\ND{\mathcal{N}}

\def\oper{\operatorname{op}}

\def\Var{\operatorname{Var}}

\def\T{\top}
\def\tr{\operatorname{tr}}

\def\CONST{\mathtt{C} \hspace{0.1em}}
\def\cond{\, \big| \,}

\def\ex{\mathrm{e}}

\def\Id{I\!\!\!I}
\def\Ind{\operatorname{1}\hspace{-4.3pt}\operatorname{I}}

\def\alp{\alpha}

\def\AA{A}

\def\BB{I\!\!B}     
\def\BB{B}

\def\CA{\cc{A}}

\def\CAt{\tilde{\CA}}
\def\CAGP{\CA_{\GP}}

\def\CONSTru{\CONST_{0}}

\def\DF{\cc{D}}

\def\DFGP{\DF_{\GP}}

\def\DFt{\tilde{\DF}}

\def\deltav{\bb{\delta}}
\def\deltavd{\deltav}

\def\dimp{p}

\def\dimtotal{\dimp^{*}}
\def\dimA{\mathtt{p}}
\def\dimAA{\dimA_{\tau}}
\def\dimB{\dimA}

\def\dimG{\dimA_{\GP}}

\def\dimq{q}

\def\dimt{\tilde{\dimA}}
\def\dimGt{\dimt_{\GP}}

\def\dimm{m}

\def\ELL{F}

\def\err{\diamondsuit}
\def\fvs{\fv}
\def\fvs{\fv^{*}}

\def\fs{f}

\def\fvs{\fv}

\def\gaussv{\bb{\gauss}}
\def\gauss{\gamma}

\def\gm{\mathtt{g}}
\def\gmc{\gm_{c}}
\def\gmb{\gm}

\def\gp{g}

\def\gvs{\gv^{*}}

\def\GP{G}

\def\GPY{\Gamma}

\def\HM{\mathbb{H}}

\def\IF{\Bbb{F}}

\def\KS{A}
\def\KSI{\mathfrak{\KS}}
\def\KSC{L}

\def\ksj{a}

\def\LL{\cc{L}}

\def\LLG{\LL_{\GP}}

\def\lambdaB{{\lambda}^{*}}		
\def\lambdaB{\lambda_{\BB}}		
\def\lambdaB{\supA}

\def\lambdav{\bb{\lambda}}

\def\nunu{\nu_{0}}

\def\PAAr{\Phi}

\def\PG{\P'}

\def\prior{\Pi}

\def\QP{Q}
\def\QPGP{\QP | \GP}

\def\rhot{t}

\def\rr{\mathtt{r}}

\def\rups{\rr_{0}}

\def\supA{\lambda}

\def\thetav{\bb{\theta}}

\def\TAU{\mathcal{T}}

\def\Thetad{\Theta^{\circ}}

\def\uvd{\uv^{\circ}}

\def\upsilonvs{\upsilonv^{*}}

\def\upsilonvc{\upsilonv'}
\def\vupsilonv{\upsilonv}

\def\Upsilond{\Upsilon_{0}}
\def\UpsilonD{\Upsilon^{d}}
\def\UpsilonC{\Upsilon^{\circ}}

\def\UV{\mathcal{U}}
\def\UVd{\UV^{\circ}}

\def\vA{\mathtt{v}}

\def\VF{\cc{V}}

\def\VP{V}

\def\VPGP{\VP | \GP}

\def\vp{\mathbf{v}}	

\def\vp{\mathsf{v}}

\def\vpB{\vp}

\def\wv{\bb{w}}
\def\wvd{\wv^{\circ}}

\def\xb{\bar{\xv}}

\def\xvd{\xv^{\circ}}
\def\xx{\mathtt{x}}
\def\xxc{\xx_{c}}

\def\xvd{\xv^{\circ}}

\def\XS{\cc{X}}
\def\XSD{\XS^{d}}

\def\Xs{X^{*}}

\def\YS{\cc{Y}}
\def\YSD{\YS^{d}}

\def\zq{z}
\def\zqc{\zq_{c}}

\def\thetitle{Bayesian inference for nonlinear inverse problems}
\def\theruntitle{Bayesian inference for nonlinear inverse problems}

\def\theabstract{
Bayesian methods are actively used for parameter identification and 
uncertainty quantification when solving nonlinear inverse problems
with random noise. 
However, there are only few theoretical results justifying the Bayesian approach.
Recent papers, see e.g. \cite{Nickl2017,lu2017bernsteinvon} and references therein,
illustrate the main difficulties and challenges in studying the properties of the 
posterior distribution in the nonparametric setup.
This paper offers a new approach for study the frequentist properties of the nonparametric 
Bayes procedures. 
The idea of the approach is to relax the nonlinear structural equation by introducing an auxiliary functional parameter 
and replacing the structural equation with a penalty and by imposing a prior on the auxiliary parameter. 
For the such extended model, we state sharp bounds on posterior concentration and on the accuracy of 
the penalized MLE and on
Gaussian approximation of the posterior, and a number of further results.
All the bounds are given in terms of effective dimension, and we show that the proposed 
calming device does not significantly affect this value. 
}

\def\kwdp{62F15}
\def\kwds{62F25}

\def\thekeywords{posterior, concentration, Gaussian approximation, Bernstein--von Mises Theorem}

\def\thankstitle{}

\aua
{Vladimir Spokoiny}
{Spokoiny, V. }
{
    Weierstrass Institute and HU Berlin,  \\
    HSE 
    and IITP RAS,
    \\
    Mohrenstr. 39, 10117 Berlin, Germany
    }
{spokoiny@wias-berlin.de}
{Weierstrass-Institute and Humboldt University Berlin}
{
  Financial support by the German Research Foundation (DFG) through the Collaborative Research Center 1294 ``Data assimilation'' is gratefully acknowledged.
  The research was supported by the Russian Science Foundation grant No. 18-11-00132.
} 

\def\rrbias{\rr_{b}}
\def\GPT{\cc{G}}
\def\fvfv{\fv \! \fv}
\def\epsVar{\mathcal{S}}
\def\SiGP{\epsVar | \GP}

\def\fblk{F}
\def\gblk{H}
\def\fgblk{A}
\def\dblk{D}
\def\IFtotal{F}
\def\IFT{\mathcal{\IFtotal}}
\def\IFTt{\tilde{\IFT}}

\def\DF{\mathcal{\dblk}}
\def\IF{\mathsf{\fblk}}
\def\IG{\mathsf{\gblk}}
\def\QF{\mathcal{\QP}}

\def\IFblk{\mathbb{\fblk}}
\def\DPblk{\mathbb{\dblk}}
\def\DPblkt{\tilde{\DPblk}}
\def\IGblk{\mathbb{\gblk}}
\def\HFblk{\mathbb{\gblk}}
\def\AFblk{\mathbb{\fgblk}}

\def\IFa{\breve{\IF}}
\def\IFTa{\breve{\IFT}}
\def\DPblka{\breve{\DPblk}}

\def\AFblka{\breve{\AFblk}}

\def\fvs{\fv^{*}}
\def\JG{\mu}
\def\GPM{\GP_{\JG}}
\def\GPUB{M}
\def\Fsmooth{\cc{F}}
\def\ProjG{\Pi}
\def\II{\mathcal{I}}

\def\upsilonvsa{\upsilonv^{\dag}}
\def\gvsa{\gv^{\dag}}

\begin{document}
\thispagestyle{empty}
{
\title{\thetitle}
\theauthors

\maketitle
\begin{abstract}
\theabstract
\end{abstract}

\ifAMS
    {\par\noindent\emph{AMS 2010 Subject Classification:} Primary \kwdp. Secondary \kwds}
    {\par\noindent\emph{JEL codes}: \kwdp}

\par\noindent\emph{Keywords}: \thekeywords
} 

\tableofcontents

\section{Introduction}
\label{Snoninverseintro}
Bayesian inference for inverse problems attracted a lot of attention in the recent literature.
We mention only few relevant papers. 
\cite{KnVaZa2011} studied minimax contraction rate for linear inverse problems, %
\cite{KnSzVa2016} discussed adaptive Bayes procedures. 
\cite{Nickl2017} studied the BvM for Schr{\"o}dinger equation, 
\cite{Nickl2019} focused on statistical inverse problems for compound Poisson processes, 
\cite{MoNiPa2017} discussed applications to X-Ray Tomography,
\cite{NiSo2017} studied posterior contraction rates for discretely observed scalar diffusions,
\cite{GuVaYa2018} considered Bayesian inverse problems with partial observations,
\cite{Trabs_2018} discussed a linear inverse problem with an unknown operator,
\cite{lu2017bernsteinvon} established BvM results for a rather general elliptic inverse problem
with an additive noise.
Nonlinearity of the model makes the study very involved and the cited results heavily used 
the recent advances in the theory of partial differential equations, inverse problems, empirical processes. 
We mention \cite{Nickl2017} and \cite{Nickl2018ConvergenceRF} as particular illustration 
of the major difficulties in the study of concentration of the penalized MLE and of posterior 
concentration. 

The main contribution of this paper is a novel approach allowing a unified study of a large class of 
nonlinear inverse problems.
The approach is based on a double relaxation by introducing an auxiliary functional parameter, 
replacing the structural equation with a penalty, and imposing an additional prior on the auxiliary parameter. 
This leads to a new model with an extended parameter set but the stochastic term is linear w.r.t. the total parameter set.  
%
This fact helps to obtain sharp finite sample bounds for concentration of the penalized 
maximum likelihood estimator (pMLE) around its population counterpart and for posterior concentration around pMLE. 
Also we establish a finite sample result about  
Gaussian approximation of the posterior with an explicit error term in the total variation distance and for the class of centrally symmetric sets around pMLE.
All the bounds are given in term of \emph{effective dimension} in place of the total parameter dimension.
This helps to compensate the increase of the parameter set and to get the right accuracy of approximation.
The approach is ``coordinate free'' and does not rely on any spectral decomposition and/or any basis representation for the target parameter and penalty term. 
We focus here on the problem of inverting an known nonlinear smooth operator from noisy discrete data
\( \Yv \) following the equation
\begin{EQA}
	\Yv 
	&=& 
	\KS(\fv) + \sigma \epsv .
\label{YAfseNI}
\end{EQA}
A forthcoming paper explains how the proposed approach called ``calming'' can be extended 
to many other models including generalized regression, nonparametric diffusion, Bayesian 
deconvolution, dimension reduction etc. 

Now we explain the idea of the method. 
For the original problem \eqref{YAfseNI}, a prior density \( \prior(\fv) \) on the target parameter \( \fv \) yields the posterior
\begin{EQA}
	\fv \cond \Yv
	& \propto &
	\exp\bigl\{ - \| \Yv - \KS(\fv) \|^{2}/(2\sigma^{2}) \bigr\} \prior(\fv) .
\label{fcYmYAf22s2}
\end{EQA}
Now denote by \( \gv \) the image function, \( \gv = \KS(\fv) \in \YSD \) and
relax the structural equation \( \gv = \KS(\fv) \)  
replacing it with a penalty \( \lambda \| \gv - \KS(\fv) \|^{2} \).
The image function \( \gv \) is modelled using a separate prior.
The proposed approach leads to the extended parameter set \( (\fv,\gv) \) which is modelled as
\begin{EQA}
	(\fv,\gv)
	& \propto &
	\exp\bigl\{ - \| \Yv - \gv \|^{2}/(2\sigma^{2}) - \lambda \| \gv - \KS(\fv) \| \bigr\} \prior(\fv) \prior(\gv) .
\label{fgpemYg22s2l}
\end{EQA}
Such a decoupling increases substantially the parameter space. 
However, by a proper choice of a \( \gv \)-prior one can keep the effective dimension of the same order as
for the original problem.
One can treat the calming approach 
as a kind of transformation of the original nonlinear problem to a linear one with an extended parameter set 
and a special prior that includes the structural penalty term. 
Our theoretical results justify the proposed method and state a number of remarkable features of the total 
and marginal posteriors. 

Certain important practical questions are not completely addressed in this paper. 
In particular, we do not discuss the question of estimation of the noise level \( \sigma \).
However, the approach continues to apply under noise misspecification.
The theoretical results only require some exponential moments of the errors, even inhomogeneous 
colored noise is allowed. 
Similarly, a proper choice of the prior should depend on the unknown regularity or smoothness of the source function \( \fvs \), and one needs a data-driven rule for practical applications;
see Section~\ref{SaltBayNI} for more discussion.
Note however, that the established Gaussian approximation of the posterior can be combined 
with existing methods of Bayesian model selection for Gaussian setup such as empirical of full Bayes methods; see e.g. 
\cite{KnSzVa2016,NiSz2016}, \cite{GiNi2019}.
We also focus on the case of Gaussian priors, an extension to more general regular priors is straightforward.

The paper is organized as follows.
Section~\ref{Snoninverse} explains the setup and the details of the proposed calming approach.
We also discuss the relations between alternative optimization and Bayesian MCMC-type procedures
for the obtained model.
Section~\ref{SBvMNI} presents the main results about the properties of the posterior distribution
including posterior concentration and Gaussian approximation of the posterior. 
Further results are given in Section~\ref{SCompNI}.
In particular, 
Section~\ref{SpMLENI} contains some important results on the properties of the penalized MLE. 
The bias induced by a prior on the source function \( \fv \) and the auxiliary prior on \( \gv \) 
is evaluated in Section~\ref{SsmoothNI}.
Bias-variance decomposition for the loss of the pMLE \( \tilde{\fv}_{\GP} \) 
and a bound on posterior contraction are presented in Section~\ref{Scontractcover}.
Section~\ref{SnonpmBayesNI} addresses the use of 
Bayesian credible sets as frequentist confidence sets. 
Section~\ref{SminmaxlinNI} discusses minimax optimality in the case of a general 
linear operator \( \KS \).
The results of this section seem to be of independent interest. 
Proofs of the main results are collected in Section~\ref{SCompNI}.
Technical statements and useful external results are put into the appendix.

\section{Non-linear Gaussian inverse problems}
\label{Snoninverse}
This section presents the proposed calming approach.
We begin with some motivation and discussion of the standard approach to
Bayesian inference in a non-linear inverse regression problems.
Consider the model
\begin{EQA}
	\Yv
	&=&
	\KS (\fv) + \sigma \epsv \in \YSD
\label{YvKSsepYS}
\end{EQA}
with a known non-linear operator \( \KS \colon \XSD \to \YSD \),
where \( \XSD,\YSD \) are discretized subspace of Hilbert spaces \( \XS, \YS \).
Examples include, e.g., the cases when \( \fv \) is a functional parameter of the Schr\"{o}dinger or 
Cald\'{e}ron equation; see \cite{Nickl2017,Nickl2018ConvergenceRF} for more details. 
Usually the subspace \( \XSD \) is spanned by 
a  given basis \( \bigl\{ \psiv_{j}\, , j=1,\ldots, \dimp \bigr\} \) in \( \XS \) and \( \dimp \leq \infty \).
This means that \( \fv \in \XSD \) can be expanded in the form
\begin{EQA}
	\fv
	&=&
	\sum_{j=1}^{\dimp} \theta_{j} \psiv_{j}
\label{f1ptjpsj}
\end{EQA}
and a vector of coefficients \( \thetav = (\theta_{j}) \in \R^{\dimp} \).
A discretized subspace \( \YSD \) is typically the vector of the solutions \( u_{\fv}(t_{i}) \) 
at given points \( t_{i} \), \( i=1,\ldots,n \).
This means that we observe the corresponding solution \( u_{\fv} \) at discrete points 
and corrupted with noise; \cite{Nickl2017,lu2017bernsteinvon}. 
In some cases, e.g. in the elliptic PDE, 
the conductivity coefficient \( \fv \) should be positive. 
Then a Gaussian prior is questionable, but everything applies to the log-transform
\( \fv \to \log \fv \), that is, the log-conductivity is modelled by a Gaussian prior.
The operator \( \KS \) has to be changed correspondingly.
It is important that such a transform does not affect smoothness properties of the operator
\( \KS \).

First we discuss a standard approach.
A Gaussian log-likelihood in the model \eqref{YvKSsepYS} w.r.t. \( \fv \) is of the form
\( - (2\sigma^{2})^{-1} \bigl\| \Yv - \KS(\fv) \bigr\|^{2} \) up to a constant term.
%
A Gaussian prior \( \fv \sim \ND(\fv_{0},\GP^{-2}) \) on \( \fv \in \XSD \) yields the penalized log-likelihood
\begin{EQA}
	L_{\GP}(\fv)
	&=&
	- (2\sigma^{2})^{-1} \bigl\| \Yv - \KS(\fv) \bigr\|^{2}
	- \frac{1}{2} \bigl\| \GP (\fv - \fv_{0}) \bigr\|^{2} \, 	.
\label{LLt2s2YkKPtt122}
\end{EQA}
The posterior is obtained by normalizing the product \( \exp L_{\GP}(\fv) \).
The nonlinear mapping \( \KS(\fv) \) in the data fidelity term creates fundamental problems 
for studying the posterior behavior.
Existing results utilise deep tools from empirical processes, PDE, inverse problem;
see e.g. \cite{Nickl2018ConvergenceRF}.
A first naive idea to avoid these technical difficulties is to use a reparametrization. 
Assume that \( \KS \) is injective and denote by \( \KSI \) its inverse:
\( \KSI \eqdef \KS^{-1} \colon \YS \to \XS \).
Then one can use another parametrization \( \gv = \KS(\fv) \in \YSD \).
This yields the classical quadratic log-likelihood
\begin{EQA}
	L(\gv)
	&=&
	- \frac{1}{2 \sigma^{2}} \bigl\| \Yv - \gv \bigr\|^{2} .
\label{Lt12s2YPt2NL}
\end{EQA}
Let now \( \fv \) be a Gaussian element in \( \XSD \) with a mean \( \fv_{0} \) and a self-adjoint covariance operator 
\( \GP^{-2} \colon \XS \to \XS \). 
It leads to a non-Gaussian prior \( \KSI(\gv) \) on \( \YS \).
The penalized log-likelihood \( L_{\GP}(\gv) \) reads
\begin{EQA}
	L_{\GP}(\gv)
	&=&
	- \frac{1}{2 \sigma^{2}} \bigl\| \Yv - \gv \bigr\|^{2} 
	- \frac{1}{2} \bigl\| \GP \bigl( \KSI(\gv) - \fv_{0} \bigr) \bigr\|^{2}
	+ \CONST.
\label{Lt12GKPt2}
\end{EQA}
Quadratic structure of the fidelity term 
\( - (2 \sigma^{2})^{-1}\bigl\| \Yv - \gv \bigr\|^{2} \)
makes the stochastic analysis of the posterior much more simple.
Further, it holds for the expected log-likelihood
\begin{EQA}
	\E L_{\GP}(\gv)
	&=&
	- \frac{1}{2 \sigma^{2}} \bigl\| \gvs - \gv \bigr\|^{2} 
	- \frac{1}{2} \bigl\| \GP \bigl( \KSI(\gv) - \fv_{0} \bigr) \bigr\|^{2}
	+ \CONST
\label{ELGC1212s2gPt}
\end{EQA}
with \( \gvs = \KS(\fvs) \).
Contraction properties of the posterior can be effectively studied if the expected penalized log-likelihood is concave and smooth.
In particular, we need that the inverse operator \( \KSI \) is twice differentiable, 
a sufficient condition for concavity reads \( - \nabla^{2} \E L_{\GP}(\gv) \geq 0 \).
Note, however, that the inverse operator \( \KSI \) and its second derivative are usually unbounded
and non-smooth. 
This makes the theoretical study extremely challenging. 


\subsection{Decoupling and calming}
Now we present a slightly different approach which allows to avoid any use of the inverse operator \( \KSI \).
The basic idea is to introduce an auxiliary parameter \( \gv \) which means a ``smooth'' approximation of the data,
and replace the structural equation \( \gv = \KS(\fv) \) by a penalty term. 
More precisely, represent \eqref{YvKSsepYS} by two identities
\( \Yv = \gv + \sigma \epsv \) and \( \gv = \KS(\fv) \).
Then the fidelity terms can be relaxed to
\begin{EQA}[c]
	\sigma^{-2} \bigl\| \Yv - \gv \bigr\|^{2} + \lambda \| \gv - \KS(\fv) \|^{2}
\label{sm2Yg2lkKf22}
\end{EQA}
with a Lagrange multiplier \( \lambda \).
Now we proceed with a couple of parameters \( (\fv, \gv) \).
Bayesian modeling assumes regular priors on both parameters \( \gv \) and \( \fv \).
In particular, one can use independent Gaussian priors \( \fv \sim \ND(\fv_{0},\GP^{-2}) \) and 
\( \gv \sim \ND(\gv_{0},\GPY^{-2}) \) yielding the posterior
\begin{EQA}
	(\fv,\gv) \cond \Yv
	& \propto &
	\exp \LLG(\fv,\gv)
\label{thYeLGthv}
\end{EQA}
with
\begin{EQA}
	\LLG(\fv,\gv)
	&=& 
	- \frac{1}{2\sigma^{2}} \bigl\| \Yv - \gv \bigr\|^{2} 
		- \frac{\lambda}{2} \| \gv - \KS(\fv) \|^{2} 
	- \frac{1}{2} \bigl\| \GP (\fv - \fv_{0}) \bigr\|^{2}
	- \frac{1}{2} \bigl\| \GPY (\gv - \gv_{0}) \bigr\|^{2} \, .
	\qquad
\label{Lth122Y1222}
\end{EQA}
This expression is quadratic in \( \gv \) for a fixed \( \fv \).
A dependence on \( \fv \) is a bit more complicated due to the structural term
\( \lambda \| \gv - \KS(\fv) \|^{2} \) which penalizes for deviations from the forward non-linear structural relation 
\( \gv = \KS (\fv) \).
It is important that this expression does not involve any inversion of the operator \( \KS \).
The original stochastic data only enters in the quadratic term \( \| \Yv - \gv \|^{2} \), this incredibly simplifies the stochastic analysis. 
The gradient and Hessian of \( \LLG(\fv,\gv) \) read as follow:
\begin{EQA}
	\frac{d}{d\fv} \LLG(\fv,\gv)
	&=& 
	- \lambda \, \bigl\{ \KS(\fv) - \gv \bigr\}^{\T} \nabla \KS(\fv) 
	- \GP^{2} (\fv - \fv_{0}) ,
	\\
	\frac{d}{d\gv} \LLG(\fv,\gv)
	&=& 
	\sigma^{-2} \bigl( \Yv - \gv \bigr) 
	+ \lambda \, \bigl\{ \KS(\fv) - \gv \bigr\}
	 - \GPY^{2} \bigl( \gv - \gv_{0} \bigr) ,
\label{Lth122Y1ppp}
\end{EQA}
and
\begin{EQA}
	\IFT_{\GP}(\fv,\gv)
	&\eqdef&
	- \nabla^{2} \E \LLG(\fv,\gv)
	\\
	&=&
	\begin{pmatrix}
	\GP^{2} + \lambda \bigl\{ \nabla \KS(\fv)^{\T} \nabla  \KS(\fv) + \deltavd^{\T} \nabla^{2} \KS(\fv) \bigr\} &
	\lambda \, \nabla \KS(\fv) 
	\\
	\lambda \, \bigl\{ \nabla \KS(\fv) \bigr\}^{\T} &
	\GPY^{2} + (\sigma^{-2} + \lambda) \Id_{\dimq}
	\end{pmatrix}
\label{G2Kzsm2PYG}
\end{EQA}
with the elasticity vector \( \deltavd = \KS(\fv) - \gv \).
The only term in the matrix 
\( \IFT_{\GP}(\fv,\gv) \) involving the second derivative of 
\( \KS(\cdot) \) appears in the \( \fvfv \)-block with the multiplicative factor \( \deltavd \) which is 
small if the structural equation \( \KS(\fv) = \gv \) is nearly fulfilled.
The proof of posterior concentration uses concavity of \( \E \LLG(\fv,\gv) \).
For this property, it suffices to check that \( \IFT_{\GP}(\fv,\gv) \geq 0 \).
A global check of this condition is difficult due to the nonlinear term \( \nabla^{2} \KS(\fv) \).
However, it is automatically fulfilled for all couples \( (\fv,\gv) \) with 
\( \deltavd = \KS(\fv) - \gv = 0 \), and hence, in a local vicinity of such points. 

In a special case when \( \KS \) is linear and \( \nabla \KS(\fv) = \KS \), it holds 
\( \nabla^{2} \KS = 0 \) and the matrix \( \IFT_{\GP}(\fv,\gv) \)
does not depend on \( (\fv,\gv) \):
\begin{EQA}
	\IFT_{\GP}
	=
	- \nabla^{2} \E \LLG(\fv,\gv)
	&=&
	\begin{pmatrix}
	\GP^{2} + \lambda \, \KS^{\T} \KS 
	& 
	\lambda \, \KS  
	\\
	\lambda \, \KS^{\T}  &
	\GPY^{2} + (\sigma^{-2} + \lambda) \Id_{\dimq} 
	\end{pmatrix} 
	=
	\begin{pmatrix}
	\IF_{\GP} &
	\lambda \, \KS 
	\\
	\lambda \, \KS^{\T} &
	\IG_{\GP}
	\end{pmatrix}\, .
\label{G2Kzsm2PYGlin}
\end{EQA}
With \( \GP^{2} > 0 \) and \( \GPY^{2} > 0 \), this matrix is also positive 
semidefinite. 
The \( \fvfv \)-block of the matrix \( \IFT_{\GP}^{-1} \) reads
\begin{EQA}
	(\IFT_{\GP}^{-1})_{\fvfv}
	&=&
	\bigl( \IF_{\GP} - \lambda^{2} \, \KS^{\T} \IG_{\GP}^{-1} \KS \bigr)^{-1}
	\\
	&=&
	\Bigl( \GP^{2} + \lambda \, \KS^{\T} \KS 
	- \lambda^{2} \KS^{\T} \bigl\{ \GPY^{2} + (\sigma^{-2} + \lambda) \Id_{\dimq} \bigr\}^{-1} \KS
	\Bigr)^{-1} , 
\label{m1IIm1ttet}
	\\
	(\IFT_{\GP}^{-1})_{\gv\gv}
	&=&
	\IG_{\GP}^{-1} + \lambda^{2} \IG_{\GP}^{-1} \KS^{\T}
	\bigl( \IF_{\GP} - \lambda^{2} \, \KS^{\T} \IG_{\GP}^{-1} \KS \bigr)^{-1} \KS \IG_{\GP}^{-1}
	\, .
\label{m1IIm1ttett}
\end{EQA}
With \( \GPY^{2} = 0 \),
\begin{EQA}
	(\IFT_{\GP}^{-1})_{\fvfv}
	&=&
	\Bigl( 
		\GP^{2} + \frac{ \sigma^{-2} \lambda}{\sigma^{-2} + \lambda} \, \KS^{\T} \KS 
	\Bigr)^{-1}
	\, .
\label{m1IIm1ttet2}
\end{EQA}
This implies the following simple result.

\begin{lemma}
\label{LeffdimNLI}
Let \( \KS(\fv) = \KS \fv \) for a linear mapping \( \KS \colon \XSD \to \YSD \).
It holds 
\begin{EQ}[rcccl]
	\bigl( \GP^{2} + \lambda \, \KS^{\T} \KS \bigr)^{-1}
	& \leq &
	\bigl( \IFT_{\GP}^{-1} \bigr)_{\fvfv}
	& \leq &
	\Bigl( 
		\GP^{2} + \frac{ \sigma^{-2} \lambda}{\sigma^{-2} + \lambda} \, \KS^{\T} \KS 
	\Bigr)^{-1} \, ,
	\\
	\bigl( \GPY^{2} + (\sigma^{-2} + \lambda) \Id_{\dimq} \bigr)^{-1}
	& \leq &
	\bigl( \IFT_{\GP}^{-1} \bigr)_{\gv\gv}
	& \leq &
	\Bigl( 
		\GPY^{2} + \sigma^{-2} \Id_{\dimq} 
	\Bigr)^{-1} \, .
\label{FGm12blDGm2HGm2}
\end{EQ}
In particular, for \( \lambda = \sigma^{-2} \), it holds
\begin{EQA}
	\frac{1}{2} \block\bigl\{ \IF_{\GP}^{-1}, \IG_{\GP}^{-1} \bigr\}
	& \leq &
	\IFT_{\GP}^{-1}
	\leq 
	2 \block\bigl\{ \IF_{\GP}^{-1}, \IG_{\GP}^{-1} \bigr\} .
\label{12blkIGm1HGIGm1}
\end{EQA}
\end{lemma}

\begin{proof}
Substituting \( \GPY^{2} = 0 \) and \( \GPY^{2} = \infty \) in \eqref{m1IIm1ttet} yields the first line of 
\eqref{FGm12blDGm2HGm2}.
The second one is obtained in the same way by manipulating with \( \GP^{2} \).
\end{proof}

One can see that 
each of  regularizations by \( \GP^{2} \) and \( \GPY^{2} \) improves 
the conditioning number of 
\( \IFT_{\GP} \).
Moreover, with \( \lambda = \sigma^{-2} \), the \( \fvfv \)-block of 
\( \IFT_{\GP}^{-1} \) behaves as \( \IF_{\GP}^{-1} \), 
that is, the estimation of \( \fv \) in the extended \( (\fv,\gv) \)-model 
yields the same rate result as in the standard estimation in the original 
\( \fv \)-model.

\subsection{Alternating optimization versus iterated Bayes }
\label{SaltBayNI}
This section discusses relations between alternating optimization for computing the pMLE
\( \tilde{\fv}_{\GP} \) and the Bayesian calculus with calming.
First we explain how the penalized MLE 
\begin{EQA}
	\bigl( \tilde{\fv}_{\GP}, \tilde{\gv}_{\GP} \bigr)
	& \eqdef &
	\argmax_{(\fv,\gv)} \LLG\bigl( \fv,\gv \bigr)
	\\
	&=& 
	\argmin_{(\fv,\gv)} \frac{1}{\sigma^{2}} \bigl\| \Yv - \gv \bigr\|^{2} 
		+ {\lambda} \| \gv - \KS(\fv) \|^{2} 
	+ \bigl\| \GP (\fv - \fv_{0}) \bigr\|^{2}
	+ \bigl\| \GPY (\gv - \gv_{0}) \bigr\|^{2} 
	\qquad
\label{tfvGtgvGLLGfv}
\end{EQA}
can be computed by an alternating procedure.
Note that the partial penalized MLE in \( \gv \) given \( \fv \) can be computed explicitly.
Indeed, the equation \( \frac{d}{d\gv} \LLG(\fv,\gv) = 0 \) yields
\begin{EQA}
	\tilde{\gv}_{0}(\fv_{0})
	& \eqdef &
	\argmax_{\gv} \LLG\bigl( \fv,\gv \bigr)
	=
	\bigl\{ (\sigma^{-2} + \lambda) \Id_{\dimq} + \GPY^{2} \bigr\}^{-1} 
	\bigl\{ \sigma^{-2} \Yv + \lambda \KS(\fv) + \GPY^{2} \gv_{0} \bigr\} .
	\qquad
\label{tthsm2lsPY}
\end{EQA}
Further, 
with \( \tilde{\gv}_{0} \) fixed, optimization w.r.t. \( \fv \) yields the equation
\begin{EQA}
	\lambda \, 
	\bigl\{ \nabla \KS(\fv) \bigr\}^{\T} \bigl\{ \KS(\fv) - \tilde{\gv}_{0} \bigr\}
	&=& 
	- \GP^{2} (\fv - \fv_{0}) .
\label{lanAPtOPt0G2}
\end{EQA}
Now we assume that \( \fv_{0} \) is chosen properly and the operator \( \KS(\fv) \) can be well approximated by its second order expansion around \( \fv_{0} \):
\( \KS(\fv) \approx \KS(\fv_{0}) + \nabla \KS(\fv_{0}) \, (\fv - \fv_{0}) \)
and 
\( \nabla \KS(\fv) \approx 
\nabla \KS(\fv_{0}) + \nabla^{2} \KS(\fv_{0}) \, (\fv - \fv_{0}) \).
This yields with \( \deltavd_{1} = \KS(\fv_{0}) - \tilde{\gv}_{0} \) an update 
\begin{EQA}
	\tilde{\fv}_{1}
	&=&
	\fv_{0} +
	\lambda \, \bigl[ 
		\GP^{2} + \lambda \bigl\{ \nabla \KS(\fv_{0})^{\T} \nabla \KS(\fv_{0}) + \deltavd_{1}^{\T} \nabla^{2} \KS(\fv_{0}) \bigr\} 
	\bigr]^{-1} 
	\\
	&& \qquad
	\times
	\nabla \KS(\fv_{0})^{\T} \bigl\{ \KS(\fv_{0}) - \tilde{\gv}_{0} \bigr\}.
	\qquad
\label{ttdG2laPTbAPff0}
\end{EQA}
%
Moreover, if the guess \( \fv_{0} \) does a good job, then the related elasticity vector \( \deltavd_{1} = \KS(\fv_{0}) - \tilde{\gv}_{0} \)
is small and the corresponding term in \eqref{ttdG2laPTbAPff0} can be ignored leading to the update
\begin{EQA}
	\tilde{\fv}_{1}
	&=&
	\fv_{0} +
	\lambda \, \bigl\{ 
		\GP^{2} + \lambda \nabla \KS(\fv_{0})^{\T} \nabla \KS(\fv_{0}) 
	\bigr\}^{-1} \nabla \KS(\fv_{0})^{\T} \bigl\{ \KS(\fv_{0}) - \tilde{\gv}_{0} \bigr\}.
	\qquad
\label{ttdG2laPTbAPff1}
\end{EQA}
The whole procedure starts from an initial guess \( \fv_{0} \) and alternates the steps 
\eqref{tthsm2lsPY} and \eqref{ttdG2laPTbAPff0} or \eqref{ttdG2laPTbAPff1}.
The technique from \cite{AASP2015} can be used to guarantee the exponential convergence of this alternating procedure
to the solution \( (\tilde{\fv}_{\GP}, \tilde{\gv}_{\GP}) \). 
However, even after many simplifications, the alternating procedure \eqref{tthsm2lsPY}--\eqref{ttdG2laPTbAPff1} 
involves computing of a high-dimensional gradient 
and inverting a large matrix.
This especially concerns the step \eqref{ttdG2laPTbAPff0} or \eqref{ttdG2laPTbAPff1}.
The described procedure is somehow related to the general 
Alternating Direction Method of Multipliers (ADMM) methods of large scale optimization;
see e.g. \cite{vono2019} and references therein.

Now we discuss the sequential Bayes procedure which mimics the same type of alternating
\eqref{tthsm2lsPY}--\eqref{ttdG2laPTbAPff1}. 
For this procedure, we only need to efficiently compute the forward operator \( \KS(\fv) \).
Again, we start with a guess \( \fv_{0} \) and compute \( \gv_{0} = \KS(\fv_{0}) \).
Now we generate from the prior \( \ND(\fv_{0},\GP^{-2}) \) and \( \ND(\gv_{0},\GPY^{-2}) \) 
and compute the corresponding posterior using the penalized log-likelihood \( \LLG(\fv,\gv) \).
Let also
\( \tilde{\fv}_{1} \) be the MAP or posterior mean after our Bayes simulations.
This step replaces \eqref{ttdG2laPTbAPff1}.
For the next step we use the prior \( \ND(\tilde{\fv}_{1},\GP^{-2}) \) and 
\( \ND(\tilde{\gv}_{1},\GPY^{-2}) \) with \( \tilde{\gv}_{1} = \KS(\tilde{\fv}_{1}) \).
We formally need a new training set \( \Yv^{(1)} \).
This allows to forget about the random nature and data-dependence of the posterior distribution 
when we use it as a prior. 
Now we repeat the Bayes step with 
\begin{EQA}
	\LLG^{(1)}(\fv,\gv)
	&=& 
	- \frac{1}{2\sigma^{2}} \bigl\| \Yv^{(1)} - \gv \bigr\|^{2} 
		- \frac{\lambda}{2} \| \gv - \KS(\fv) \|^{2} 
	- \frac{1}{2} \bigl\langle \fv - \tilde{\fv}_{1} \bigr\rangle_{\GP}
	- \frac{1}{2} \bigl\langle \gv - \tilde{\gv}_{1} \bigr\rangle_{\GPY} \, .
	\qquad
\label{Lth122Y12221}
\end{EQA}
To accelerate Bayes computations and to improve the numerical performance, one can also update the prior covariance
at each step using the posterior covariance from the previous step. 
For the theoretical study, it is important that a new sample \( \Yv^{(k)} \) is used for each step \( k \).
It is worth mentioning that the Bayes procedure is gradient free and does not involve computing or inverting any matrix;
cf. \cite{NeSp2016}. 
For a given sample \( \fv^{(m)} \) from the prior \( \ND(\fv_{0},\GP^{-2}) \) we have to compute the forward operator \( \KS(\fv^{(m)}) \), this step computationally most expensive. 
Afterward, computing the penalized log-likelihood is straightforward.
The procedure is scalable in dimension and sample size. 
Our theoretical result claim that the posterior mimics the second order optimization routine and
hence, ensures a convergence after logarithmic number of steps. 

A practical implementation of the method requires to specify the priors \( \ND(\fv_{0},\GP^{-2}) \)
for the target parameter \( \fv \) and \( \ND(\gv_{0},\GPY^{-2}) \) 
for the auxiliary parameter \( \gv \).
First we comment on the choice of \( \gv \)-prior. 
It has to fulfill the coordination condition \nameref{GGref} below that recommends to apply 
the Gaussian measure on \( \YSD \) obtained by a linear mapping \( \nabla \KS(\fv_{0}) \)
from the Gaussian prior for \( \fv \).
This yields
\( \gv_{0} = \KS(\fv_{0}) \) and 
\( \GPY^{2} = \nabla \KS(\fv_{0}) \,  \GP^{2} \, \nabla \KS(\fv_{0})^{\T} \).
A proper choice of the precision operator \( \GP^{-2} \) for the \( \fv \)-prior depends on the regularity 
of the target function \( \fv \).
A data-driven choice can be performed by the empirical or full Bayes approach.
We refer to \cite{KnSzVa2016,NiSz2016,SnVa2015,belitser2017} for description and 
justification of this approach for linear inverse problems
or \cite{GiNi2019} for the nonlinear elliptic case.
Results of this paper on Gaussian approximation of the posterior justify the applicability of 
these methods in the general nonlinear case.

\section{Bernstein -- von Mises Theorem}
\label{SBvMNI}
This section discusses general properties of the posterior for the proposed method. 
The calming approach extends the parameter space: 
the target source function \( \fv \) 
is accomplished by the nuisance parameter \( \gv \) from the image space.
This auxiliary parameter can be even infinite dimensional.
The proposed procedure defines the estimate \( (\tilde{\fv}_{\GP},\tilde{\gv}_{\GP}) \) of the couple 
\( (\fvs,\gvs) \) by maximizing \( \LLG(\fv,\gv) \) from \eqref{Lth122Y1222}.
The question under study is whether 
\( \tilde{\fv}_{\GP} \) is a proper estimator of the target function \( \fvs \).
The related question is whether the posterior on \( \fv \) obtained as the marginal of the joint posterior 
of \( (\fv,\gv) \) possesses some contraction properties and can be approximated by a Gaussian measure.
We show below that a proper choice of the prior on \( \gv \) allows to keep the total effective dimension
essentially the same as in the case of the original problem before calming.
One can say that the calming approach helps to avoid most of very complicated empirical process study 
at no additional costs. 
Our results are nearly sharp, stated for finite samples under mild and reasonable conditions.

\subsection{Conditions} 
\label{Scondgllo} 
This section lists the conditions which appear to be sufficient for stating the concentration and BvM results.
%
Consider linear discretized subspaces \( \XSD \subset \XS \) 
and similarly \( \YSD \subset \YS \).
We will assume that the \( \fv \)-priors \( \ND(\fv_{0},\GP^{-2}) \) is concentrated on 
\( \XSD \) while the \( \gv \)-prior \( \ND(\gv_{0},\GPY^{-2}) \) on \( \YSD \).
By \( \dimq \) we denote the dimension of \( \YSD \) while \( \dimp \) stands for 
the dimension of \( \XSD \) and \( \dimtotal \) for the total dimension of \( \XSD \times \YSD \).
For notational simplicity we assume \( \dimtotal < \infty \), however,
the study extends to the case of \( \dimtotal = \infty \) in a straightforward way.
It is convenient to introduce the full parameter 
\( \upsilonv = (\fv,\gv) \in \UpsilonD = \XSD \times \YSD \)  
and denote \( \LLG(\upsilonv) \eqdef \LLG(\fv,\gv) \)
\begin{EQA}
	\LLG(\upsilonv)
	&=& 
	- \frac{1}{2\sigma^{2}} \bigl\| \Yv - \gv \bigr\|^{2} 
		- \frac{\lambda}{2} \| \gv - \KS(\fv) \|^{2} 
	- \frac{1}{2} \bigl\| \GP (\fv - \fv_{0}) \bigr\|^{2}
	- \frac{1}{2} \bigl\| \GPY (\gv - \gv_{0}) \bigr\|^{2} \, .
	\qquad
\label{LLth122Y1}
\end{EQA}
It also holds with \( \gvs = \KS(\fvs) \)
\begin{EQA}
	\E \LLG(\upsilonv)
	&=& 
	- \frac{1}{2\sigma^{2}} \bigl\| \gvs - \gv \bigr\|^{2} 
	- \frac{\lambda}{2} \| \gv - \KS(\fv) \|^{2} 
	- \frac{1}{2} \bigl\| \GP (\fv - \fv_{0}) \bigr\|^{2}
	- \frac{1}{2} \bigl\| \GPY (\gv - \gv_{0}) \bigr\|^{2} 
	\qquad
\label{ELLth122Y1}
\end{EQA}
up to a constant term which we omit. 
Given a local set \( \Upsilond \), define
\begin{EQA}[rcccl]
	\tilde{\upsilonv}_{\GP}
	& \eqdef &
	\argmax_{\upsilonv \in \Upsilond} \LLG(\upsilonv) ,
\label{tsGetsGthLL}
	\qquad
\label{tsGetsGthELL}
	\upsilonvs_{\GP}
	& \eqdef &
	\argmax_{\upsilonv \in \Upsilond} \E \LLG(\upsilonv) .
\end{EQA}
Below we heavily use that the data \( \Yv \) only 
enter in the fidelity term \( \sigma^{-2} \| \Yv - \gv \|^{2}/2 \)
and therefore, the stochastic term linearly depends on \( \gv \) and free of
\( \fv \).
%
%
Let also \( \IFT_{\GP}(\upsilonv) = - \nabla^{2} \E \LLG(\upsilonv) \) be the 
negative Hessian of \( \E \LLG(\upsilonv) \).
It can be written in the block form:
\begin{EQA}
	\IFT_{\GP}(\upsilonv)
	& \eqdef &
	- \nabla^{2} \E \LLG(\fv,\gv)
	\\
	&=&
	\begin{pmatrix}
	\GP^{2} + \lambda \bigl\{ \nabla \KS(\fv)^{\T} \nabla \KS(\fv) + \deltavd^{\T} \nabla^{2} \KS(\fv) \bigr\} &
	\lambda \, \nabla \KS(\fv) 
	\\
	\lambda \, \bigl\{ \nabla \KS(\fv) \bigr\}^{\T} &
	\GPY^{2} + (\sigma^{-2} + \lambda) \Id_{\dimq}  
	\end{pmatrix} 
\label{G2Kzsm2PYGfg}
\end{EQA}
with \( \deltavd = \deltavd(\upsilonv) = \gv - \KS(\fv) \).
It is worth mentioning that the \( \gv\gv \)-block of \( \IFT_{\GP}(\upsilonv) \) 
is constant and does not depend on \( \upsilonv \). 
We use the notations
\begin{EQA}
\label{c}
	\IF_{\GP}(\upsilonv)
	& \eqdef &
	\GP^{2} + \lambda \bigl\{ 
		\nabla \KS(\fv)^{\T} \nabla \KS(\fv) + \deltavd^{\T} \nabla^{2} \KS(\fv) 
	\bigr\},
	\\
	\IG_{\GP}
	& \eqdef &
	\GPY^{2} + (\sigma^{-2} + \lambda) \Id_{\dimq} \, 
\label{DPGP2HPGP2GG2}
\end{EQA}
for the diagonal blocks of \( \IFT_{\GP}(\upsilonv) \).
%
Note also that the the auxiliary parameter \( \gv \) enters in \( \IF_{\GP}(\upsilonv) \) 
only through the elasticity vector \( \deltavd \). 
For the points \( \upsilonv \) with \( \deltavd = \deltavd(\upsilonv) = 0 \), 
the term with the second derivative of \( \KS(\fv) \) vanishes leading to the matrices 
\begin{EQA}
	\IFa_{\GP}(\fv)
	& \eqdef &
	\GP^{2} + \lambda \nabla \KS(\fv)^{\T} \nabla \KS(\fv) , 
	\\
	\IFTa_{\GP}(\fv)
	& \eqdef &
	\begin{pmatrix}
	\IFa_{\GP}(\fv) &
	\lambda \, \nabla \KS(\fv) 
	\\
	\lambda \, \bigl\{ \nabla \KS(\fv) \bigr\}^{\T} &
	\IG_{\GP}  
	\end{pmatrix} \, .
\label{DPGP2HPGP2GG2m2}
\end{EQA}
Our results make systematic use of the matrices \( \IFT(\upsilonv) \) 
for \( \upsilonv = \upsilonvs_{\GP} \) and
\( \upsilonv = \tilde{\upsilonv}_{\GP} \). 
We denote \( \IFT_{\GP} = \IFT_{\GP}(\upsilonvs_{\GP}) \),
\( \DF_{\GP} = \sqrt{\IFT_{\GP}} \), and
\( \DFt_{\GP} = \sqrt{\IFT_{\GP}(\tilde{\upsilonv})} \).
We also show that one can use \( \IFTa_{\GP}(\fv) \) and its blocks in place of \( \IFT_{\GP}(\upsilonv) \) for those \( \upsilonv \).
Such a replacement corresponds to a local approximation of the original nonlinear model by a linear one.
A bound on the inverse of \( \IFTa_{\GP}(\fv) \) in terms of the blocks 
\( \IFa_{\GP}(\fv) \) and \( \IG_{\GP} \) is very useful; see Lemma~\ref{LeffdimNLI} below.  

Our conditions rely on the inverse of \( \IFT_{\GP}(\upsilonv) \) and, in particular, 
on its \( \fvfv \) and \( \gv\gv \) blocks.
The use of the block inversion yields
\begin{EQA}
	\IFT_{\GP}^{-1}(\upsilonv)
	&=&
	\begin{pmatrix}
		\IFblk_{\GP}^{-1}(\upsilonv) & \AFblk_{\GP}(\upsilonv)
		\\
		\AFblk_{\GP}^{\T}(\upsilonv) & \IGblk_{\GP}^{-1}(\upsilonv)
	\end{pmatrix}	
\label{IFm1uDAATHm2}
\end{EQA}
with
\begin{EQ}[rclcl]
	\IFblk_{\GP}^{-1}(\upsilonv)
	& \eqdef &
	\bigl\{ \IFT_{\GP}^{-1}(\upsilonv) \bigr\}_{\fvfv}
	&=&
	\Bigl\{ \IF_{\GP}(\upsilonv) - \lambda^{2} \, \nabla \KS^{\T}(\fv) \, \IG_{\GP}^{-1} \, \nabla \KS(\fv) \Bigr\}^{-1} ,
	\\
	\IGblk_{\GP}^{-1}(\upsilonv)
	& \eqdef &
	\bigl\{ \IFT_{\GP}^{-1}(\upsilonv) \bigr\}_{\gv\gv}
	&=&
	\IG_{\GP}^{-1} + \lambda^{2} \, 
	\IG_{\GP}^{-1} \nabla \KS^{\T}(\fv) \, \IFblk_{\GP}^{-1}(\upsilonv) \, \nabla \KS(\fv) \, \IG_{\GP}^{-1},
\label{HG2IFGueeIII}
\end{EQ}
and
\begin{EQA}
	\AFblk_{\GP}(\upsilonv)
	& \eqdef &
	\bigl\{ \IFT_{\GP}^{-1}(\upsilonv) \bigr\}_{\fv\gv}
	= 
	- \lambda \, \IFblk_{\GP}^{-1}(\upsilonv) \,
	\nabla \KS^{\T}(\fv) \, \IG_{\GP}^{-1} .
\label{HG2IFGueeIIIA}
\end{EQA}
By \( \IFT(\upsilonv) \) we denote a similar matrix corresponding to the non-penalized
log-likelihood \( \LL(\upsilonv) \). 
It can be formally obtained by letting 
\( \GP = \GPY = 0 \).
Also denote \( \IFT = \IFT(\upsilonvs_{\GP}) \) and define the corresponding square root \( \DF = \sqrt{\IFT} \).
Now we are about to state our conditions.

\begin{description}
\item[\label{LLref} \( \bb{(\LL)} \)]
  \textit{\( \Upsilond \) is an open and convex subset in 
  \( \UpsilonD = \XSD \times \YSD \).
  The function \( \E \LLG(\upsilonv) \) is concave in 
  \( \upsilonv \in \Upsilond \).
  }
\end{description}

A condition of a global concavity of \( \E \LLG(\upsilonv) \) on \( \UpsilonD \) can be restrictive and difficult to check because of the squared norm of the nonlinear term \( \gv - \KS(\fv) \).
That is why we state the condition for a local set \( \Upsilond \) and restrict 
the prior to this set. 
Inspection of the proof reveals that this condition can be relaxed to the condition that the function
\( - \E \LLG(\upsilonv) \) can be bounded from below by a quadratic function inside of a certain 
elliptic set around \( \upsilonvs_{\GP} \) and linear outside of this ball. 

The stochastic component \( \zeta(\upsilonv) = \LLG(\upsilonv) - \E \LLG(\upsilonv) \) of the process 
\( \LLG(\upsilonv) \) is linear in \( \gv \) and free of \( \fv \) by construction.
More precisely, \eqref{LLth122Y1} implies \( \nabla \zeta = (\nabla_{\fv} \zeta,\nabla_{\gv} \zeta) \)
with 
\begin{EQA}
	\nabla_{\fv} \zeta
	& \eqdef &
	\frac{d \zeta}{d\fv}(\upsilonv)
	\equiv
	0,
	\\
	\nabla_{\gv} \zeta
	& \eqdef &
	\frac{d \zeta}{d\gv}(\upsilonv)
	\equiv
	\sigma^{-2} (\Yv - \E \Yv) 
	=
	\sigma^{-1} \epsv
\label{nanafvgvsm2e}
\end{EQA}
for \( \epsv \eqdef \sigma^{-1} (\Yv - \E \Yv) \).
Our likelihood function is built for a standard Gaussian homogeneous noise 
\( \epsv = \sigma^{-1} (\Yv - \E \Yv) \) in \eqref{YAfseNI}.
However, our results apply for a non-Gaussian inhomogeneous possible correlated noise \( \epsv \), 
and we impose rather mild assumption on its distribution. 
Namely, we only require some exponential moment of \( \epsv \).

\begin{description}
\item[\( \bb{(E\epsVar)}\)\label{ED0ref}]
  \emph{There exist a positive \( \dimq \times \dimq \) symmetric matrix \( \epsVar \),
    and constants \( \gmb > 0 \), \( \nunu \ge 1 \) such that
    \( \Var(\epsv) \le \epsVar^{2} \) and 
  }
\begin{EQA}[c]
\label{expzetaclocGP}
    \sup_{\uv \in \YSD} 
    \log \E \exp\biggl\{
      	\lambda \frac{\langle \uv, \epsv \rangle}{\| \epsVar \uv \|}
    \biggr\} \le
    \frac{\nunu^{2} \lambda^{2}}{2}, \qquad 
    |\lambda| \le \gmb.
\end{EQA}
\end{description}
Under correct noise specification \( \epsv \sim \ND(0, \Id_{\dimp}) \), one can take 
\( \epsVar = \Id_{\dimp} \).
Condition \nameref{ED0ref} with \( \gmb = \infty \) means sub-Gaussian errors.
In fact, we only need a deviation bound for the quadratic form
\( \epsv^{\T} \BB \epsv \) for a specific matrix \( \BB \); see \eqref{uTzDGvTxm12z} below.
Condition \nameref{ED0ref} is sufficient but not necessary. 
%

A proper choice of the priors on \( \fv \) and, especially, on \( \gv \) 
in the calming approach is crucial for our results. 
Namely, we assume that the prior \( \ND(\gv_{0},\GPY^{-2}) \) 
on the image \( \gv = \KS(\fv) \) is well coordinated with the prior \( \ND(\fv_{0},\GP^{-2}) \) 
on source function \( \fv \).

\begin{description}
\item[\( \bb{(\GP|\GPY)} \)\label{GGref}]
It holds \( \gv_{0} = \KS(\fv_{0}) \) and 
there exists a constant \( \CONST_{\GP|\GPY} \) such that 
\begin{EQA}
	\bigl\| \GPY \bigl\{ \KS(\fvs) - \gv_{0} \bigr\} \bigr\|^{2}
	& \leq &
	\CONST_{\GP|\GPY} \, \bigl\| \GP (\fvs - \fv_{0}) \bigr\|^{2} \, 
\label{ff0GAfg0Ga}
\label{Gff0ff0Gagg0}
\end{EQA}
and with \( \IF_{\GP}(\upsilonv) \) from \eqref{DPGP2HPGP2GG2}
\begin{EQA}
	\tr \bigl( \sigma^{2} \GPY^{2} + \Id_{\dimq} \bigr)^{-1}
	& \leq &
	\CONST_{\GP|\GPY} \, 
	\tr \bigl\{ \IF_{\GP}^{-1}(\upsilonvs_{\GP}) \, \IF(\upsilonvs_{\GP}) \bigr\} .
\label{trs2G2Iqm1CGG}
\end{EQA}
\end{description}
%
This condition explains the choice of the prior for the auxiliary parameter \( \gv \).
Effectively this prior can be obtained from the prior on the source function \( \fv \) 
by the linear mapping \( \nabla \KS(\fv_{0}) \).

Apart the basic conditions \nameref{LLref}, 
\nameref{ED0ref}, 
\nameref{GGref} we need some local smoothness properties of the expected log-likelihood 
\( \E \LL(\upsilonv) \).
Let \( \Upsilond\) be a local subset of \( \UpsilonC \). 
We only need that this set contains the concentration set 
\( \CAGP(\rr_{\GP}) \) of the estimate \( \tilde{\upsilonv}_{\GP} \); see Theorem~\ref{TtifvsGPeff} below.
%
Our results assume that the function \( \E \LL(\upsilonv) \) is three or four times Gateaux differentiable.
Results on pMLE only involve the third derivative. 
For the most advanced results about Gaussian approximation of the posterior on the class 
of centrally symmetric sets, we require four Gateaux derivatives.
Let \( F(\upsilonv) \eqdef - \E \LL(\upsilonv) \).
Define for each \( \upsilonv \in \Upsilond\), 
and any \( \uv \in \UpsilonD \), the directional derivative
\begin{EQA}
    \partial_{\uv}^{m} F(\upsilonv)
    & \eqdef &
      \frac{d^{m}}{d\rhot^{m}} F(\upsilonv + \rhot \uv) \biggr|_{\rhot=0} \, ,
    \quad 
    m = 3, 4.
\label{d3utd4utd4dt4c}
\end{EQA}
Clearly the value \( \partial_{\uv}^{m} F(\upsilonv) \) is proportional to \( \| \uv \|^{3}\), 
while 
\( \partial_{\uv}^{m} F(\upsilonv) \asymp \| \uv \|^{4} \).
Note that all the quadratic terms from the expression \eqref{ELLth122Y1} of \( \E \LL(\upsilonv) \)
cancel in \( \delta_{m} \), only the structural term \( {\lambda} \| \gv - \KS(\fv) \|^{2}/2 \) matters.
It  obviously holds for \( \uv = (\alphav,\betav) \)
\begin{EQA}
    \partial_{\uv} \| \gv - \KS(\fv) \|^{2} 
    & = &
    2 \bigl\{ \gv - \KS(\fv) \bigr\}^{\T} 
    \bigl\{ \betav - \partial_{\alphav} \KS(\fv) \bigr\} ,
    \\
    \partial^{2}_{\uv} \| \gv - \KS(\fv) \|^{2}
    & = &
    2 \bigl\| \betav - \partial_{\alphav} \KS(\fv) \bigr\|^{2}
    - 2 \bigl\{ \gv - \KS(\fv) \bigr\}^{\T}
    \partial^{2}_{\alphav} \KS(\fv),
    \\
    \partial^{3}_{\uv} \| \gv - \KS(\fv) \|^{2} 
    & = &
    - 6 \bigl\{ \betav - \partial_{\alphav} \KS(\fv) \bigr\}^{\T}
    	\partial^{2}_{\alphav} \KS(\fv) 
    - 2 \bigl\{ \gv - \KS(\fv) \bigr\}^{\T}
    \partial^{3}_{\alphav} \KS(\fv) ,
\label{d3utd4utd4dt4cm}
    \\
    \partial^{4}_{\uv} \| \gv - \KS(\fv) \|^{2} 
    & = &
    6 \bigl\| \partial^{2}_{\alphav} \KS(\fv) \bigr\|^{2}
    - 8 \bigl\{ \betav - \partial_{\alphav} \KS(\fv) \bigr\}^{\T}
    \partial^{3}_{\alphav} \KS(\fv)
    - 2 \bigl\{ \gv - \KS(\fv) \bigr\}^{\T}
    \partial^{4}_{\alphav} \KS(\fv) .
\end{EQA}
Therefore, the function \( F(\upsilonv) \) inherits the smoothness properties of 
the operator \( \KS(\fv) \).
We assume the following condition.

\begin{description}
\item[\( \bb{(\LL_{0})}\)\label{LL0ref}]
   \textit{The functions \( \partial^{3}_{\alphav} \KS(\fv) \) and 
   \( \partial^{4}_{\alphav} \KS(\fv) \) are well defined and for 
   specific sets \( \Upsilond , \UVd \subset \UpsilonD \), it holds}
\begin{EQA}
	\sup_{\fv \in \Upsilond, \alphav \in \UVd} \bigl| \partial^{m}_{\alphav} \KS(\fv) \bigr|
	& \leq &
	\CONST_{m}(\Upsilond,\UVd) \, ,
	\quad m=3,4. 
\label{finUdainVdCm34}
\end{EQA}   

\end{description}

Checking this condition can be tricky in general situation.
However, a number of results are available for particular cases; see 
e.g. in \cite{YAJIMA200481}, \cite{DANCONA20094552} and references therein 
for the case of an elliptic operator \( \KS \). 

\subsection{Effective dimension}
This section discusses the central notion of \emph{effective dimension}.
We first present the definition inspired by \cite{SP2013_rough} in context of penalized 
ML estimation. 
Define 
\begin{EQA}
	\BB_{\SiGP} 
	& \eqdef &
	\sigma^{-2} \, \epsVar \, \IGblk_{\GP}^{-1} \epsVar ,
\label{BV1Gtm1t}
\end{EQA}
where \( \epsVar \) is from \nameref{ED0ref} and \( \IGblk_{\GP}^{-1} \) is the 
\( \gv\gv \)-block of \( \IFT_{\GP}^{-1}(\upsilonvs_{\GP}) \).
One can also use the inverse of \( \IFTa_{\GP}(\upsilonvs_{\GP}) \) from \eqref{DPGP2HPGP2GG2m2} 
in place of \( \IFT_{\GP}(\upsilonvs_{\GP}) \) ignoring the term with the second derivative.
Also define 
\begin{EQA}[rcccl]
    \dimA_{\SiGP} 
    & \eqdef &
    \tr \BB_{\SiGP} \, , 
    \qquad
    \lambda_{\SiGP}
    & \eqdef & 
    \| \BB_{\SiGP} \|   \, . 
\label{lamGPDPVPfis}
\end{EQA}
Here \( \| \BB \| \) means the operator norm or the maximal eigenvalue of \( \BB \).
These values are important because they enter in the definition of the upper quantile function \( \zq^{2}(\BB_{\SiGP}, \xx)\) 
for \( \sigma^{-2} \epsv^{\T} \IGblk_{\GP}^{-1} \epsv \).
In \cite{SP2013_rough} a similarly defined quantity \( \dimA_{\SiGP} \) was called 
the \emph{effective dimension} 
in context of penalized MLE; see also Theorem~\ref{TtifvsGPeff} below.
In Bayesian framework we introduce a slightly different definition 
of effective dimension
which mimics the posterior distribution rather than the distribution of the pMLE.
%
Let \( \IFT(\upsilonv) = - \nabla^{2} \E \LL(\upsilonv) \) be the negative Hessian 
of the non-penalized log-likelihood
\( \LL(\upsilonv) = - \| \Yv - \gv \|^{2}/(2\sigma^{2}) - 
\lambda \| \gv - \KS(\fv) \|^{2} \).
The \emph{local effective dimension} \( \dimG(\upsilonv) \) at the point \( \upsilonv \in \Upsilond \) is given by
\begin{EQA}
    \dimG(\upsilonv)
    & \eqdef &
    \tr \bigl\{ \IFT(\upsilonv) \IFT_{\GP}^{-1}(\upsilonv) \bigr\} .
\label{dimGtdef}
\end{EQA}
Again, one can use here \( \IFTa(\upsilonv) \) and \( \IFTa_{\GP}(\upsilonv) \) for any \( \upsilonv \)
from a local vicinity of the point \( \upsilonvs_{\GP} \).
%
Condition \nameref{GGref} means that the such defined effective dimension of the full problem is 
of the same order as the dimension of the original problem for the parameter \( \fv \).
To make this point more clear, we consider a special case of a linear operator 
\( \KS(\fv) = \KS \fv \).
For simplicity also assume \( \lambda = \sigma^{-2} \).
Effective dimension for the original problem without calming is given by 
\begin{EQA}
	\dimA_{\fv}
	&=&
	\tr\bigl( \IF_{\GP}^{-1} \IF \bigr)
	=
	\tr \bigl\{ (\GP^{2} + \sigma^{-2} \KS^{\T} \KS)^{-1} \sigma^{-2} \KS^{\T} \KS \bigr\} ;
\label{detDGm2D2tG2sm2}
\end{EQA}
cf. \cite{SpPa2019}.
Similarly define
\begin{EQA}
	\dimA_{\gv}
	&=&
	\tr\bigl( \IG_{\GP}^{-1} \IG \bigr)
	=
	\tr \bigl((\sigma^{-2} + \lambda)^{-1} \GPY^{2} + \Id_{\dimq} \bigr)^{-1} .
\label{detDGm2D2tG2sm2}
\end{EQA}
By Lemma~\ref{LeffdimNLI}, it holds \( \IFT_{\GP}^{-1} \leq 2 \block\bigl\{ \IF_{\GP}^{-1},\IG_{\GP}^{-1} \bigr\} \) and hence
\begin{EQA}
	\dimG
	&=&
	\tr\bigl( \IFT_{\GP}^{-1} \IFT \bigr)
	\leq 
	2 \tr \bigl( \IF_{\GP}^{-1} \IF \bigr) + 2 \tr\bigl( \IG_{\GP}^{-1} \IG \bigr)
	=
	2 \dimA_{\fv} + 2 \dimA_{\gv} \, .
\label{pGGFm1F22pf2pg}
\end{EQA}
Condition \nameref{GGref} yields that
\begin{EQA}
	\tr \bigl( (\sigma^{-2} + \lambda)^{-1} \GPY^{2} + \Id_{\dimq} \bigr)^{-1}
	& \lesssim &
	\tr \bigl( \IF_{\GP}^{-1} \IF \bigr) .
\label{m1G2sm2nKf0T0m2}
\end{EQA}
This particularly means that the effective dimension is not increased in order
after calming. 
For a general nonlinear operator \( \KS \), similar bounds apply with \( \KS \)
replaced by the gradient \( \nabla \KS(\fv) \) at the point \( \fv = \fvs_{\GP} \).
This can be seen if we replace \( \IFT_{\GP} \) with \( \IFTa_{\GP} \).

\subsection{Main results}
\label{SmainresultsNI}
Without explicitly mentioned, we assume that the conditions 
\nameref{LLref},
\nameref{ED0ref}, 
\nameref{GGref}, 
and
\nameref{LL0ref}
are fulfilled.
We only specify the requirements on the local subset \( \Upsilond \) from
condition \nameref{LL0ref}. 
The proposed calming approach suggests to consider the couple \( \upsilonv = (\fv,\gv) \) and 
the corresponding penalized log-likelihood \( \LLG(\upsilonv) \).
Define 
\begin{EQA}
	\tilde{\upsilonv}_{\GP}
	=
	\bigl( \tilde{\fv}_{\GP},\tilde{\gv}_{\GP} \bigr)
	& \eqdef &
	\argmax_{\upsilonv \in \UpsilonC} \LLG(\upsilonv) ,
	\\
	\upsilonvs_{\GP}
	=
	\bigl( \fvs_{\GP},\gvs_{\GP} \bigr)
	& \eqdef &
	\argmax_{\upsilonv \in \UpsilonC} \E \LLG(\upsilonv) .
\label{tfGtgGafXPsg}
\end{EQA}
The joint posterior \( \vupsilonv_{\GP} \cond \Yv \) is defined by normalizing the exponent
\( \exp \LLG(\upsilonv) \).

Remind that \( \HFblk_{\GP}^{-1} \) is the \( \gv\gv \)-block of 
\( \IFT_{\GP}^{-1}(\upsilonvs_{\GP}) \); see \eqref{IFm1uDAATHm2}.
We now apply the deviation bound of Theorem~\ref{LLbrevelocroB}
to the quadratic form \( \sigma^{-2} \epsv^{\T} \IGblk_{\GP}^{-1} \epsv \).
Under condition \nameref{ED0ref}, for any \( \xx > 0 \) there exists a random set \( \Omega(\xx)\) with \( \P\bigl( \Omega(\xx) \bigr) \geq 1 - \CONST \ex^{-\xx}\) such that on this set
\begin{EQA}
    \sigma^{-2} \epsv^{\T} \IGblk_{\GP}^{-1} \epsv
    & \leq &
    \zq^{2}(\BB_{\SiGP},\xx) ,
\label{uTzDGvTxm12z}
\end{EQA}
where \( \BB_{\SiGP} = \sigma^{-2} \, \epsVar \, \HFblk_{\GP}^{-1} \epsVar \) and 
for any matrix \( \BB \)
\begin{EQA}
	\zq(\BB,x)
	&=&
	\sqrt{\tr \BB} + \sqrt{2 \lambda_{\max}(\BB) \, \xx} .
\label{zzxxppdBlroBB}
\end{EQA}
Here we assume \( \gm = \infty \) in \nameref{ED0ref}.
Otherwise the formula for \( \zq(\BB_{\SiGP},\xx) \) is more involved.

Our first result claims that the posterior distribution \( \fv_{\GP} \cond \Yv \)
of the target parameter \( \fv \) is nearly Gaussian.
More precisely, we present two finite sample bounds on the accuracy of Gaussian approximation. 
The first bound is limited to the class \( \cc{B}_{s}(\XSD) \) of centrally symmetric Borel sets in \( \XSD \),
while the second one is in total variation distance. 
Our results rely on smoothness of function \( \E \LL(\upsilonv) \) in terms of the third and fourth Gateaux derivatives of 
\( \| \gv - \KS(\fv) \|^{2} \).
Let \( \DF^{2}(\upsilonv) = \IFT(\upsilonv) = - \nabla^{2} \E \LL(\upsilonv) \).
Given a local set \( \Upsilond \) in \( \UpsilonD \), define
\begin{EQA}
    \delta_{m}(\Upsilond,\rr)
    & \eqdef &
    \sup_{\upsilonv \in \Upsilond} \,\, 
          \sup_{\uv \colon \| \DF(\upsilonv) \uv \| \leq \rr} \, \, 
      \Bigl| \partial^{m}_{\uv} \| \gv - \KS(\fv) \|^{2} \Bigr|,
    \quad
    m=3,4.
\label{d3sBsurCGrGu}
\end{EQA}
We suppress the argument \( \Upsilond \) and write simply \( \delta_{m}(\rups)\). 
Below in Theorem~\ref{TtifvsGPeff} we show that this bound together with local smoothness and concavity 
of the expected log-likelihood \( \E \LLG(\upsilonv) \) allow to establish sharp concentration bounds 
for the total estimator \( \tilde{\upsilonv}_{\GP} \)
in an elliptic set \( \CAGP(\rr_{\GP}) \) around \( \upsilonvs_{\GP} \).
For our main result we only need that \( \CAGP(\rr_{\GP}) \subset \Upsilond \).
We also denote 
\( \DPblkt_{\GP} = \IFblk_{\GP}^{1/2}(\tilde{\upsilonv}_{\GP}) \), 
where \( \IFblk_{\GP}^{-1}(\upsilon) \) is the \( \fvfv \)-block of \( \IFT_{\GP}^{-1}(\upsilonv) \).

\begin{theorem}
\label{PrhoQPBvMf}
Let the local set \( \Upsilond \) contain the set \( \CAGP(\rr_{\GP}) \) of Theorem~\ref{TtifvsGPeff}.
Let also, for some fixed values \( \rups \) and \( \xx > 0 \), it hold
\begin{EQA}
\label{LmgfquadELGPf}
    \err(\rups)
    & \eqdef &
    4 \delta_{3}^{2}(\rups) + 4 \delta_{4}(\rups)
    \leq 
    1/2,
    \\
    \CONSTru
    & \eqdef &
    1 - 3 \rups^{-2} \delta_{3}(\rups) 
    \geq 1/2.
    \\
    \CONSTru \rups 
    & \geq &
    2 \sqrt{\dimG(\upsilonv)} + \sqrt{\xx} \, ,
    \qquad
    \upsilonv \in \Upsilond \, 
\label{CONSTruAxx1f}
\end{EQA}
with \( \dimG(\upsilonv) \) from \eqref{dimGtdef}.
Let \eqref{uTzDGvTxm12z} hold on the random set \( \Omega(\xx)\).
Then on this set for any centrally symmetric Borel set \( A \)
\begin{EQ}[rcl]
    \P\bigl( \fv_{\GP} - \tilde{\fv}_{\GP} \in A \cond \Yv \bigr)
    & \geq &
    \frac{1 - \err(\rups) }
     	 {\bigl\{ 1 + \err(\rups) + \rho(\rups) \bigr\}} \,
     		\PG\bigl( \DPblkt_{\GP}^{-1} \gammav \in A \bigr) - \rho(\rups) \, ,
    \\
    \P\bigl( \fv_{\GP} - \tilde{\fv}_{\GP} \in A \cond \Yv \bigr)
    & \leq & 
    \frac{1 + \err(\rups)}{\bigl\{ 1 - \err(\rups) \bigr\} \bigl( 1 - \ex^{-\xx} \bigr)} 
    \PG\bigl( \DPblkt_{\GP}^{-1} \gammav \in A \bigr) + \rho(\rups) \, ,
\label{1p1m1a1emxf}
\end{EQ}
where \( \gammav \) is standard Gaussian in \( \XSD \), \( \PG \) is the conditional distribution
given \( \DPblkt_{\GP} \), and
\begin{EQA}
    \rho(\rups)
    & \leq &
    \frac{1}{1 - \err(\rups)} 
    \frac{\exp\bigl\{- (\dimGt + \xx)/2\bigr\}}
    	 {1 - \exp\{ -(\dimGt + \xx)/2 \}} \, 
\label{rhoQPubnBvmf}
\end{EQA}
with \( \dimGt = \dimG(\tilde{\upsilonv}_{\GP}) \); see \eqref{dimGtdef}.
For any Borel set \( A \), the bounds from \eqref{1p1m1a1emxf} apply with \( \delta_{3}(\rups) \)
in place of \( \err(\rups) \).
\end{theorem}

All the previous results are finite-sample with explicit error terms.
Now we introduce a large-sample parameter \( n \).
Typical example is 
inverse of noise energy \( n = \sigma^{-2} \). 

\begin{corollary}
\label{CnonparBvMmarg}
Suppose that \( \rups \) satisfies the conditions \eqref{LmgfquadELGPf} and \eqref{CONSTruAxx1f} with 
\( \xx = 2 \log n \).
It holds on \( \Omega(\xx)\) 
\begin{EQA}
\label{1p1m1a1emxrr}
	\sup_{A \in \cc{B}_{s}(\XSD)}
	\left| 
		\P\bigl( \fv_{\GP} - \tilde{\fv}_{\GP} \in A \cond \Yv \bigr)
		- \PG\bigl( \DPblkt_{\GP}^{-1} \gammav \in A \bigr)  
	\right|
	& \leq &
	\CONST \bigl\{ \err(\rups) 
	+ 1/n \bigr\}
	\\
\label{1p1m1a1emxrras}
	\sup_{A \in \cc{B}(\XSD)}
	\left| 
		\P\bigl( \fv_{\GP} - \tilde{\fv}_{\GP} \in A \cond \Yv \bigr)
		- \PG\bigl( \DPblkt_{\GP}^{-1} \gammav \in A \bigr)  
	\right|
	& \leq &
	\CONST \bigl\{ \delta_{3}(\rups) 
	+ 1/n \bigr\} .
\end{EQA}
\end{corollary}

Comparison of two bounds of Corollary~\ref{CnonparBvMmarg} reveals
that the use of symmetric credible sets improves the accuracy of Gaussian approximation
from \( \delta_{3}(\rups) \) to \( \err(\rups) \asymp \delta_{3}^{2}(\rups) + \delta_{4}(\rups) \).
In typical regular cases, \( \delta_{3}(\rups) \asymp \sqrt{\rups^{3}/n} \) and 
\( \delta_{4}(\rups) \asymp {\rups^{2}/n} \)
yielding 
\( \err(\rups) \asymp \rups^{3}/n \).
The choice \( \xx = 2 \log n \) and 
\( \rups = \CONST \bigl( \sqrt{\dimG} + \sqrt{\log n} \bigr) \) yields 
\( \rho(\rups) \leq 1/n \) in \eqref{1p1m1a1emxf},
and the only leading term in the error of approximation is 
\( \err(\rups) \asymp \dimG^{3}/n \),
and this is the guaranteed approximation error in the BvM approximation under symmetricity.
The bound in TV-distance ensures an error \( \delta_{3}(\rups) \asymp \sqrt{\dimG^{3}/n} \);
cf. \cite{SpPa2019}.

Gaussian approximation of the posterior allows to derive a number of corollaries about posterior behavior. 
First we state the concentration result.
Let \( \QP \) be a linear mapping from \( \XSD \) to \( \R^{\dimm} \) for some \( \dimm \leq \dimp \).
A canonical choice is an identity mapping or a projector on some subspace of \( \XSD \).
We are interested in the posterior distribution \( \QP \fv_{\GP} \cond \Yv \), in particular, in its concentration
set. 
For ease of notation, we state a large sample result. 
Combination of the Gaussian approximation bound of Theorem~\ref{PrhoQPBvMf} and the large deviation bound for 
Gaussian quadratic forms from Theorem~\ref{TexpbLGA} yields with \( \zq(\BB,\xx) \) from 
\eqref{zzxxppdBlroBB} the following corollary. 

\begin{corollary}
\label{CBvMNI}
Under the conditions of Theorem~\ref{PrhoQPBvMf}, it holds on \( \Omega(\xx) \) with \( \xx = \log n \)
\begin{EQA}
	\P\bigl( \| \QP (\fv_{\GP} - \tilde{\fv}_{\GP}) \| > \zq(\QP \, \DPblkt_{\GP}^{-2} \, \QP^{\T},\xx) \cond \Yv \bigr)
	& \leq &
	\CONST \bigl\{ \err(\rups) 
	+ 1/n \bigr\}. 
\label{QpostconcNI}
\end{EQA}
\end{corollary}

A more general fact about concentration of the full parameter will be given in Corollary~\ref{CBvMNIT}.
A further combination of this bound with the concentration results for the pMLE \( \tilde{\upsilonv}_{\GP} \) can be used
to get the contraction results for \( \QP (\fv_{\GP} - \fvs) \cond \Yv \); see Section~\ref{Scontractcover}. 

\section{Further results}
\label{SCompNI}
This section presents more results about the properties of the penalized MLE 
\( \tilde{\upsilonv}_{\GP} \) and of the posterior \( \vupsilonv_{\GP} \cond \Yv \)
and its marginal \( \fv_{\GP} \cond \Yv \).

\subsection{Properties of the pMLE \( \tilde{\upsilonv}_{\GP} \) and \( \tilde{\fv}_{\GP} \)}
\label{SpMLENI}
First we present a result about large deviation bound for the penalized MLE 
\( \tilde{\upsilonv}_{\GP} = (\tilde{\fv}_{\GP},\tilde{\gv}_{\GP}) \) defined by maximizing \( \LLG(\upsilonv) \).
Remind the notation \( \IFT_{\GP}(\upsilonv) = - \nabla^{2} \E \LLG(\upsilonv) \), 
\( \upsilonvs_{\GP} = \argmax_{\upsilonv \in \Upsilond} \E \LLG(\upsilonv) \),
\( \IFT_{\GP} = \IFT_{\GP}(\upsilonvs_{\GP}) \), and \( \DF_{\GP} = \sqrt{\IFT_{\GP}} \).
We show that the estimator \( \tilde{\upsilonv}_{\GP} \) concentrates in an elliptic vicinity of \( \upsilonvs_{\GP} \)
of the form
\begin{EQA}
	\CAGP(\rr) 
	& \eqdef &
	\bigl\{ \upsilonv \colon \| \DF_{\GP} (\upsilonv - \upsilonvs_{\GP}) \| \leq \rr \bigr\} 
\label{AGuvDFGsuuSr}
\end{EQA}
for a proper choice of \( \rr \).

With \( \partial^{3}_{\uv} \) from \eqref{d3utd4utd4dt4c}, 
define for each \(\rr > 0 \) 
\begin{EQA}
    \delta_{3,\GP}(\rr)
    & \eqdef &
    \sup_{\upsilonv \colon \| \DF_{\GP} (\upsilonv - \upsilonvs_{\GP}) \| \leq \rr} \,\, 
    \sup_{\uv \colon \| \DF_{\GP} \uv \| \leq \rr}
      \Bigl| \partial^{3}_{\uv} \| \gv - \KS(\fv) \|^{2} \Bigr|.
\label{d3sBsurCGrGu3u}
\end{EQA}
This value is finite under~\nameref{LL0ref} provided that \( \{\upsilonv \colon \| \DF_{\GP} (\upsilonv - \upsilonvs_{\GP}) \| \leq \rr\} \subseteq \Upsilond \) and  \( \{\uv \colon \| \DF_{\GP} \uv \| \leq \rr\} \subseteq \UVd \).

\begin{theorem}
\label{TtifvsGPeff}
Let \eqref{uTzDGvTxm12z} hold on a random set \(\Omega(\xx)\) with
\(\P\bigl( \Omega(\xx) \bigr) \geq  1 - \ex^{-\xx}\).
Let also \(\rr_{\GP}\) be such that \( \Upsilond \) contain the set 
\( \CAGP(\rr_{\GP}) \) from \eqref{AGuvDFGsuuSr} and  
\begin{EQ}[rcl]
    \frac{3 \delta_{3,\GP}(\rr_{\GP})}{\rr_{\GP}^{2}} 
    & \leq &
    \rho \leq 1/2, 
    \\
    (1 - \rho) \rr_{\GP}
    & \geq &
    \zq(\BB_{\SiGP},\xx) .
\label{3rGm2d3rG12}
\end{EQ}
Then on \(\Omega(\xx)\), the estimate 
\( \tilde{\upsilonv}_{\GP} \) belongs to the set \( \CAGP(\rr_{\GP}) \), that is, 
\begin{EQA}
	\bigl\| \DF_{\GP} \bigl( \tilde{\upsilonv}_{\GP} - \upsilonvs_{\GP} \bigr) \bigr\|
	& \leq &
	\rr_{\GP}
	=
	(1 - \rho)^{-1} \zq(\BB_{\SiGP},\xx) .
\label{1mrjm1zBVGx}
\end{EQA}
\end{theorem}
  
Note that the concentration of pMLE can be stated under a weaker condition on the effective dimension, namely, 
we only need that \( \delta_{3,\GP}(\rr_{\GP}) \ll \rr_{\GP}^{2} \)
which results in \( \dimA_{\GP} \ll n \); cf. \cite{SpPa2019}.

Due to the concentration result of Theorem~\ref{TtifvsGPeff}, the estimate 
\(\tilde{\upsilonv}_{\GP}\) lies with a dominating probability in a local vicinity 
of the point \(\upsilonvs_{\GP}\).
Now one can use a quadratic approximation for the penalized log-likelihood process 
\( \LL_{\GP}(\upsilonv) \) to establish an expansion for the penalized MLE
\(\tilde{\upsilonv}_{\GP}\) and for the excess \(\LL_{\GP}(\tilde{\upsilonv}_{\GP}) - \LL_{\GP}(\upsilonvs_{\GP})\).
Remind that \( \nabla \zeta = \bigl( 0, \sigma^{-1} \epsv \bigr) \) is the score vector in the full model.

\begin{theorem}
\label{TFiWititG}
Under the conditions of Theorem~\ref{TtifvsGPeff}, it holds on \(\Omega(\xx)\) 
\begin{EQA}
    \bigl\| \DF_{\GP} \bigl( \tilde{\upsilonv}_{\GP} - \upsilonvs_{\GP} \bigr) - \DF_{\GP}^{-1} \nabla \zeta \bigr\|^{2}
    & \leq &
    4 \delta_{3,\GP}(\rr_{\GP}) .
\label{DGttGtsGDGm13rG}
    \\
    \biggl| 
    \LL_{\GP}(\tilde{\upsilonv}_{\GP}) - \LL_{\GP}(\upsilonvs_{\GP}) 
    - \frac{1}{2} \bigl\| \DF_{\GP}^{-1} \nabla \zeta \bigr\|^{2}
    \biggr|
    & \leq &
    \delta_{3,\GP}(\rr_{\GP}) ,
    \\
    \biggl| 
    \LL_{\GP}(\tilde{\upsilonv}_{\GP}) - \LL_{\GP}(\upsilonvs_{\GP}) 
    - \frac{1}{2} \bigl\| \DF_{\GP} (\tilde{\upsilonv}_{\GP} - \upsilonvs_{\GP}) \bigr\|^{2}
    \biggr|
    & \leq &
    \delta_{3,\GP}(\rr_{\GP}) ,
\label{3d3Af12DGttG}
\end{EQA}
and also, for any \( \upsilonv \in \CAGP(\rr_{\GP}) \),
\begin{EQA}
    \biggl| 
    \LL_{\GP}(\tilde{\upsilonv}_{\GP}) - \LL_{\GP}(\upsilonv) 
    - \frac{1}{2} \bigl\| \DFt_{\GP} (\tilde{\upsilonv}_{\GP} - \upsilonv) \bigr\|^{2}
    \biggr|
    & \leq &
    \delta_{3,\GP}(\rr_{\GP}) ,
\label{3d3Af12DGttt}
\end{EQA}
where the random matrix \(\DFt_{\GP}^{2} = \IFT_{\GP}(\tilde{\upsilonv}_{\GP})\) fulfills
on \( \Omega(\xx) \) for some universal constant \( \CONST \)
\begin{EQA}
    \bigl\| \DF_{\GP}^{-1} \bigl( \DFt_{\GP}^{2} - \DF_{\GP}^{2} \bigr) \DF_{\GP}^{-1} \bigr\|
    & \leq &
    \CONST \rr_{\GP}^{-2} \delta_{3,\GP}(\rr_{\GP}) .
\label{DPGPm1Cd3rG}
\end{EQA}
\end{theorem}

Similarly to Theorem~\ref{TtifvsGPeff}, the results of Theorem~\ref{TFiWititG} are 
meaningful if \( \delta_{3,\GP}(\rr_{\GP}) \ll \rr_{\GP}^{2} \), that is, if 
\( \dimA_{\SiGP} \) is significantly smaller than \( n \).

Special structure of the score vector \( \nabla \zeta = \sigma^{-1} (0,\epsv) \) 
together with the CLT for the standardized score \( \epsVar^{-1/2} \epsv \) and \eqref{DGttGtsGDGm13rG}, 
can be used to establish asymptotic normality of \( \tilde{\upsilonv}_{\GP} - \upsilonvs_{\GP} \), and hence, 
of \( \tilde{\fv}_{\GP} - \fvs_{\GP} \). 
Also, the expansion \eqref{DGttGtsGDGm13rG} for the pMLE \( \tilde{\upsilonv}_{\GP} \) can be used to extend
the deviation bound of Theorem~\ref{TtifvsGPeff} to \( \bigl\| \QF \bigl( \tilde{\upsilonv}_{\GP} - \upsilonvs_{\GP} \bigr) \bigr\| \)
for any linear mapping \( \QF \to \R^{m} \). 
We only present a result for the \( \fv \)-component of \( \upsilonv \).
Remind that \( \DPblk_{\GP}^{-2} \) and 
\( \AFblk_{\GP} \) are the \( \fvfv \) and \( \fv\gv \)-blocks of the matrix 
\( \IFT_{\GP}^{-1} \),
and \( \epsVar^{2} \) bounds the variance of \( \epsv \); see \nameref{ED0ref}.

\begin{corollary}
\label{CBvMNIML}
Let \( \QF \upsilonv = \QP \fv \). 
Define \( \BB_{\QPGP} = \sigma^{-2} \QP \, \AFblk_{\GP} \, \epsVar^{2} \, \AFblk_{\GP}^{\T} \, \QP^{\T} \).
Under the conditions of Theorem~\ref{TtifvsGPeff}, it holds on \( \Omega_{1}(\xx) \) with 
\( \P\bigl( \Omega_{1}(\xx) \bigr) \leq 2/n \) 
\begin{EQA}
	\| \QP (\tilde{\fv}_{\GP} - \fvs_{\GP}) \| 
	& \leq &
	\zq(\BB_{\QPGP},\xx) + \| \QP \, \DPblk_{\GP}^{-1} \| \sqrt{\delta_{3,\GP}(\rr_{\GP})} .
\label{QfGtfGzBQGxm1}
\end{EQA}
\end{corollary}

\subsection{The use of \( \IFTa_{\GP} \) and \( \IFa_{\GP} \)}
The presented results involve the inverse of 
\( \IFT_{\GP} = \IFT_{\GP}(\upsilonvs_{\GP}) 
= - \nabla^{2} \LL_{\GP}(\upsilonvs_{\GP}) \).
An issue in the analysis of this matrix is the term with the second derivative 
of the operator \( \KS \).
It appears that this term can be omitted and one can use the matrix \( \IFTa_{\GP} \)
from \eqref{DPGP2HPGP2GG2m2}
instead of \( \IFT_{\GP} \) in most of results.
This also concerns the inverse \( \IFT_{\GP} \) and its blocks \( \DPblk_{\GP}^{-2} = \IFblk_{\GP}^{-1} \),
\( \IGblk_{\GP}^{-1} \), and \( \AFblk_{\GP} \).
All these objects can be replaced by its breve-version based on \( \IFTa_{\GP} \).
The main reason is that the elasticity vector \( \deltavd = \gv - \KS(\fv) \) nearly 
vanishes at \( \upsilonvs_{\GP} \).

Moreover, the norm \( \bigl\| \QP \, \DPblka_{\GP}^{-2} \QP^{\T} \bigr\| \) can be replaced 
by the norm \( \bigl\| \QP \, \IFa_{\GP}^{-1} \QP^{\T} \bigr\| \) up to a constant factor,
where \( \IFa_{\GP} = \GP^{2} + \sigma^{-2} \nabla \KS^{\T} \nabla \KS(\fvs_{\GP}) \). 
The latter corresponds to an approximation of the original nonlinear operator \( \KS \) by its 
linearization at \( \fvs_{\GP} \).
Similarly, the block \( \AFblk_{\GP} \) in the matrix 
\( \BB_{\QPGP} = \sigma^{-2} \QP \, \AFblk_{\GP} \, \epsVar^{2} \, \AFblk_{\GP}^{\T} \, \QP^{\T} \)
can be replaced by its breve-analog \( \AFblka_{\GP} \) based on \( \IFTa_{\GP} \).
For any two matrices \( \BB \) and \( \breve{\BB} \), denote 
\begin{EQA}
	\Delta(\BB,\breve{\BB})
	& \eqdef &
	\bigl\| \BB^{-1/2} \breve{\BB} \, \BB^{-1/2} - \Id \bigr\| . 
\label{DeltBbB}
\end{EQA}
Obviously, the bound \( \Delta(\BB,\breve{\BB}) \leq \rho < 1 \) implies for any 
linear operator \( \QP \) and any \( \xx > 0 \)
\begin{EQA}
	\biggl| \frac{\bigl\| \QP \, \breve{\BB} \QP^{\T} \bigr\|}
		 {\bigl\| \QP \, \BB \QP^{\T} \bigr\|} - 1
	\biggr|
	\leq 
	\rho,
	& \quad &
	\biggl| \frac{\zq( \QP \, \breve{\BB} \QP^{\T},\xx)}
		 {\zq(\QP \, \BB \QP^{\T},\xx)} - 1
	\biggr|
	\leq 
	\sqrt{\rho}. 
\label{zqrQBbBzqm1}
\end{EQA}

\begin{lemma}
\label{LIFTaIFa}
It holds under \eqref{3rGm2d3rG12}
\begin{EQA}
	\Delta\bigl( \IFT_{\GP}^{-1},\IFTa_{\GP}^{-1} \bigr)
	\leq 
	\rho,
	& \quad &
	\Delta\bigl( \IFT_{\GP}^{-1},\IFTt_{\GP}^{-1} \bigr)
	\leq 
	\rho,
	\quad
	\Delta\bigl( \IFTa_{\GP}^{-1},\block\bigl\{ \IFa_{\GP}^{-1},\IG_{\GP}^{-1} \bigr\} \bigr)
	\leq 
	1/2. 
\label{12IFGm1GFm1G1}
\end{EQA}
Similar bounds hold for the blocks \( \IFblk_{\GP}^{-1} = \DPblk_{\GP}^{-2} \) 
and \( \AFblk_{\GP} \, \epsVar^{2} \, \AFblk_{\GP}^{\T} \) and their breve analogs. 
\end{lemma}

\subsection{BvM for the total parameter}
By Theorem~\ref{TtifvsGPeff}, on the set \( \Omega(\xx) \), 
the pMLE \( \tilde{\upsilonv}_{\GP} \) concentrates on the elliptic vicinity set 
\( \CAGP(\rr_{\GP}) = \bigl\{ \upsilonv \colon \| \DF_{\GP} (\upsilonv - \upsilonvs_{\GP}) \| \leq \rr_{\GP} \bigr\} \) of \( \tilde{\upsilonv}_{\GP} \) of \( \upsilonvs_{\GP} \),
where \( \rr_{\GP} \leq 2 \zq(\BB_{\SiGP},\xx) \); see \eqref{uTzDGvTxm12z}.
Our first result describes the concentration properties of the 
full posterior \( \upsilonv \cond \Yv \) and of the marginal posterior
\( \fv \cond \Yv \).
Given some \( \rups \), introduce an elliptic set 
of the form
\begin{EQA}
	\CAt(\rups) 
	&=& 
	\bigl\{\upsilonv \colon \|\DFt \upsilonv\| \leq \rups \bigr\} 
	\subset \UpsilonD \, ,
\label{CAr0CAr0fDfr0}
\end{EQA}
with \( \DFt^{2} = \IFT(\tilde{\upsilonv}_{\GP}) \).
%
First we bound from above the random quantity
\begin{EQA}
    \rho(\rups)
    & \eqdef &
    \frac{\int \Ind\bigl( \| \DFt \uv \| > \rups \bigr)
          \exp \bigl\{ \LL_{\GP}(\tilde{\upsilonv}_{\GP} + \uv) \bigr\} d \uv}
       	 {\int \Ind\bigl( \| \DFt \uv \| \leq \rups \bigr) 
          \exp \bigl\{ \LL_{\GP}(\tilde{\upsilonv}_{\GP} + \uv) \bigr\} d \uv} \, .
\label{rhopipoprDGP}
\end{EQA}
Obviously \( \P\bigl( \upsilonv_{\GP} - \tilde{\upsilonv}_{\GP} \not\in \CAt(\rups) \cond \Yv \bigr) \le \rho(\rups)\).
Therefore, small values of \( \rho(\rups)\) indicate a concentration of 
\( \upsilonv_{\GP} - \tilde{\upsilonv}_{\GP} \cond \Yv\) on the set \( \CAt(\rups)\).
We show that the choice \( \rups^{2} \geq \CONST \dimGt \)
ensures the desirable concentration, where 
\( \dimGt = \dimG(\tilde{\upsilonv}_{\GP}) \) is the total effective dimension.
Let \( \Upsilond \) be an open subset of \( \Upsilon \) that contains 
the concentration set 
\( \CAGP(\rr_{\GP}) = \bigl\{ \upsilonv \colon \| \DF_{\GP} (\upsilonv - \upsilonvs_{\GP}) \| \leq \rr_{\GP} \bigr\} \) of \( \tilde{\upsilonv}_{\GP} \);
see Theorem~\ref{TtifvsGPeff}.
Let \( \delta_{m}(\rups) = \delta_{m}(\Upsilond,\rups) \) be defined in \eqref{d3sBsurCGrGu}.

\begin{theorem}
\label{PrhoQPBvM}
Let conditions of Theorem~\ref{TtifvsGPeff} be satisfied. 
Let, for some fixed values \( \rups\) and \( \xx > 0\), it hold
\begin{EQA}
\label{LmgfquadELGP}
    \err(\rups)
    & \eqdef &
    4\delta_{3}^{2}(\rups) + 4 \delta_{4}(\rups)
    \leq 
    1/2,
    \\
    \CONSTru
    & \eqdef &
    1 - 3 \rups^{-2} \delta_{3}(\rups) 
    \geq 1/2.
    \\
    \CONSTru \rups 
    & \geq &
    2 \sqrt{\dimG(\upsilonv)} + \sqrt{\xx} \, ,
    \qquad
    \upsilonv \in \Upsilond \, . 
\label{CONSTruAxx1}
\end{EQA}
Then, on the random set \( \Omega(\xx)\) on which \eqref{uTzDGvTxm12z} holds, 
the quantity \( \rho(\rups)\) from \eqref{rhopipoprDGP} fulfills 
\begin{EQA}
    \rho(\rups)
    & \leq &
    \frac{1}{1 - \err(\rups)} 
    \frac{\exp\bigl\{- (\dimGt + \xx)/2\bigr\}}
    	 {1 - \exp\{ -(\dimGt + \xx)/2 \}} \, .
\label{rhoQPubnBvm}
\end{EQA}
\end{theorem}

\begin{remark}
The result of Theorem~\ref{PrhoQPBvM} is meaningful only if the condition \( \err(\rups) \leq 1/2 \) is fulfilled.
This condition poses certain constraints on the smoothness of the operator \( \KS \) and the underlying function 
\( \fs \) and also on the effective dimension of the problem.
\cite{SP2013_rough,SpPa2019} showed that for typical examples like generalized regression or log-density models, 
this condition can be rewritten as 
\( \dimA_{\GP}^{3} \ll n \) where \( n \) is the sample size.
\end{remark}

\begin{remark}
The set \( \bigl\{ \uv \colon \| \DF(\upsilonv) \uv \| \leq \rr \bigr\} \) in \eqref{d3sBsurCGrGu}
can be large in directions where the eigenvalues of \( \IFT(\upsilonv) \) are small. 
One can slightly modify the definition of this set.
Let us fix \( \GP_{0}^{2} \) such that \( \GP_{0}^{2} \leq \GP^{2} \) and define \( \GP_{1}^{2} = \GP^{2} - \GP_{0}^{2} \).
Informally we split the penalty term 
\( \| \GP \fv \|^{2} = \| \GP_{0} \fv \|^{2} + \| \GP_{1} \fv \|^{2} \).
The \( \GP_{0} \) part is used to upgrade \( \IF_{\GP}(\upsilonv) \) while \( \GP_{1} \)-penalty replaces 
the original \( \GP \)-penalty.
This leads to a slight increase of the effective dimension from \( \tr \bigl( \IF \IF_{\GP}^{-1} \bigr) \)
to \( \tr \bigl( \IF_{\GP_{0}} \IF_{\GP_{1}}^{-1} \bigr) \) but the set 
\( \bigl\{ \uv \colon \| \DF(\upsilonv) \uv \| \leq \rr \bigr\} \) is replaced with the set
\( \bigl\{ \uv \colon \| \DF_{\GP_{0}}(\upsilonv) \uv \| \leq \rr \bigr\} \).
\end{remark}

%
The concentration result can be restated in the form that the centered posterior 
\( \upsilonv_{\GP} - \tilde{\upsilonv}_{\GP} \cond \Yv \) concentrates on the random set 
\( \CAt(\rups) \) from \eqref{CAr0CAr0fDfr0}.
%
Now we aim to show that, after restricting to this set, the posterior can be well approximated
by a Gaussian distribution \( \ND\bigl( \tilde{\upsilonv}_{\GP}, \DFt_{\GP}^{-2} \bigr)\).
In what follows we use that \( \tilde{\upsilonv}_{\GP} \) is random on the original probability space, however, it can be considered as fixed under the posterior measure. 
By \( \PG \) we denote a standard normal distribution of a random vector 
\( \gammav \in \UpsilonD \) given \( \DFt_{\GP} = \DF_{\GP}(\tilde{\upsilonv}_{\GP}) \). 
In our results we distinguish between the class \( \cc{B}_{s}(\UpsilonD) \) of centrally symmetric Borel sets and the class \( \cc{B}(\UpsilonD) \) of all Borel sets in \( \UpsilonD \).

\begin{theorem}
\label{TnonparBvMm}
Let the conditions of Theorem~\ref{PrhoQPBvM} hold and \( \rho(\rups) \) satisfy~\eqref{rhoQPubnBvm}.
It holds on the set \( \Omega(\xx)\) from Theorem~\ref{TtifvsGPeff} 
for any centrally symmetric Borel set \(A \in \cc{B}_{s}(\UpsilonD)\) 
\begin{EQA}
    \P\bigl( \upsilonv_{\GP} - \tilde{\upsilonv}_{\GP} \in A \cond \Yv \bigr)
    & \geq &
    \frac{1 - \err(\rups) }
     	 {\bigl\{ 1 + \err(\rups) + \rho(\rups) \bigr\}} \,
     		\PG\bigl( \DFt_{\GP}^{-1} \gammav \in A \bigr) - \rho(\rups) \, ,
    \\
    \P\bigl( \upsilonv_{\GP} - \tilde{\upsilonv}_{\GP} \in A \cond \Yv \bigr)
    & \leq & 
    \frac{1 + \err(\rups)}{\bigl\{ 1 - \err(\rups) \bigr\} \bigl( 1 - \ex^{-\xx} \bigr)} 
    \PG\bigl( \DFt_{\GP}^{-1} \gammav \in A \bigr) + \rho(\rups) \, ,
\label{1p1m1a1emx}
\end{EQA}
For any measurable set \( A \in \cc{B}(\UpsilonD) \), similar bounds 
hold with \( \delta_{3}(\rups) \) in place of \( \err(\rups) \).
\end{theorem}

The first result of the theorem for 
can be represented in the form
\begin{EQA}
    && \nquad
    \left| 
    {\P\bigl( \upsilonv_{\GP} - \tilde{\upsilonv}_{\GP} \in A \cond \Yv \bigr)}
    - {\PG\bigl( \DFt_{\GP}^{-1} \gammav \in A \bigr)} 
    \right|
    \\
    & \lesssim &
    {\PG\bigl( \DFt_{\GP}^{-1} \gammav \in A \bigr)}
    \bigl\{ \err(\rups) + \ex^{-\xx} \bigr\} + \rho(\rups).
\label{PvtGttGAcY}
\end{EQA}
The second statement of the theorem for any \( A \in \cc{B}(\UpsilonD) \) allows
to bound the distance in total variation between the posterior and its Gaussian approximation
\( \DFt_{\GP}^{-1} \gammav \).

We also present a large sample bound.
The next result extends Corollary~\ref{CBvMNI}.

\begin{corollary}
\label{CBvMNIT}
\label{CnonparBvMm}
Suppose that \( \rups \) satisfies the conditions \eqref{LmgfquadELGP} and \eqref{CONSTruAxx1} with 
\( \xx = 2 \log n \).
It holds on \( \Omega(\xx)\) 
\begin{EQA}
\label{1p1m1a1emxrr}
	\sup_{A \in \cc{B}_{s}(\UpsilonD)}
	\left| 
		\P\bigl( \upsilonv_{\GP} - \tilde{\upsilonv}_{\GP} \in A \cond \Yv \bigr)
		- \PG\bigl( \DFt_{\GP}^{-1} \gammav \in A \bigr)  
	\right|
	& \leq &
	\CONST \bigl\{ \err(\rups) 
	+ 1/n \bigr\}
	\\
\label{1p1m1a1emxrras}
	\sup_{A \in \cc{B}(\UpsilonD)}
	\left| 
		\P\bigl( \vupsilonv_{\GP} - \tilde{\upsilonv}_{\GP} \in A \cond \Yv \bigr)
		- \PG\bigl( \DFt_{\GP}^{-1} \gammav \in A \bigr)  
	\right|
	& \leq &
	\CONST \bigl\{ \delta_{3}(\rups) 
	+ 1/n \bigr\} .
\end{EQA}
%
Let also \( \QF \) be a linear mapping from \( \UpsilonD \) to \( \R^{m} \).
Then it holds on \( \Omega(\xx) \) 
\begin{EQA}
	\P\bigl( \| \QF (\vupsilonv_{\GP} - \tilde{\upsilonv}_{\GP}) \| > \zq(\QF \, \DFt_{\GP}^{-2} \, \QF^{\T},\xx) \cond \Yv \bigr)
	& \leq &
	\CONST \bigl\{ \err(\rups) 
	+ 1/n \bigr\}. 
\label{QpostconcNIT}
\end{EQA}
\end{corollary}

\begin{proof}
The result follows from \eqref{PvtGttGAcY} with \( \xx = \log n \) and from the deviation bound of Theorem~\ref{TexpbLGA}.
\end{proof}



\subsection{Smoothness and bias}
\label{SsmoothNI}

Now we discuss the bias induced by the double penalization 
\( \bigl\| \GP \fv \bigr\|^{2} + \bigl\| \GPY \gv \bigr\|^{2} \).
Without loss of generality assume \( \upsilonv_{0} = 0 \).
This is just a reparametrization that helps to simplify our notation. 
Define \( \GPT^{2} = \block\bigl\{ \GP^{2},\GPY^{2} \bigr\} \) and 
\begin{EQA}
	\| \GPT \upsilonv \|^{2}
	&=&
	\bigl\| \GP \fv \bigr\|^{2} + \bigl\| \GPY \gv \bigr\|^{2} \, .
\label{222GGYfg}
\end{EQA}
The concentration set \( \CAGP(\rr_{\GP}) \) becomes smaller when \( \GPT^{2} \) increases. 
In particular, if \( \GPT^{2} \) is large then \( \tilde{\upsilonv}_{\GP} \) concentrates 
in a small vicinity of \( \upsilonvs_{\GP} \).
At the same time, 
penalization \( \| \GPT \upsilonv \|^{2} \)  
yields some estimation bias measured by 
\( \E \LL_{\GP}(\upsilonvs_{\GP}) - \E \LL_{\GP}(\upsilonvs) \) and 
\( \upsilonvs_{\GP} - \upsilonvs \).
The bias is not critical if the underlying truth \( \upsilonvs = (\fvs,\gvs) \) 
with \( \gvs = \KS(\fvs) \) is ``smooth'', that is, 
\( \| \GPT \upsilonvs \|^{2} \) is not too big.
%
Smoothness properties of the source function \( \fvs \) is reflected by the prior covariance 
\( \GP^{2} \)
via the penalty \( \bigl\| \GP \fvs \bigr\|^{2} \). 
Similarly, smoothness of the image \( \gvs = \KS(\fvs) \) has to be reflected 
by the prior choice in terms of \( \bigl\| \GPY \gvs \bigr\|^{2} \). 
Effectively we require that \( \GPY^{2} \) is selected in a way that 
the roughness penalties \( \bigl\| \GP \fvs \bigr\|^{2} \)
and \( \bigl\| \GPY \gvs \bigr\|^2 \) are of the same order; see \nameref{GGref}.
Define 
\begin{EQA}
\label{Lth122Y1222c}
	\ELL(\fv,\gv)
	& \eqdef &
	\sigma^{-2} \| \gvs - \gv \|^{2} + \lambda \| \gv - \KS(\fv) \|^{2}
	+ \bigl\| \GP \fv \bigr\|^{2} 
	+ \bigl\| \GPY \gv \bigr\|^{2} \, ,
	\qquad
	\\
	(\fvs_{\GP},\gvs_{\GP})
	&=&
	\argmin_{(\fv,\gv) \in \UpsilonC} \ELL(\fv,\gv).	
\label{ELfgs2fGgGE}
\end{EQA}
First we show that \( (\fvs,\gvs) \) is a reasonable value for minimizing 
the functional \( \ELL(\fv,\gv) \) from \eqref{Lth122Y1222c}.
Indeed, the fidelity term 
\( \sigma^{-2} \bigl\| \gvs - \gv \bigr\|^{2} \)
as well as the structural term 
\( \lambda \| \gv - \KS(\fv) \|^{2} \) vanish wenn \( \fv = \fvs \) and 
\( \gv = \gvs \),
and only penalty terms \( \bigl\| \GP \fvs \bigr\|^{2} \)
and \( \bigl\| \GPY \gvs \bigr\|^{2} \) are still active
in the value \( \ELL(\fvs,\gvs) \).
So, smoothness of \( \fvs \) and \( \gvs \) make the value 
\( \ELL(\fvs,\gvs) \) sufficiently small. 
Definition of \( (\fvs_{\GP},\gvs_{\GP}) \) implies
\( \ELL(\fvs,\gvs) \geq \ELL(\fvs_{\GP},\gvs_{\GP}) \).
This yields in particular
\begin{EQA}
	\sigma^{-2} \bigl\| \gvs - \gvs_{\GP} \bigr\|^{2}
	+ \lambda \bigl\| \KS(\fvs_{\GP}) - \gvs_{\GP} \bigr\|^{2}
	& \leq &
	\| \GPT \upsilonvs \|^{2} \, .
\label{GaGPff0sm2}
\end{EQA}
The calming approach suggests to use \( \fvs_{\GP} \) as a proxi for 
\( \fvs \).
This could be possible if \( \KS(\fvs) \approx \KS(\fvs_{\GP}) \).
The next result justifies this relation. 
For simplicity we set \( \lambda = \sigma^{-2} \).

\begin{theorem}
\label{TLbiasGPGPY}
It holds
\begin{EQA}
	\E \LL_{\GP}(\upsilonvs_{\GP}) - \E \LL_{\GP}(\upsilonvs)
	& \leq &
	\frac{1}{2} \| \GPT \upsilonvs \|^{2} 
	\leq 
	\frac{\CONST_{\GP|\GPY} + 1}{2} \| \GP \fvs \|^{2}\, .
\label{122GttsELG}
\end{EQA}
\end{theorem}

\begin{proof}
The result follows directly from \eqref{GaGPff0sm2} and \eqref{Gff0ff0Gagg0}.
\end{proof}

The next result presents an explicit sharp bound on the bias induced by the penalty
\( \| \GPT \upsilonvs \|^{2} \) and based on the local quadratic approximation 
of the expected log-likelihood.
Remind the notation \( \IFT_{\GP} = \IFT_{\GP}(\upsilonvs_{\GP}) \),
\( \DF_{\GP} = \sqrt{\IFT_{\GP}} \), 
\( \DF = \sqrt{\IFT} \) and
\( \DPblk_{\GP}^{-2} = \IFblk_{\GP}^{-1} \)
for the \( \fvfv \)-block of \( \IFT_{\GP}^{-1} \) and
the matrix \( \IFa_{\GP} \) is defined as
\( \IFa_{\GP} = \GP^{2} + \sigma^{-2} \nabla \KS^{\T} \nabla \KS (\fvs_{\GP}) \).

\begin{theorem}
\label{TbiasGP}
Let \( \| \GPT \upsilonvs \|^{2} \leq \rrbias^{2}/2 \) for some \( \rrbias \) 
such that \( \delta_{3,\GP}(\rrbias)/\rrbias^{2} \leq 1/2 \). 
Then 
\begin{EQA}
\label{fr33GrGd3rELG21}
	\Bigl| 
		\E \LL_{\GP}(\upsilonvs_{\GP}) - \E \LL_{\GP}(\upsilonvs) - \bigl\| \DF_{\GP} (\upsilonvs_{\GP} - \upsilonvs) \bigr\|^{2}/2 
	\Bigr|
	& \leq &
	\delta_{3,\GP}(\rrbias),
\end{EQA}
and 
\begin{EQA}
	\bigl\| 
		\DF^{-1} \DF_{\GP}^{2} \bigl( \upsilonvs - \upsilonvs_{\GP} - \DF_{\GP}^{-2} \GPT^{2} \upsilonvs \bigr) 
	\bigr\|^{2}
	& \leq & 
	4 \delta_{3,\GP}(\rrbias) \, .
	\qquad
\label{fr33GrGd3rELG22}
\end{EQA}
Moreover, for any linear mapping \( \QF \) in \( \R^{\dimtotal} \), it holds
\begin{EQA}
	\| \QF (\upsilonvs_{\GP} - \upsilonvs) \|
	& \leq &
	\bigl\| \QF \, \DF_{\GP}^{-2} \QF^{\T} \bigr\|^{1/2} 
	\left( \| \GPT \upsilonvs \| + 2  \sqrt{\delta_{3,\GP}(\rrbias)} \right) \, .
\label{22GtsQDGm2}
\end{EQA}
If \( \QF \upsilonv = \QP \fv \), then
\begin{EQA}
	\| \QP (\fvs_{\GP} - \fvs) \|
	& \lesssim &
	\bigl\| \QP \, \IFa_{\GP}^{-1} \QP^{\T} \bigr\|^{1/2} 
	\left( \| \GP \fvs \| + \sqrt{\delta_{3,\GP}(\rrbias)} \right) \, .
\label{22GtsQDGm23}
\end{EQA}
\end{theorem}

\subsection{Bias-variance decomposition and posterior conraction}
\label{Scontractcover}

In this section, we bring together all the previous results  
to bound the accuracy of estimation \( \tilde{\fv}_{\GP} - \fvs \) and 
the posterior deviations \( \fv_{\GP} - \fvs \).
We use that \( \fv \) is a subvector of \( \upsilonv \) and 
fix some linear mapping \( \QF \colon \UpsilonD \to \R^{m} \) such that 
\( \QF \) only depends on \( \fv \)-subvector and thus, 
\( \QF \upsilonv = \QP \fv \) where \( \QP \) is a linear mapping from \( \XSD \) to \( \R^{m} \).
One can apply \( \QF \upsilonv = \IF^{1/2} \fv \) for prediction and 
\( \QF \upsilonv = \fv \) for estimation.
We first aim at stating an analog of classical bias-variance decomposition of the loss 
\( \| \QP (\tilde{\fv}_{\GP} - \fvs) \| 
= \| \QP (\tilde{\fv}_{\GP} - \fvs_{\GP}) + \QP (\fvs_{\GP} - \fvs) \| \). 
%
Remind the notation 
\( \BB_{\QPGP} = \sigma^{-2} \QP \, \AFblk_{\GP} \, \epsVar^{2} \, \AFblk_{\GP}^{\T} \, \QP^{\T} \),
where \( \AFblk_{\GP} \) is the \( \fv\gv \)-block of the matrix 
\( \IFT_{\GP}^{-1} \)
and \( \epsVar^{2} \) bounds the variance of \( \epsv \); see \nameref{ED0ref}.
Combination of the probabilistic bound \eqref{QfGtfGzBQGxm1} for \( \| \QP (\tilde{\fv}_{\GP} - \fvs_{\GP}) \| \)
and the bound \eqref{22GtsQDGm23} on the bias term \( \| \QP (\fvs_{\GP} - \fvs) \| \) yields
the following bound for the loss of \( \tilde{\fv}_{\GP} \).

\begin{theorem}
\label{TestlosspMLE}
Let \( \QF \upsilonv = \QP \fv \).
On a random set \( \Omega_{1}(\xx) \) with \( \P\bigl( \Omega_{1}(\xx) \bigr) \leq 2 \ex^{-\xx} \),
it holds for some fixed constant \( \CONST \)
\begin{EQA}
	\| \QP (\tilde{\fv}_{\GP} - \fvs) \|
	& \leq &
	\CONST \bigl\| \QP \, \IFa_{\GP}^{-1} \QP^{\T} \bigr\|^{1/2} \, \, \| \GP \fvs \| 
	+ \zq(\BB_{\QPGP},\xx) .
\label{22GtsQDGm2loss}
\end{EQA}
\end{theorem}

%
In this result we omited the small \( \delta_{3} \)-term from \eqref{22GtsQDGm23} assuming that it is hidden in the bias term.
The term \( \zq(\BB_{\QPGP},\xx) \) in \eqref{22GtsQDGm2loss} provides a sharp bound for 
the norm of the stochastic component \( \QP (\tilde{\fv}_{\GP} - \fvs_{\GP}) \) 
of \( \tilde{\fv}_{\GP} \), while 
\( \bigl\| \QP \, \IFa_{\GP}^{-1} \QP^{\T} \bigr\|^{1/2} \, \, \| \GP \fvs \| \)
is the order of the bias. 
The bias-variance trade-off corresponds to the relation
\begin{EQA}
	\bigl\| \QP \, \IFa_{\GP}^{-1} \QP^{\T} \bigr\|^{1/2} \, \, \| \GP \fvs \|
	& \asymp &
	\zq(\BB_{\QPGP},\xx) .
\label{DtsGtsrGafstar}
\end{EQA}
Suppose that smoothness properties of \( \fvs \) can be measured by \( \| \GP \fvs \|^{2} \),
namely,
\( \| \GP \fvs \|^{2} \leq \GPUB\).
Then the \emph{bias-variance trade-off} writes as
\begin{EQA}
	\GPUB \, \bigl\| \QP \, \IFa_{\GP}^{-1} \QP^{\T} \bigr\|  
	& \asymp &
	\zq^{2}(\BB_{\QPGP},\xx) .
\label{DtsGtsrGa}
\end{EQA}

Now we aim at bounding the distance between the support of the posterior and the true value \( \upsilonvs \).
The difference \( \fv_{\GP} - \fvs \) can be decomposed as
\begin{EQA}
	\fv_{\GP} - \fvs 
	& = &
	\bigl( \fv_{\GP} - \tilde{\fv}_{\GP} \bigr) + \bigl( \tilde{\fv}_{\GP} - \fvs \bigr) .
\label{vtvGtsGGtss}
\end{EQA}
Theorem~\ref{TestlosspMLE} provides a deviation bound for the value \( \| \QP (\tilde{\fv}_{\GP} - \fvs) \| \). 
Further we apply 
Corollary~\ref{CBvMNI} to bound the quantity 
\( \| \QP (\fv_{\GP} - \tilde{\fv}_{\GP}) \| \) given \( \Yv \).
%
We conclude by the following result.

\begin{theorem}
\label{Ccontactionrate}
It holds on \( \Omega_{1}(\xx) \) for some fixed 
\( \CONST \)
\begin{EQA}
	\P\Bigl( 
		\bigl\| \QP (\fv_{\GP} - \fvs) \bigr\| 
		\geq 
		\CONST \bigl\| \QP \, \IFa_{\GP}^{-1} \QP^{\T} \bigr\|^{1/2} \, \| \GP \fvs \| + \zq(\BB_{\QPGP},\xx) + 
		\zq\bigl( \QP \, \IF_{\GP}^{-1} \QP^{\T},\xx \bigr)
		\, \cond \Yv 
	\Bigr)
	& \leq &
	2 \, \ex^{-\xx} .
\label{PDvtGtsCrGY}
\end{EQA}
\end{theorem}

A prior ensuring the bias-variance trade-off leads to the optimal contraction rate which 
corresponds to the optimal penalty choice in penalized maximum likelihood estimation;
see below Section~\ref{SminmaxlinNI} for the details in a linear case.

\subsection{Nonparametric Bayes}
\label{SnonpmBayesNI}
One of the main questions of nonparametric Bayes approach is whether one can use Bayesian credible sets as \emph{frequentist confidence sets}. 
Corollary~\ref{CnonparBvMm} suggests to consider credible sets of the form
\begin{EQA}
	\CA_{\QPGP}(\rr)
	& \eqdef &
	\bigl\{ \fv \colon \bigl\| \QP (\tilde{\fv}_{\GP} - \fv) \bigr\| \leq \rr \bigr\} ,
\label{SAQPrtcQttr}
\end{EQA}
and \( \rr = \rr_{\alp} \) is fixed to ensure 
\begin{EQA}
	\PG\bigl( \bigl\| \QP \, \DPblkt_{\GP}^{-1} \gammav \bigr\| > \rr_{\alp} \bigr)
	& = &
	\alp 
\label{alPGQDtm1gra}
\end{EQA}
with \( \gammav \) standard normal.
Our results allow to reduce this question to reliability of pMLE-based confidence sets.

\begin{theorem}
\label{CThonestCS}
Let \( \BB_{\QPGP} = \sigma^{-2} \QP \, \AFblk_{\GP} \, \epsVar^{2} \, \AFblk_{\GP}^{\T} \, \QP^{\T} \) 
and let \( \epsv \) be asymptotically normal in the sense that 
\begin{EQA}
	\sup_{\zq > 0} \Bigl| 
		\P\bigl( \sigma^{-1} \bigl\| \QP \, \AFblk_{\GP} \epsv \bigr\| \leq \zq \bigr)
		- \P\bigl( \bigl\| \BB_{\QPGP}^{1/2} \gammav \bigr\| \leq \zq \bigr)
	\Bigr|
	&=&
	o(1)
\label{szQDGm2nzzVQG}
\end{EQA}
with \( \gammav \in \R^{\dimp} \) standard normal.
Assume also the ``small modeling bias'' condition 
\begin{EQA}
	\frac{\bigl\| \QP \, \IFa_{\GP}^{-1} \QP^{\T} \bigr\| \, \, \| \GP \fvs \|^{2}}{\tr\bigl( \BB_{\QPGP} \bigr)}
	&=&
	o(1) .
\label{QDGm2G2ts2V2}
\end{EQA} 
Then it holds
\begin{EQA}
	\sup_{\zq > 0} \Bigl| 
		\P\bigl( \bigl\| \QP (\tilde{\fv}_{\GP} - \fvs) \bigr\| \leq \zq \bigr)
		- \P\bigl( \bigl\| \BB_{\QPGP}^{1/2} \gammav \bigr\| \leq \zq \bigr)
	\Bigr|
	&=&
	o(1) .
\label{01zPPQVzzts}
\end{EQA}
Moreover, if \( \BB_{\QPGP} \leq \QP \, \DPblk_{\GP}^{-2} \, \QP^{\T} \), then
on the set \( \Omega(\xx) \), it holds with \( \rr_{\alp} \) from \eqref{alPGQDtm1gra}
\begin{EQA}
	\P\bigl( \fvs \in \CA_{\QPGP}(\rr_{\alp}) \bigr)
	& \leq &
	\alp + o(1) . 
\label{Palpo1AGra}
\end{EQA}
\end{theorem}

The condition \( \BB_{\QPGP} \leq \QP \, \DPblk_{\GP}^{-2} \QP^{\T} \) means that 
the variance of pMLE \( \tilde{\fv}_{\GP} \) 
is not larger than the variance \( \DPblkt_{\GP}^{-2} \approx \DPblk_{\GP}^{-2} \) of the posterior
\( \fv_{\GP} \cond \Yv \).
Under correct noise specification \( \epsVar = \Id_{\dimq} \),
this condition is always fulfilled.
Indeed, by Theorem~\ref{TFiWititG}  
\begin{EQA}
	\Var \bigl( \tilde{\upsilonv}_{\GP} \bigr) 
	& \approx &
	\Var\bigl( \DF_{\GP}^{-2} \nabla \zeta \bigr) 
	=
	\sigma^{-2} \DF_{\GP}^{-2} \block\bigl( 0,\Id_{\dimp} \bigr) \DF_{\GP}^{-2} .
\label{VartvGVarDG2}
\end{EQA}
For the posterior covariance, Theorem~\ref{TnonparBvMm} implies
\( \Var\bigl( \vupsilonv_{\GP} \cond \Yv \bigr) \approx \DF_{\GP}^{-2} \) which is 
obviously larger than the one in \eqref{VartvGVarDG2}. 
The same holds for \( \QF \, \tilde{\upsilonv}_{\GP} \) any linear mapping \( \QF \).

We conclude that the ``small bias'' condition \eqref{QDGm2G2ts2V2} 
together with some regularity constraints 
ensures frequentist validity of the credible sets. 

\subsection{Minimax risk for a linear operator \( \KS \)}
\label{SminmaxlinNI}
Consider as a special case a linear Gaussian inverse problem 
\( \Yv = \KS \fvs + \sigma \epsv \) for a homogeneous noise \( \epsv \) with \( \Var(\epsv) = \Id_{\dimp} \),
small noise level \( \sigma \),
and a smooth linear operator \( \KS \).
Of course, there is no need to apply the calming approach in this setup,
one can proceed directly with the log-likelihood 
\( - \| \Yv - \KS \fvs \|^{2}/(2 \sigma^{2}) \). 
We, however, show that calming still applies and does not change 
essentially the results. 
This issue is important because in the general situation, we locally approximate 
the underlying model by a linear one. 

Let \( \GP^{2} \) be fixed. 
We consider the class of ``smooth'' functions defined by 
\begin{EQA}
	\Fsmooth_{\GP}(\GPUB)
	& \eqdef &
	\bigl\{ \fv \colon \| \GP \fv \|^{2} \leq \GPUB\bigr\} 
\label{FGfGf21G}
\end{EQA}
for \( \GPUB > 0 \).
Also consider the class of Gaussian priors \( \ND(0,\JG^{-1} \GP^{-2}) \) for some \( \JG \geq 1 \).
Denote \( \IF = \sigma^{-2} \KS^{\T} \KS \),
\( \GPM^{2} = \JG \GP^{2} \), 
\begin{EQA}
	\IF_{\JG}
	& \eqdef &
	\IF + \GPM^{2}
	=
	\sigma^{-2} \KS^{\T} \KS + \JG \GP^{2}, 
\label{IFGF2G2sm2}
\end{EQA}
and consider the penalized MLE \( \tilde{\fv}_{\JG} \)
\begin{EQA}
	\tilde{\fv}_{\JG}
	&=&
	\argmin_{\fv} \bigl\{ \sigma^{-2} \| \Yv - \KS \fv \|^{2} + \| \GPM \fv \|^{2} \bigr\}
	\\
	&=&
	\bigl( \KS^{\T} \KS + \sigma^{2} \GPM^{2} \bigr)^{-1} \KS^{\T} \Yv
	=
	\sigma^{-2} \IF_{\JG}^{-1} \KS^{\T} \Yv .
\label{tfGsm2IFGm1AT}
\end{EQA}
First we show that 
the penalized MLE \( \tilde{\fv}_{\JG} \) is nearly minimax over the set \( \Fsmooth_{\GP}(\GPUB) \)
for a proper choice of \( \JG \). 
We split the result into lower and upper bounds.

\begin{theorem}
\label{TminimaxupperGP}
Let \( \Yv = \KS \fvs + \sigma \epsv \) with \( \epsv \) obeying the condition \nameref{ED0ref} for 
\( \epsVar^{2} = \Id \).
Given \( \GP^{2} \) and \( \QP \), let us fix a value \( \JG \) such that 
\begin{EQA}
	\GPUB \JG \| \QP \, \IF_{\JG}^{-1} \QP^{\T} \|
	&=&
	\tr \bigl( \QP \, \IF_{\JG}^{-1} \QP^{\T} \bigr) 
\label{JFJm1trIGm1u}
\end{EQA}
with \( \IF_{\JG} \) from \eqref{IFGF2G2sm2}.
Then it holds for the penalized MLE \( \tilde{\fv}_{\JG} \) 
\begin{EQA}
	\sup_{\fvs \in \Fsmooth_{\GP}(\GPUB)}
	\P \biggl( 
		\bigl\| \QP \bigl( \tilde{\fv}_{\JG} - \fvs \bigr) \bigr\|^{2} 
		> 3 \tr \bigl( \QP \, \IF_{\JG}^{-1} \QP^{\T} \bigr)
	\biggr)
	& \leq &
	\ex^{- \GPUB \JG/2} .
\label{JFJm1trIGm1222u}
\end{EQA}
\end{theorem}

\begin{remark}
A value \( \JG \) ensuring \eqref{JFJm1trIGm1u} exists and unique because 
\( \tr(\QP \, \IF_{\JG}^{-1} \QP^{\T}) \) decreases to zero while 
\( \GPUB \JG \| \QP \, \IF_{\JG}^{-1} \QP^{\T} \| \) increases from zero to \( \GPUB \) as 
\( \JG \to \infty \).
The relation \eqref{JFJm1trIGm1u} is a specification of the general bias-variance trade-off from
\eqref{DtsGtsrGa}. 
\end{remark}

Now we present a lower bound for the risk for the special case of \( \QP = \Id \).

\begin{theorem}
\label{TminimaxGP}
Given \( \GP^{2} \), let us fix a value \( \JG \) such that 
\begin{EQA}
	\GPUB \JG \| \IF_{\JG}^{-1} \|
	&=&
	\tr (\IF_{\JG}^{-1}) .
\label{JFJm1trIGm1}
\end{EQA}
Let also
\begin{EQA}
	\tr \bigl( \IF_{\JG}^{-1} / \| \IF_{\JG}^{-1} \| \bigr)
	& \leq &
	\CONST \tr \bigl( \IF_{\JG}^{-2} / \| \IF_{\JG}^{-2} \| \bigr) 
\label{trIFJm1FG222}
\end{EQA} 
with a fixed \( \CONST \).
Then with \( \CONST_{1} > 0 \) and \( \CONST_{2} > 0 \) depending on \( \CONST \) only
\begin{EQA}
	\inf_{\hat{\fv}} \sup_{\fv \in \Fsmooth_{\GP}(\GPUB)} 
	\P \biggl( \bigl\| \hat{\fv} - \fvs \bigr\|^{2} > \CONST_{1} \, \tr(\IF_{\JG}^{-1}) \biggr)
	& \geq &
	1 - \ex^{- \CONST_{2} \GPUB \JG}  .
\label{JFJm1trIGm1222}
\end{EQA}
\end{theorem}

Combining the results of Theorem~\ref{TminimaxupperGP} and \ref{TminimaxGP}
yields that the penalized MLE \( \tilde{\fv}_{\JG} \) with \( \JG \) selected by
\eqref{JFJm1trIGm1} is nearly minimax over the class 
\( \Fsmooth_{\GP}(\GPUB) \).

\subsubsection{Commutative case}
For the purpose of comparing with the existing rate results, 
we now specify the results for the case when \( \KS^{\T} \KS \) and \( \GP^{2} \) commute.
By a proper orthonormal basis transform we can reduce the model to the case when both
\( \KS^{\T} \KS \) and \( \GP^{2} \) are diagonal which leads to a sequence space model; see \cite{KnVaZa2011,KnSzVa2016}.
In addition, we reorder the basis in a way that the eigenvalues \( \gp_{j}^{2} \) of 
\( \GP^{2} \) increase with \( j \).
By \( \ksj_{j}^{2} \) we denote the eigenvalues of \( \KS^{\T} \KS \) for the same basis.
%
To simplify the study we suppose that \( \ksj_{j}^{2} \) decrease as \( j \to \infty \).
This means that smoothness properties of \( \KS \) are coordinated with that of \( \fvs \)
in terms of the eigenvalue decomposition of \( \GP^{2} \).


%
Below we focus on the case when 
\( \sum_{j} \gp_{j}^{2} < \infty \) and \( \gp_{j}^{2} \) increase polynomially with \( j \)
yielding for some \( \CONST \geq 1 \) and each \( J \)
\begin{EQA}
	\CONST^{-1} 
	\leq 
	\frac{1}{J \gp_{J}^{-2}} \sum_{j \geq J} \gp_{j}^{-2}
	& \leq &
	\CONST \,  .
\label{sumjJgjm2C}
\end{EQA}	
A typical example is given by \( \gp_{j}^{2} = j^{2s} \) for \( s > 1/2 \).

\begin{theorem}
\label{TcommutAG}
Suppose that 
	\( \GP^{2} \) and \( \KS^{\T} \KS \) commute;
	the eigenvalues \( \gp_{j}^{2} \) of \( \GP^{2} \) grow polynomially; see \eqref{sumjJgjm2C}, and 
	the eigenvalues \( \ksj_{j}^{2} \) are non-increasing.
Let \( J \) be the smallest index such that
\begin{EQA}
	J \gp_{J}^{2}
	& \geq &
	\GPUB \, \sigma^{-2} \ksj_{J}^{2} .
\label{JgJ2sm2ajJ2}
\end{EQA}
Then with \( \JG = J/\GPUB \), the penalized MLE \( \tilde{\fv}_{\JG} \) is nearly minimax
over \( \Fsmooth_{\GP}(\GPUB) \) and the minimax risk is of order 
\( \tr(\IF_{\JG}^{-1}) \asymp \GPUB \, \gp_{J}^{-2} \).
The prior \( \ND(0,\GP_{\JG}^{-2}) \) leads to the posterior contraction of order 
\( \tr(\IF_{\JG}^{-1}) \asymp \GPUB \, \gp_{J}^{-2} \) as well.
\end{theorem}

\begin{example}
\label{ExminiAGJm2}
Let \( \gp_{j}^{2} = j^{2s} \). 
Let also \( \ksj_{j}^{2} = \KSC j^{-2\alp} \).
The case \( \alp = 0 \) corresponds to the direct problem, while \( \alp > 0 \) 
leads to a linear inverse setup.
The equation \( J \gp_{J}^{2} = \GPUB \, \sigma^{-2} \KSC J^{-2\alp} \) leads to 
\( J = (\GPUB \, \KSC/\sigma^{2})^{1/(2s + 2\alpha +1)} \) and 
\( \tr(\IF_{\JG}^{-1}) \asymp \GPUB \, J^{-s} = \GPUB \, (\GPUB \, \KSC/\sigma^{2})^{-s/(2s + 2\alpha +1)} \).
\end{example}

\subsubsection{Non-commutative case}
If \( \KS^{\T} \KS \) and \( \GP^{2} \) do not commute, the study is a bit more involved,
because we cannot reduce it to the sequence space model. 
However, in many application, even in the non-commutative case, the operator 
\( \KS \) transfers smooth functions into smooth ones. 
This property can be used to extend the result of Theorem~\ref{TcommutAG}.
Again, by a basis change, one can reduce the study to the case when 
the operator \( \GP^{2} \) is diagonal with 
increasing eigenvalues \( \gp_{j}^{2} \).
Denote by \( \II_{j} \) the subspace spanned by the first \( j \) eigenvectors of \( \GP^{2} \).
Let also \( \ProjG_{j} \) be the projector on \( \II_{j} \), and \( \ProjG_{j^{c}} \) be 
the orthogonal projector.
Given \( \JG \), consider \( \IF_{\JG} = \sigma^{-2} \KS^{\T} \KS + \JG \GP^{2} \).
If \( \GP^{2} > 0 \), then this matrix is well posed.
For each fixed \( j \), consider the diagonal blocks 
\( \IF_{\JG,j} = \ProjG_{j} \IF_{\JG} \ProjG_{j}^{\T} \) and 
\( \IF_{\JG,j^{c}} = \ProjG_{j^{c}} \IF_{\JG} \ProjG_{j^{c}}^{\T} \)
corresponding to the subspaces \( \II_{\JG} \) and \( \II_{\JG}^{c} \).
In addition, we impose a regularity condition that for \( \JG,J \) large enough,
it holds
\begin{EQA}
	\CONST^{-1} \, \block\bigl( \IF_{\JG,J}^{-1}, \IF_{\JG,J^{c}}^{-1} \bigr)
	& \leq &
	\IF_{\JG}^{-1}
	\leq 
	\CONST \, \block\bigl( \IF_{\JG,J}^{-1}, \IF_{\JG,J^{c}}^{-1} \bigr) .
\label{Cm1FJJm1cc1}
\end{EQA}
To get a link between smoothness of \( \fvs \) and smoothness of \( \KS \), we define 
for \( \KS_{j} = \ProjG_{j} \KS \)
\begin{EQA}
	\ksj_{j}^{-2}
	& \eqdef &
	\bigl\| \bigl( \KS_{j}^{\T} \KS_{j} \bigr)^{-1} \bigr\| .
\label{kjm2AjTAjm1}
\end{EQA}
Obviously, \( \ksj_{j}^{2} \) decrease with \( j \).
Now we proceed as in the commutative case defining the index \( J \) by the relation 
\( \GPUB \, \sigma^{-2} \ksj_{J}^{2} = J \gp_{J}^{2} \) and letting \( \JG = J/\GPUB \). 

\begin{theorem}
\label{TnoncommutAG}
Suppose that eigenvalues \( \gp_{j}^{2} \) of \( \GP^{2} \) grow polynomially; see \eqref{sumjJgjm2C}. 
Define \( J \) as the smallest index such that \( J \gp_{J}^{2} \geq \GPUB \, \sigma^{-2} \ksj_{J}^{2} \)
with \( \ksj_{j}^{2} \) from \eqref{kjm2AjTAjm1}. 
Set \( \JG = J /\GPUB \) and suppose that \eqref{Cm1FJJm1cc1} holds.
Then the penalized MLE \( \tilde{\fv}_{\JG} \) is nearly minimax
over \( \Fsmooth_{\GP}(\GPUB) \) and the minimax risk is of order 
\( \tr(\IF_{\JG}^{-1}) \asymp \GPUB \, \gp_{J}^{-2} \).
The prior \( \ND(0,\GP_{\JG}^{-2}) \) leads to the posterior contraction of order 
\( \tr(\IF_{\JG}^{-1}) \asymp \GPUB \, \gp_{J}^{-2} \) as well.
\end{theorem}


%

\subsubsection{Calming in linear models}
Now we repeat the calculus for the calming approach with \( \lambda = \sigma^{-2} \). 
Denote \( \gvs = \KS \fvs \).
Let \( \GP^{2} \) be fixed and \( \GP_{\JG}^{2} = \JG \GP^{2} \).
Define also \( \GPY^{2} \) by the equation 
\( \KS^{\T} \GPY^{2} \KS = \GP^{2} \) yielding 
\( \| \GPY \gvs \|^{2} = \| \GP \fvs \|^{2} \).
Set \( \GPY_{\JG}^{2} = \JG \GPY^{2} \).
It is straightforward to check that the condition \nameref{GGref} is automatically fulfilled
with \( \CONST_{\GP|\GPY} = 1 \).

Consider the full \( \upsilonv \)-model
\begin{EQA}
	\LL_{\JG}(\upsilonv)
	&=&
	- \frac{1}{\sigma^{2}} \| \Yv - \gv \|^{2} - \frac{\lambda}{2} \| \gv - \KS \fv \|^{2}
	- \frac{1}{2} \| \GP_{\JG} \fv \|^{2} - \frac{1}{2} \| \GPY_{\JG} \, \gv \|^{2} .
\label{?}
\end{EQA}
The corresponding total Hessian is defined as  
\( \IFT_{\JG} = - \nabla^{2} \LL_{\JG}(\upsilonv) \).
By Lemma~\ref{LeffdimNLI} it holds
\( \IFT_{\JG}^{-1} \asymp \block\bigl( \IF_{\JG}^{-1},\IG_{\JG}^{-1} \bigr) \).
If \( \QF \upsilonv = \QP \fv \), then 
\( \tr \bigl( \QF \IFT_{\JG}^{-1} \QF^{\T} \bigr) 
\asymp \tr\bigl( \QP \IF_{\JG}^{-1} \QP^{\T} \bigr) \).
Now the results of Theorem~\ref{TestlosspMLE} through Theorem~\ref{CThonestCS} 
state the same risk bounds as in the linear case in terms of 
\( \tr\bigl( \QP \IF_{\JG}^{-1} \QP^{\T} \bigr) \) and 
\( \| \QP \IF_{\JG}^{-1} \QP^{\T} \| \).

\section{Proofs of the main results}
\label{SproofsNI}
This section collects the proofs of our results.


\subsection{Proof of Theorem~\ref{TtifvsGPeff}}
The idea of the proof is to show that for each \( \uv \) with \( \| \DF_{\GP} \uv \| = \rr_{\GP} \), 
the derivative of the function \( \LL_{\GP}(\upsilonvs_{\GP} + t \uv) \) in \( t \) is negative for 
\( |t| \geq 1 \). 
This yields that the point of maximum of \( \LL_{\GP}(\upsilonv) \) cannot be outside of 
\( \CA_{\GP}(\rr_{\GP}) \).
Let us fix any \( \uv \) with \( \| \DF_{\GP} \uv \| \leq \rr \).
We use the decomposition
\begin{EQA}
	\LL_{\GP}(\upsilonvs_{\GP} + t \uv) - \LL_{\GP}(\upsilonvs_{\GP})
	&=&
	\bigl\langle \nabla \zeta, \uv \bigr\rangle \, t + \E \LL_{\GP}(\upsilonvs_{\GP} + t \uv) - \E \LL_{\GP}(\upsilonvs_{\GP}) .
\label{LGtsGtuLGts}
\end{EQA}
With \( f(t) = \E \LL_{\GP}(\upsilonvs_{\GP} + t \uv) \), it holds
\begin{EQA}
	\frac{d}{dt} \LL_{\GP}(\upsilonvs_{\GP} + t \uv)
	&=&
	\bigl\langle \nabla \zeta, \uv \bigr\rangle + f'(t) .
\label{frddtLtGstu}
\end{EQA}
The bound \eqref{uTzDGvTxm12z} implies on \( \Omega(\xx) \)
in view of \( \nabla \zeta = (0, \sigma^{-1} \epsv) \) and 
\( \| \DF_{\GP}^{-1} \nabla \zeta \|^{2} = \sigma^{-2} \epsv^{\T} \IGblk_{\GP}^{-1} \epsv \)
\begin{EQA}
	\bigl| \bigl\langle \nabla \zeta, \uv \bigr\rangle \bigr|
	&=&
	\bigl| \langle \DF_{\GP}^{-1} \nabla \zeta, \DF_{\GP} \uv \rangle \bigr|
	\leq 
	\rr \, \zq(\BB_{\VPGP},\xx) .
\label{uTnzrzBGx}
\end{EQA}
By definition of \( \upsilonvs_{\GP} \), it also holds \( f'(0) = 0 \).
Condition \nameref{LL0ref} implies 
\begin{EQA}
	\bigl| f'(t) - t f''(0) \bigr|
	&=&
	\bigl| f'(t) - f'(0) - t f''(0) \bigr|
	\leq 
	3 t^{2} \delta_{3,\GP}(\rr_{\GP}) .
\label{fptfp0fpttfpp13}
\end{EQA}
For \( t = 1 \), we obtain
\begin{EQA}
	f'(1) 
	& \leq &
	f''(0) + 3 \delta_{3,\GP}(\rr_{\GP})
	=
	- \langle \DF_{\GP}^{2} \uv, \uv \rangle + 3 \delta_{3,\GP}(\rr_{\GP})
	=
	- \rr_{\GP}^{2} + 3 \delta_{3,\GP}(\rr_{\GP}) .
\label{fp1fpp13d3rG}
\end{EQA}
If \( 3 \delta_{3,\GP}(\rr_{\GP}) \leq \rho \rr_{\GP}^{2} \) for \( \rho < 1 \), 
then \( f'(1) < 0 \).
Concavity of \( f(t) \) and \( f'(0) = 0 \) imply that \( f'(t) \) decreases in 
\( t \) for \( t > 1 \).
Further, on \( \Omega(\xx) \) by \eqref{uTnzrzBGx}
\begin{EQA}
	\frac{d}{dt} \LL_{\GP}(\upsilonvs_{\GP} + t \uv) \cond_{t=1}
	& \leq &
	\bigl\langle \nabla \zeta, \uv \bigr\rangle - \rr_{\GP}^{2} + 3 \delta_{3,\GP}(\rr_{\GP})
	\\
	& \leq &
	\rr_{\GP} \, \zq(\BB_{\GP},\xx) - \rr_{\GP}^{2} + 3 \delta_{3,\GP}(\rr_{\GP})
	\leq 
	\rr_{\GP} \, \zq(\BB_{\GP},\xx) - (1 - \rho) \rr_{\GP}^{2}
	< 0
\label{ddtLGtstu33}
\end{EQA}
for \( \rr_{\GP} > (1 - \rho)^{-1} \zq(\BB_{\GP},\xx) \).
As \( \frac{d}{dt} \LL_{\GP}(\upsilonvs_{\GP} + t \uv) \) decreases with \( t \geq 1 \) together with \( f'(t) \) due to \eqref{frddtLtGstu}, the same applies to all such \( t \).
This implies the assertion.

\subsection{Proof of Theorem~\ref{TFiWititG} and Corollary~\ref{CBvMNIML}}
To show \eqref{3d3Af12DGttt}, we use that \( \tilde{\upsilonv}_{\GP} \in \CA_{\GP}(\rr_{\GP}) \) and
\( \nabla \LL_{\GP}(\tilde{\upsilonv}_{\GP}) = 0 \).
Therefore,
\begin{EQA}
	\LL_{\GP}(\tilde{\upsilonv}_{\GP} + \uv) - \LL_{\GP}(\tilde{\upsilonv}_{\GP})
	&=&
	\LL_{\GP}(\tilde{\upsilonv}_{\GP} + \uv) - \LL_{\GP}(\tilde{\upsilonv}_{\GP})
	 - \langle \nabla \LL_{\GP}(\tilde{\upsilonv}_{\GP}), \uv \rangle .
\label{LGtttGLGtTtGW}
\end{EQA}
Let us fix any \( \upsilonv \in \CA_{\GP}(\rr_{\GP}) \) and \( \uv \) with 
\( \| \DF_{\GP} \uv \| \leq \rr \), and consider
\begin{EQA}
	f(t)
	=
	f(t,\uv)
	& \eqdef &
	\LL_{\GP}(\upsilonv + t \uv) - \LL_{\GP}(\upsilonv) - \langle \nabla \LL_{\GP}(\upsilonv), \uv \rangle \, t .
\label{fuELtuELGW}
\end{EQA}
As the stochastic term of \( \LL(\upsilonv) \) and thus, of \( \LL_{\GP}(\upsilonv) \) is linear in \( \upsilonv \), it cancels in this expression,
and it suffices to consider the deterministic part \( \E \LL_{\GP}(\upsilonv) \). 
Obviously \( f(0) = 0 \), \( f'(0) = 0 \).
Moreover, \( f''(0) = \langle \nabla^{2} \E \LL_{\GP}(\upsilonv) \, \uv, \uv \rangle 
= - \langle \DF_{\GP}^{2}(\upsilonv) \, \uv, \uv \rangle < 0 \).
Taylor expansion of the third order implies 
\begin{EQA}
	\bigl| f(1) - \frac{1}{2} f''(0) \bigr|  
	& \leq &
	\bigl| \delta_{3}(\upsilonvc,\uv) \bigr| \, , 
	\quad
	\upsilonvc \in [\upsilonv,\upsilonv + \uv] .
\label{fu12un2fud2W}
\end{EQA}
In particular, for any \( \upsilonv \in \CA_{\GP}(\rr_{\GP}) \) 
\begin{EQA}
	\Bigl| 
		\E \LL_{\GP}(\upsilonvs_{\GP}) - \E \LL_{\GP}(\upsilonv) - \frac{1}{2} \bigl\| \DF_{\GP} (\upsilonv - \upsilonvs_{\GP}) \bigr\|^{2} 
	\Bigr|
	& \leq &
	\delta_{3,\GP}(\rr_{\GP}) .
\label{d3GrGELGtsG12}
\end{EQA}
We now use that by Theorem~\ref{TtifvsGPeff}, \( \uv = \upsilonvs_{\GP} - \tilde{\upsilonv}_{\GP} \) fulfills 
\( \| \DF_{\GP} \uv \| \leq \rr_{\GP} \) on \( \Omega(\xx) \).
Therefore, for \( \upsilonv \in \CA_{\GP}(\rr_{\GP}) \)
\begin{EQA}
	&& \nquad
	\Bigl| 
		\LL_{\GP}(\upsilonv) - \LL_{\GP}(\tilde{\upsilonv}_{\GP}) 
		- \frac{1}{2} \| \DFt_{\GP} (\upsilonv - \tilde{\upsilonv}_{\GP}) \|^{2}
	\Bigr|
	\\
	&=&
	\Bigl| 
		\LL_{\GP}(\upsilonv) - \LL_{\GP}(\tilde{\upsilonv}_{\GP}) 
		- \langle \nabla \LL_{\GP}(\tilde{\upsilonv}_{\GP}), \upsilonv - \tilde{\upsilonv}_{\GP} \rangle 
		- \frac{1}{2} \| \DFt_{\GP} (\upsilonv - \tilde{\upsilonv}_{\GP}) \|^{2}
	\Bigr|
	\leq 
	\delta_{3,\GP}(\rr_{\GP}) .
\label{23GrG122LGt}
\end{EQA}
The result \eqref{3d3Af12DGttt} follows. 
Further, as \( \tilde{\upsilonv}_{\GP} \in \CA_{\GP}(\rr_{\GP}) \), it holds
\begin{EQA}
	&& \nquad
	\LL_{\GP}(\tilde{\upsilonv}_{\GP}) - \LL_{\GP}(\upsilonvs_{\GP}) - \frac{1}{2} \| \DF_{\GP}^{-1} \nabla \zeta \|^{2}
	=
	\max_{\upsilonv \in \CA_{\GP}(\rr_{\GP})} 
	\Bigl\{ 
		\LL_{\GP}(\upsilonv) - \LL_{\GP}(\upsilonvs_{\GP}) - \frac{1}{2} \| \DF_{\GP}^{-1} \nabla \zeta \|^{2} 
	\Bigr\}
	\\
	&=&
	\max_{\upsilonv \in \CA_{\GP}(\rr_{\GP})} 
	\Bigl\{ \bigl\langle \upsilonv - \upsilonvs_{\GP}, \nabla \zeta \bigr\rangle
	+ \E \LL_{\GP}(\upsilonv) - \E \LL_{\GP}(\upsilonvs_{\GP}) - \frac{1}{2} \| \DF_{\GP}^{-1} \nabla \zeta \|^{2} 
	\Bigr\}
	\\
	& \leq &
	\max_{\upsilonv \in \CA_{\GP}(\rr_{\GP})} 
	\Bigl\{ \bigl\langle \DF_{\GP} (\upsilonv - \upsilonvs_{\GP}), \DF_{\GP}^{-1} \nabla \zeta \bigr\rangle
		- \frac{1}{2} \| \DF_{\GP} (\upsilonv - \upsilonvs_{\GP}) \|^{2} 
		- \frac{1}{2} \| \DF_{\GP}^{-1} \nabla \zeta \|^{2} 
	\Bigr\} + \delta_{3,\GP}(\rr_{\GP})
	\\
	& \leq &
	\max_{\upsilonv \in \CA_{\GP}(\rr_{\GP})} 
	\Bigl\{ 
		- \frac{1}{2} \| \DF_{\GP} (\upsilonv - \upsilonvs_{\GP}) - \DF_{\GP}^{-1} \nabla \zeta \|^{2} 
	\Bigr\} + \delta_{3,\GP}(\rr_{\GP})
	\leq  
	\delta_{3,\GP}(\rr_{\GP})
\label{d3G1212222B} 
\end{EQA}
and similarly 
\( \LL_{\GP}(\tilde{\upsilonv}_{\GP}) - \LL_{\GP}(\upsilonvs_{\GP}) - \frac{1}{2} \| \DF_{\GP}^{-1} \nabla \zeta \|^{2} \geq - \delta_{3,\GP}(\rr_{\GP}) \).
This two-sided bound yields as \eqref{DGttGtsGDGm13rG} as \eqref{3d3Af12DGttG}.

The last statement \eqref{DPGPm1Cd3rG} of the theorem follows directly from Lemma~\ref{LellUVD2w}
with \( \QP = \DF_{\GP} \) and \( f(\upsilonv) = \E \LL_{\GP}(\upsilonv) \).

Now we show \eqref{QfGtfGzBQGxm1}.
By \eqref{DGttGtsGDGm13rG}, it holds on \( \Omega(\xx) \)
\begin{EQA}
	\| \QF (\tilde{\upsilonv}_{\GP} - \upsilonvs_{\GP}) \|
	& \leq &
	\bigl\| \QF (\tilde{\upsilonv}_{\GP} - \upsilonvs_{\GP}) - \QF \, \DF_{\GP}^{-2} \nabla \zeta \bigr\|
	+ \bigl\| \QF \, \DF_{\GP}^{-2} \nabla \zeta \bigr\|
	\\
	& \leq &
	\bigl\| 
		\QF \, \DF_{\GP}^{-1} \bigl\{ \DF_{\GP} \bigl( \tilde{\upsilonv}_{\GP} - \upsilonvs_{\GP} \bigr) 
		- \DF_{\GP}^{-1} \nabla \zeta \bigr\} 
	\bigr\| + \bigl\| \QF \, \DF_{\GP}^{-2} \nabla \zeta \bigr\|
	\\
	& \leq &
	4 \delta_{3,\GP}(\rr_{\GP}) \bigl\| \QF \, \DF_{\GP}^{-1} \bigr\| 
	+ \bigl\| \QF \, \DF_{\GP}^{-2} \nabla \zeta \bigr\| 
\label{DFm1QDm2nz}
\end{EQA}
Now the special structure of \( \QF \) and \( \nabla \zeta \) implies
\( \| \QF \, \DF_{\GP}^{-1} \|^{2} = \| \QF \, \DF_{\GP}^{-2} \QF^{\T} \| 
= \| \QP \, \DPblk_{\GP}^{-2} \QP^{\T} \| \)
and similarly
\( \QF \, \DF_{\GP}^{-2} \nabla \zeta = \sigma^{-1} \QP \, \AFblk_{\GP} \epsv \).
Finally we apply the bound \eqref{PxivbzzBBroB} or \eqref{PxivbzzBBroBinf} of Theorem~\ref{LLbrevelocroB} and note that 
for the matrix \( \AA = \sigma^{-1} \QP \, \AFblk_{\GP} \epsVar \), it holds
\( \zq(\AA \AA^{\T},\xx) = \zq(\AA^{\T} \AA,\xx) \).

\subsection{Proof of Theorem~\ref{PrhoQPBvM}}
Let \( \tilde{\upsilonv}_{\GP} = \argmax_{\upsilonv} \LL_{\GP}(\upsilonv) \) be the penalized MLE of the parameter \( \upsilonv \).
We aim at bounding from above the quantity
\begin{EQA}
    \rho(\rups)
    &=&
    \frac{\int_{\| \DFt \uv \| > \rups} \exp\bigl\{ \LL_{\GP}(\tilde{\upsilonv}_{\GP}+\uv) \bigr\} d \uv}
    {\int_{\| \DFt \uv \| \leq \rups} \exp \bigl\{ \LL_{\GP}(\tilde{\upsilonv}_{\GP}+\uv) \bigr\} d \uv}  
    \, 
\label{rhopiDGP}
\end{EQA}
with 
\( \DFt^{2} = \IFT(\tilde{\upsilonv}_{\GP}) \) for \( \IFT(\upsilonv) = - \nabla^{2} \E \LL(\upsilonv) \) and 
\( \DF(\upsilonv) = \sqrt{\IFT(\upsilonv)} \).

\paragraph{Step 1}

The use of \( \nabla \LL_{\GP}(\tilde{\upsilonv}_{\GP}) = 0 \) allows to represent
\begin{EQA}
    \rho(\rups)
    &=&
    \frac{\int_{\| \DFt \uv \| > \rups} \exp \bigl\{ 
    	  \LL_{\GP}(\tilde{\upsilonv}_{\GP}+\uv) - \LL_{\GP}(\tilde{\upsilonv}_{\GP}) \bigr\} d \uv}
    	 {\int_{\| \DFt \uv \| \leq \rups} \exp \bigl\{ 
	 	  \LL_{\GP}(\tilde{\upsilonv}_{\GP}+\uv) - \LL_{\GP}(\tilde{\upsilonv}_{\GP}) \bigr\} d \uv}
    \\
    &=&
    \frac{\int_{\| \DFt \uv \| > \rups} \exp \bigl\{ 
    		\LL_{\GP}(\tilde{\upsilonv}_{\GP} + \uv) - \LL_{\GP}(\tilde{\upsilonv}_{\GP}) 
			- \bigl\langle \nabla \LL_{\GP}(\tilde{\upsilonv}_{\GP}), \uv \bigr\rangle
		  \bigr\} d \uv}
    {\int_{\| \DFt \uv \| \leq \rups} \exp \bigl\{ 
    		\LL_{\GP}(\tilde{\upsilonv}_{\GP} + \uv) - \LL_{\GP}(\tilde{\upsilonv}_{\GP}) 
    		- \bigl\langle \nabla \LL_{\GP}(\tilde{\upsilonv}_{\GP}), \uv \bigr\rangle 
		   \bigr\} d \uv}
    \, .
\label{rhopiDGPt}
\end{EQA}
Now we study this expression for any possible value \( \upsilonv \) from the concentration set of \( \tilde{\upsilonv}_{\GP} \).
Consider \( f(\upsilonv) = \E \LL_{\GP}(\upsilonv) \).
As the stochastic term of \( \LL(\upsilonv) \) and thus, of \( \LL_{\GP}(\upsilonv) \) is linear in 
\( \upsilonv \), it holds
\begin{EQA}
	\LL_{\GP}(\upsilonv + \uv) - \LL_{\GP}(\upsilonv) - \bigl\langle \nabla \LL_{\GP}(\upsilonv), \uv \bigr\rangle 
	&=&
	f(\upsilonv + \uv) - f(\uv) - \bigl\langle \nabla f(\upsilonv), \uv \bigr\rangle.
\label{fuELtuELG}
\end{EQA}
Therefore, it suffices to bound the ratio
\begin{EQA}
	\rho(\rups,\upsilonv)
	& \eqdef &
    \frac{\int \Ind\bigl( \| \DF(\upsilonv) \uv \| > \rups \bigr) \exp \bigl\{ 
    		f(\upsilonv + \uv) - f(\uv) - \bigl\langle \nabla f(\upsilonv), \uv \bigr\rangle \bigr\} d \uv}
    	 {\int \Ind\bigl( \| \DF(\upsilonv) \uv \| \leq \rups \bigr)  \exp \bigl\{ 
	 		f(\upsilonv + \uv) - f(\uv) - \bigl\langle \nabla f(\upsilonv), \uv \bigr\rangle \bigr\} d \uv}
	\qquad
\label{rhopifUV}
\end{EQA}
uniformly in \( \upsilonv \) from the set 
\( \bigl\{ \upsilonv \colon \bigl\| \DF_{\GP} \bigl( \upsilonv - \upsilonvs_{\GP} \bigr) \bigr\| 
\leq \rr_{\GP} \bigr\} \); see Theorem~\ref{TtifvsGPeff}.

\paragraph{Step 2}
First we present some bounds for the denominator of \( \rho(\upsilonv) \).
Lemma~\ref{Lintfxupp2} yields
\begin{EQA}
	&& 
	\nquad
	\int_{\| \DF(\upsilonv) \uv \| \leq \rups} \exp \bigl\{ 
		f(\upsilonv + \uv) - f(\uv) - \bigl\langle \nabla f(\upsilonv), \uv \bigr\rangle
	\bigr\} \, d\uv
	\\
	& \geq &
	\bigl( 1 - \err(\rups) \bigr)
	\int_{\| \DF(\upsilonv) \uv \| \leq \rups} \exp \Bigl( - \frac{\| \DF_{\GP}(\upsilonv) \uv \|^{2}}{2} \Bigr) \, d\uv ,
	\\
	&& 
	\nquad
	\int_{\| \DF(\upsilonv) \uv \| \leq \rups} \exp \bigl\{ 
		f(\upsilonv + \uv) - f(\uv) - \bigl\langle \nabla f(\upsilonv), \uv \bigr\rangle
	\bigr\} \, d\uv
	\\
	& \leq &  
	\bigl( 1 + \err(\rups) \bigr)
	\int_{\| \DF(\upsilonv) \uv \| \leq \rups} \exp \Bigl( - \frac{\| \DF_{\GP}(\upsilonv) \uv \|^{2}}{2} \Bigr) \, d\uv ,
	\qquad
\label{1eGAAiDu22dT31}
\end{EQA}
where \( \DF_{\GP}^{2}(\upsilonv) = \IFT_{\GP}(\upsilonv) = - \nabla^{2} f(\upsilonv) \) and 
\( \err(\rups) \) is given by \eqref{LmgfquadELGP}.
%
Moreover, after a proper normalization, the integral 
\( \int_{\| \DF(\upsilonv) \uv \| \leq \rups} \exp \Bigl( - {\| \DF_{\GP}(\upsilonv) \uv \|^{2}}/{2} \Bigr) \, d\uv \)
can be viewed as the probability of the Gaussian event.
Namely
\begin{EQA}
	\frac{\det \DF_{\GP}(\upsilonv)}{(2\pi)^{\dimp/2}}
	\int_{\| \DF(\upsilonv) \uv \| \leq \rups} \exp \Bigl( - \frac{\| \DF_{\GP}(\upsilonv) \uv \|^{2}}{2} \Bigr) \, d\uv
	&=&
	\P\bigl( \bigl\| \DF(\upsilonv) \DF_{\GP}^{-1}(\upsilonv) \gammav \bigr\| \leq \rups \bigr)
\label{dtDGt2pip2Dtm1r0}
\end{EQA}
for a standard normal \( \gammav \in \R^{\dimp} \). 
The choice \( \rups \geq \sqrt{\dimA_{\GP}(\upsilonv)} + \sqrt{2\xx} \) yields by Corollary~\ref{Cchi2p}
\begin{EQA}
	\P\bigl( \bigl\| \DF(\upsilonv) \DF_{\GP}^{-1}(\upsilonv) \gammav \bigr\| \leq \rups \bigr) 
	& \geq &
	1 - \ex^{-\xx} .
\label{PDtDGm1tgar0}
\end{EQA}
If the error term \( \err(\rups) \) is small, we obtain a sharp 
bound for the integral in the denominator of \( \rho(\rups,\upsilonv) \) from 
\eqref{rhopifUV}.

\paragraph{Step 3}
Now we bound the integral on the exterior of \( \UVd = \bigl\{ \uv \colon \| \DF(\upsilonv) \uv \| \leq \rups \bigr\} \).
Linearity of stochastic term in 
\( \LL_{\GP}(\upsilonv) = \LL(\upsilonv) - \| \GPT \upsilonv \|^{2}/2 \) 
and quadraticity of the penalty term imply
\begin{EQA}
	\LL_{\GP}(\upsilonv + \uv) - \LL_{\GP}(\upsilonv) - \bigl\langle \nabla \LL_{\GP}(\upsilonv), \uv \bigr\rangle
	&=&
	\E \LL(\upsilonv + \uv) - \E \LL(\upsilonv) - \bigl\langle \nabla \E \LL(\upsilonv), \uv \bigr\rangle
	- \frac{1}{2} \| \GP \uv \|^{2} \, .
\label{t22Gu2uTvnLGt}
\end{EQA}
Now we apply Lemma~\ref{CELLDtuvLt} with \( f(\upsilonv + \uv) = \E \LL(\upsilonv + \uv) \). 
This function is concave and it holds 
\( - \bigl\langle \nabla^{2} f(\upsilonv) \uv, \uv \bigr\rangle = \| \DF(\upsilonv) \uv \|^{2} \).
The bound \eqref{LGtuLGtr221622} yields for any \( \uv \) with \( \| \DF(\upsilonv) \uv \| = \rr > \rups \)
\begin{EQA}
	&& \nquad
	\LL_{\GP}(\upsilonv + \uv) - \LL_{\GP}(\upsilonv) - \bigl\langle \nabla \LL_{\GP}(\upsilonv), \uv \bigr\rangle 
	=
	f(\upsilonv + \uv) - f(\upsilonv) - \bigl\langle \nabla f(\upsilonv), \uv \bigr\rangle - \| \GP \uv \|^{2}/2
	\\
	& \leq &
	- \CONSTru (\| \DF(\upsilonv) \uv \| \rups - \rups^{2}/2) - \| \GP \uv \|^{2}/2
	\\
	&=&
	- \CONSTru (\| \DF(\upsilonv) \uv \| \rups - \rups^{2}/2) - \| \DF_{\GP}(\upsilonv) \uv \|^{2}/2 
	+ \| \DF(\upsilonv) \uv \|^{2}/2 .
\label{LGttLGT22r22}
\end{EQA}
with \( \CONSTru = 1 - 3 \rups^{-2} \delta_{3}(\rups) \geq 1/2 \) and
\( \DF_{\GP}^{2}(\upsilonv) = \DF^{2}(\upsilonv) + \GP^{2} \).

Now we can use the result about Gaussian integrals from Section~\ref{SGaussintegr}.
With \( \TAU = \DF(\upsilonv) \DF_{\GP}^{-1}(\upsilonv) \), it holds by Lemma~\ref{TGaussintext}
\begin{EQA}
	&& \nquad
	\frac{\det \DF_{\GP}(\upsilonv)}{(2\pi)^{\dimp/2}}
	\int \Ind\bigl( {\| \DF(\upsilonv) \uv \| > \rups} \bigr) \exp \bigl\{ 
		\LL_{\GP}(\upsilonv + \uv) - \LL_{\GP}(\upsilonv) - \bigl\langle \nabla \LL_{\GP}(\upsilonv), \uv \bigr\rangle 
	\bigr\} \, d\uv
	\\
	& \leq &
	     \E \Bigl\{ \exp\Bigl(
        - \CONSTru \rups \| \TAU \gaussv \| + \frac{\CONSTru \rups^{2}}{2}
        + \frac{1}{2} \| \TAU \gaussv \|^{2} 
      \Bigr) \Ind\bigl( \| \TAU \gaussv \| > \rups \bigr) \Bigr\}
    \leq 
    \CONST \ex^{ - (\dimA_{\GP}(\upsilonv) + \xx)/2} .
\label{dDG2pp2LttG2x2}
\end{EQA}
Putting together the resuluts of Step 1 through Step 3 yields the statement
about \( \rho(\rups) \).

\subsection{Proof of Theorem~\ref{TnonparBvMm} and Corollary~\ref{CnonparBvMm}} 
We proceed similarly to the proof of Theorem~\ref{PrhoQPBvM}.
Fix any centrally symmetric set \( A \).
First we restrict the posterior probability to the set 
\( \CAt(\rups) = \{ \uv \colon \| \DFt \uv \| \leq \rups \} \).
Then we apply the quadratic approximation of the log-likelihood function \( \LL(\upsilonv) \).
Denote \( A(\rups) = A \cap \CAt(\rups) \).
Obviously, \( A(\rups) \) is centrally symmetric as well. 
Further, 
\begin{EQA}
	\P\bigl( \vupsilonv_{\GP} - \tilde{\upsilonv}_{\GP} \in A \cond \Yv \bigr)
	&=& 
	\frac{\int_{A} \exp\bigl\{ \LL_{\GP}(\tilde{\upsilonv}_{\GP} + \uv) \bigr\} d\uv}
		 {\int_{\R^{\dimp}} \exp\bigl\{ \LL_{\GP}( \tilde{\upsilonv}_{\GP} + \uv) \bigr\}d\uv}
	\\
	& \leq &
	\frac{\int_{A(\rups)} 
	\exp\bigl\{ \LL_{\GP}( \tilde{\upsilonv}_{\GP} + \uv) - \LL_{\GP}(\tilde{\upsilonv}_{\GP}) 
				- \bigl\langle \nabla \LL_{\GP}(\tilde{\upsilonv}_{\GP}), \uv \bigr\rangle \bigr\} d\uv}
		{\int_{\| \DFt \uv \| \leq \rups} 
	\exp\bigl\{ \LL_{\GP}( \tilde{\upsilonv}_{\GP} + \uv) - \LL_{\GP}(\tilde{\upsilonv}_{\GP}) 
				- \bigl\langle \nabla \LL_{\GP}(\tilde{\upsilonv}_{\GP}), \uv \bigr\rangle \bigr\}d\uv}
	+ \rho(\rups) .
\label{rorLGttGuuG}
\end{EQA}
Now we apply the bounds from the proof of Theorem~\ref{PrhoQPBvM} yielding
\begin{EQA}
	\P\bigl( \vupsilonv_{\GP} - \tilde{\upsilonv}_{\GP} \in A \cond \Yv \bigr)
	& \leq & 
	\frac{\bigl\{ 1 + \err(\rups) \bigr\} 
		  	\int_{A(\rups)} \exp\bigl\{ - \| \DFt_{\GP} \uv \|^{2}/2 \bigr\} d\uv}
		 {\bigl\{ 1 - \err(\rups) \bigr\}
			\int_{\| \DFt \uv \| \leq \rups} \exp\bigl\{ - \| \DFt_{\GP} \uv \|^{2}/2 \bigr\} \, d\uv}	
	+ \rho(\rups) 
	\\
	& \leq &
	\frac{\bigl\{ 1 + \err(\rups) \bigr\} 
			\P\bigl( \DFt_{\GP}^{-1} \gammav \in A \bigr)}
		 {\bigl\{ 1 - \err(\rups) \bigr\}
		 \P\bigl( \| \DFt \DFt_{\GP}^{-1} \gammav \| \leq \rups \bigr)}	
	+ \rho(\rups) .
\label{1mdi1pdimDGu22}
\end{EQA}
This implies the upper estimate for the posterior probability.
Now we prove the lower bound.
It holds in a similar way that 
\begin{EQA}
	\P\bigl( \vupsilonv_{\GP} - \tilde{\upsilonv}_{\GP} \in A \cond \Yv \bigr)
	&=& 
	\frac{\int_{A} \exp\bigl\{ \LL_{\GP}(\tilde{\upsilonv}_{\GP} + \uv) \bigr\} d\uv}
		 {\int_{\R^{\dimp}} \exp\bigl\{ \LL_{\GP}( \tilde{\upsilonv}_{\GP} + \uv) \bigr\}d\uv}
	\\
	& \geq &
	\frac{\int_{A(\rups)} 
			\exp\bigl\{ \LL_{\GP}( \tilde{\upsilonv}_{\GP} + \uv) - \LL_{\GP}(\tilde{\upsilonv}_{\GP}) 
				- \bigl\langle \nabla \LL_{\GP}(\tilde{\upsilonv}_{\GP}), \uv \bigr\rangle \bigr\} d\uv}
		 {\bigl( \int_{\| \DFt \uv \| \leq \rups} + \int_{\| \DFt \uv \| > \rups} \bigr)
	\exp\bigl\{ \LL_{\GP}( \tilde{\upsilonv}_{\GP} + \uv) - \LL_{\GP}(\tilde{\upsilonv}_{\GP}) 
				- \bigl\langle \nabla \LL_{\GP}(\tilde{\upsilonv}_{\GP}), \uv \bigr\rangle \bigr\}d\uv}
	\\
	& \geq & 
	\frac{\bigl\{ 1 - \err(\rups) \bigr\} 
			\P\bigl( \DFt_{\GP}^{-1} \gammav \in A(\rups) \bigr)}
		 {\bigl\{ 1 + \err(\rups) \bigr\} \P\bigl( \| \DFt \DFt_{\GP}^{-1} \gammav \| \leq \rups \bigr) +
			\CONST \ex^{ - (\dimt_{\GP} + \xx)/2}
		 } 
	\\
	& \geq & 
	\frac{\bigl\{ 1 - \err(\rups) \bigr\} 
			\bigl\{ \P\bigl( \DFt_{\GP}^{-1} \gammav \in A \bigr) - \rho(\rups) \bigr\}}
		 {\bigl\{ 1 + \err(\rups) \bigr\} \P\bigl( \| \DFt \DFt_{\GP}^{-1} \gammav \| \leq \rups \bigr) +
			\CONST \ex^{ - (\dimt_{\GP} + \xx)/2}
		 } 	\, .
\label{rorLGttGuuGl}
\end{EQA}
For the case of an arbitrary possibly non-symmetric \( A \), the proof is similar with the use of 
\eqref{d324d4efppm3} instead of \eqref{4d324d4efppm}.
%

\subsection{Proof of Theorem~\ref{TbiasGP} and Lemma~\ref{LIFTaIFa}} 
The definition of \( \upsilonvs \) and \( \upsilonvs_{\GP} \) implies that
\begin{EQA}
	\E \LL_{\GP}(\upsilonvs_{\GP})
	\geq 
	\E \LL_{\GP}(\upsilonvs),
	&\quad &
	\E \LL(\upsilonvs_{\GP})
	\leq 
	\E \LL(\upsilonvs) .
\label{ELGttsGELELs}
\end{EQA}
As \( \E \LL_{\GP}(\upsilonv) = \E \LL(\upsilonv) - \| \GPT \upsilonv \|^{2}/2 \), it follows that
\begin{EQA}
	\E \LL_{\GP}(\upsilonvs_{\GP}) - \E \LL_{\GP}(\upsilonvs)
	& \leq & 
	\frac{1}{2} \bigl\| \GPT \upsilonvs \bigr\|^{2}
	- \frac{1}{2} \bigl\| \GPT \upsilonvs_{\GP} \bigr\|^{2}
	\leq 
	\frac{1}{2} \bigl\| \GPT \upsilonvs \bigr\|^{2} .
\label{12ELGSmGts22}
\end{EQA}
The bound \eqref{d3GrGELGtsG12} with \( \upsilonv = \upsilonvs \) implies the first statement of \eqref{fr33GrGd3rELG22}. 

Further we show that \( \| \GPT \upsilonvs \| \leq \rrbias/2 \) implies \( \| \DF_{\GP} (\upsilonvs_{\GP} - \upsilonvs) \| \leq  \rrbias \).
Indeed, suppose the opposite inequality.
Define \( \uv = \rrbias \DF_{\GP} (\upsilonvs - \upsilonvs_{\GP}) / \| \DF_{\GP} (\upsilonvs_{\GP} - \upsilonvs) \| \), so that 
\( \| \uv \| = \rrbias \). 
The function \( f(t) = \E \LL_{\GP}(\upsilonvs_{\GP}) - \E \LL_{\GP}(\upsilonvs_{\GP} + t \uv) \) is convex in \( t \) and 
\( \upsilonvs_{\GP} + t \uv \in \Thetad \) for \( |t| \leq 1 \).
Using the approximation \eqref{d3GrGELGtsG12} for \( \upsilonv = \upsilonvs_{\GP} + \uv \) implies 
\begin{EQA}
	\E \LL_{\GP}(\upsilonvs_{\GP}) - \E \LL_{\GP}(\upsilonvs_{\GP} + t \uv)
	& \geq & 
	\frac{\rrbias^{2} - \delta_{3,\GP}(\rrbias)}{2} 
	\geq 
	\frac{\rrbias^{2}}{4}
\label{tts24ELGEL}
\end{EQA}
and concavity of \( \E \LL_{\GP}(\upsilonv) \) together with \( \nabla \E \LL_{\GP}(\upsilonvs_{\GP}) = 0 \) implies  
\begin{EQA}
	\E \LL_{\GP}(\upsilonvs_{\GP}) - \E \LL_{\GP}(\upsilonvs_{\GP} + t \uv)
	& \geq &
	\frac{\rrbias^{2}}{4}
\label{24t1ELGELGtu}
\end{EQA}
for \( t \geq 1 \). 
This contradicts to the bounds \eqref{12ELGSmGts22} and \( \| \GPT \upsilonvs \|^{2} \leq \rrbias^{2}/2 \).

Now for any \( \upsilonv \) with \( \| \DF_{\GP} (\upsilonvs_{\GP} - \upsilonv) \| \leq \rrbias \)
\begin{EQA}
	\Bigl| 
		\E \LL_{\GP}(\upsilonvs_{\GP}) - \E \LL_{\GP}(\upsilonv) - \frac{1}{2} \bigl\| \DF_{\GP} (\upsilonv - \upsilonvs_{\GP}) \bigr\|^{2}
	\Bigr|
	& \leq &
	\delta_{3,\GP}(\rrbias) .
\label{BB3GrGLGD2}
\end{EQA}
Further we use that \( \upsilonvs = \argmax \E \LL(\upsilonv) \) and 
\( \E \LL_{\GP}(\upsilonv) = \E \LL(\upsilonv) - \| \GPT \upsilonv \|^{2}/2  \).
By \eqref{BB3GrGLGD2} in view of \( \| \DF_{\GP} (\upsilonvs_{\GP} - \upsilonvs) \| \leq  \rrbias \)
and \( \DF_{\GP}^{2} = \DF^{2} + \GPT^{2} \)
\begin{EQA}
	\E \LL(\upsilonvs) - \E \LL_{\GP}(\upsilonvs_{\GP})
	&=&
	\max_{\upsilonv \in \CA_{\GP}(\rrbias)} 
	\bigl\{ \E \LL_{\GP}(\upsilonv) + \frac{1}{2} \| \GPT \upsilonv \|^{2} 
		- \E \LL_{\GP}(\upsilonvs_{\GP}) \bigr\}
	\\
	& \leq & 
	\max_{\upsilonv \in \CA_{\GP}(\rrbias)} 
	\Bigl\{ 
		- \frac{1}{2} \bigl\| \DF_{\GP} (\upsilonv - \upsilonvs_{\GP}) \bigr\|^{2} + \frac{1}{2} \| \GPT \upsilonv \|^{2} 
	\Bigr\} + \delta_{3,\GP}(\rrbias)
	\\
	& = & 
	\max_{\upsilonv \in \CA_{\GP}(\rrbias)} 
	\Bigl\{ 
		- \frac{1}{2} \bigl\| \DF \upsilonv - \DF^{-1} \DF_{\GP}^{2} \upsilonvs_{\GP} \bigr\|^{2} 
		+ \frac{1}{2} \| \DF^{-1} \DF_{\GP}^{2} \upsilonvs_{\GP} \|^{2} 
	\Bigr\} + \delta_{3,\GP}(\rrbias)
\label{12d34GDm1DG2}
\end{EQA}
A similar inequality holds from below with another sign for \( \delta_{3,\GP} \)-term 
yielding for the maximizer \( \upsilonvs \) the bound
\begin{EQA}
	\bigl\| \DF \upsilonvs - \DF^{-1} \DF_{\GP}^{2} \upsilonvs_{\GP} \bigr\|^{2}
	& \leq &
	4 \delta_{3,\GP}(\rrbias) .
\label{43GrbDtsm1}
\end{EQA}
Equivalently, using again \( \DF_{\GP}^{2} = \DF^{2} + \GPT^{2} \)
\begin{EQA}
	\bigl\| \DF^{-1} \DF_{\GP}^{2} (\upsilonvs - \upsilonvs_{\GP}) - \DF^{-1} \GPT^{2} \upsilonvs \bigr\|^{2}
	& \leq &
	4 \delta_{3,\GP}(\rrbias) .
\label{Dm1D2ttssGm1}
\end{EQA}
This implies for any linear \( \QF \) by Cauchy-Schwarz inequality
\begin{EQA}
	\bigl\| \QF (\upsilonvs - \upsilonvs_{\GP}) \bigr\|
	& \leq &
	\bigl\| \QF \, \DF_{\GP}^{-2} \GPT^{2} \upsilonvs \bigr\|
	+ 2 \bigl\| \QF \, \DF_{\GP}^{-2} \, \DF \bigr\| \sqrt{\delta_{3,\GP}(\rrbias)}
	\\
	& \leq &
	\bigl\| \QF \, \DF_{\GP}^{-2} \GPT \bigr\| \, \| \GPT \upsilonvs \| 
	+ 2 \bigl\| \QF \, \DF_{\GP}^{-2} \, \DF \bigr\| \sqrt{\delta_{3,\GP}(\rrbias)} \, .
\label{2QFDFm2Dd3Grb}
\end{EQA}
It remains to note that \( \GPT^{2} \leq \DF_{\GP}^{2} \) and 
thus \( \bigl\| \QF \, \DF_{\GP}^{-2} \GPT \bigr\|^{2} \leq \bigl\| \QF \, \DF_{\GP}^{-2} \QF^{\T} \bigr\| \) and similarly \( \bigl\| \QF \, \DF_{\GP}^{-2} \DF \bigr\|^{2} \leq \bigl\| \QF \, \DF_{\GP}^{-2} \QF^{\T} \bigr\| \).
If \( \QF \upsilonv = \QP \fv \), then we can use 
\( \QF \, \DF_{\GP}^{-2} \QF^{\T} = \QP \, \DPblk_{\GP}^{-2} \QP^{\T} \).
In addition, by \nameref{GGref} it holds \( \| \GPT \upsilonvs \| \lesssim \| \GP \fvs \| \).

For proving Lemma~\ref{LIFTaIFa} we use , where 
\( \upsilonvsa_{\GP} = (\fvs_{\GP}, \gvsa_{\GP}) \) 
with \( \gvsa_{\GP} = \KS(\fvs_{\GP}) \).
In view of \eqref{12ELGSmGts22}, it holds 
\begin{EQA}
	\sigma^{-2} \| \gvs - \gvs_{\GP} \|^{2} + \lambda \| \gvs_{\GP} - \KS(\fvs_{\GP}) \|^{2} 
	+ \| \GPT \upsilonvs_{\GP} \|^{2} 
	& \leq &
	\| \GPT \upsilonvs \|^{2} 
\label{GT2usgvsAfGssm2}
\end{EQA}
yielding for \( \lambda \geq \sigma^{-2} \)
\begin{EQA}
	\sigma^{-2} \| \gvs_{\GP} - \gvsa_{\GP} \|^{2}
	& \leq &
	\| \GPT \upsilonvs \|^{2} . 
\label{Gs2sm2gsGgsa2}
\end{EQA}
This, similarly to the above, this point belongs to local vicinity of \( \upsilonvs_{\GP} \).
This allows to apply the result of Lemma~\ref{LellUVD2w} for 
\( \xv = \upsilonvs_{\GP} \) and 
\( \uv = \upsilonvsa_{\GP} - \upsilonvs_{\GP} \).

\subsection{Proof of Theorem~\ref{CThonestCS}}
Note first that by definition,
it holds for the true parameter \( \upsilonvs \):
\begin{EQA}
	\P\bigl( \upsilonvs \in \CA_{\QPGP}(\rr) \bigr)
	&=&
	\P\bigl( \bigl\| \QP (\tilde{\upsilonv}_{\GP} - \upsilonvs) \bigr\| \leq \rr \bigr) .
\label{PtsiSAQrQtGtsr}
\end{EQA}
The Fisher expansion \eqref{DGttGtsGDGm13rG} 
\( \tilde{\upsilonv}_{\GP} - \upsilonvs_{\GP} \approx \DFGP^{-2} \nabla \zeta \) 
of Theorem~\ref{TFiWititG} combined with the CLT 
\( \VP^{-1} \nabla_{\gv} \zeta \tow \gammav \) for a standard normal \( \gammav \) reduces the latter question to Gaussian probability
\begin{EQA}
	\P\bigl( \bigl\| \QP (\tilde{\upsilonv}_{\GP} - \upsilonvs) \bigr\| \leq \rr \bigr)
	& \approx &
	\P\Bigl( \bigl\| \QP \bigl(\DFGP^{-2} \VF \gammav + \upsilonvs_{\GP} - \upsilonvs\bigr) \bigr\| \leq \rr \Bigr) .
\label{PPQQDGm2ttsGtts}
\end{EQA}
By Gaussian comparison Theorem~\ref{Tgaussiancomparison3}, 
the impact of the bias \( \upsilonvs_{\GP} - \upsilonvs \)
is negligible under the undersmoothing condition \( \| \QP (\upsilonvs_{\GP} - \upsilonvs) \|^{2} \ll \tr \bigl( \QP \DFGP^{-2} \VF \bigr)^{2} \). 
Combining with Theorem~\ref{TestlosspMLE} yields 
in view of \( \DFGP^{-2} \VF \leq \DFGP^{-1} \)
\begin{EQA}
	1 - \alp
	&=&
	\PG\bigl( \bigl\| \QP \DFt_{\GP}^{-1} \gammav \bigr\| \leq \rr_{\alp} \bigr)
	\approx
	\P\bigl( \bigl\| \QP \DFGP^{-1} \gammav \bigr\| \leq \rr_{\alp} \bigr)
	\leq 
	\P\Bigl( \bigl\| \QP \DFGP^{-2} \VF \gammav \bigr\| \leq \rr_{\alp} \Bigr)
	\\
	& \approx &
	\P\bigl( \bigl\| \QP (\tilde{\upsilonv}_{\GP} - \upsilonvs) \bigr\| \leq \rr_{\alp} \bigr),
\label{1maPGPQm2111}
\end{EQA}
that is, the credible set \( \CA_{\QPGP}(\rr_{\alp}) \) is an asymptotically valid confidence set.

\subsection{Proof of Theorem~\ref{TminimaxupperGP} and \ref{TminimaxGP}}
Let a linear mapping \( \QP \) be fixed, and let \( \JG \) satisfy \eqref{JFJm1trIGm1u}.
It holds with \( \GPM^{2} = \JG \GP^{2} \) 
\begin{EQA}
	\QP \bigl( \tilde{\fv}_{\JG} - \fvs \bigr)
	&=&
	\sigma^{-1} \QP \, \IF_{\JG}^{-1} \KS^{\T} \epsv + 
	\sigma^{-2} \QP \bigl( \IF_{\JG}^{-1} \KS^{\T} \KS - \Id_{\dimp} \bigr) \fvs
	\\
	&=&
	\sigma^{-1} \QP \, \IF_{\JG}^{-1} \KS^{\T} \epsv - \QP \, \IF_{\JG}^{-1} \GPM^{2} \fvs .
\label{sm121QtfGfs}
\end{EQA}
In view of \( \sigma^{-2} \KS^{\T} \KS \leq \IF_{\JG} \)
\begin{EQA}
	\sigma^{-2} \QP \, \IF_{\JG}^{-1} \KS^{\T} \KS \, \IF_{\JG}^{-1} \QP^{\T} 
	& \leq &
	\QP \, \IF_{\JG}^{-1} \, \QP^{\T} .
\label{QPIFGPQPTGP}
\end{EQA}
This yields by Theorem~\ref{TexpbLGA} with \( \BB = \BB_{\JG} = \QP \, \IF_{\JG}^{-1} \, \QP^{\T} \)
\begin{EQA}
	\P \Bigl( 
		\| \sigma^{-1} \QP \, \IF_{\JG}^{-1} \KS^{\T} \epsv \| 
		> \zq\bigl( \BB_{\JG},\xx \bigr) 
	\Bigr)
	& \leq &
	\ex^{-\xx} ,
\label{QPIFGPQPTGPfs}
\end{EQA}
where 
\begin{EQA}
	\zq\bigl( \BB_{\JG},\xx \bigr)
	& \leq &
	\sqrt{\tr \BB_{\JG}} + \sqrt{2 \xx \| \BB_{\JG} \|} .
\label{zBlaxstrGxs2x}
\end{EQA}
Further, by \( \GPM^{2} \leq \IF_{\JG} \) and \( \| \GP \fvs \|^{2} \leq \GPUB\)
\begin{EQA}
	\| \QP \, \IF_{\JG}^{-1} \GPM^{2} \fvs \|^{2}
	& \leq &
	\| \QP \, \IF_{\JG}^{-1} \, \GPM^{2} \, \IF_{\JG}^{-1} \,\QP^{\T} \| \, \| \GPM \fvs \|^{2} 
	\leq 
	\JG \, \| \BB_{\JG} \| \, \GPUB .
\label{QPIFGPQPTGPfs}
\end{EQA}
This yields
\begin{EQA}
	\P \Bigl( 
		\bigl\| \QP \bigl( \tilde{\fv}_{\JG} - \fvs \bigr) \bigr\| 
		> \zq\bigl( \BB_{\JG},\xx \bigr) + \sqrt{\GPUB \JG \, \| \BB_{\JG} \| }
	\Bigr)
	& \leq &
	\ex^{-\xx}  .
\label{QPIFGPQPTGPfs}
\end{EQA}
Suppose that \( \GPUB\JG \, \| \BB_{\JG} \| \leq \tr \BB_{\JG} \).
With \( \xx = \GPUB\JG /2 \) and 
\( \zq(\BB_{\JG},\xx) = \sqrt{\tr \BB_{\JG}} + \sqrt{2 \xx \| \BB_{\JG} \|} \), we obtain
\begin{EQA}
	\P \biggl( 
		\bigl\| \QP \bigl( \tilde{\fv}_{\JG} - \fvs \bigr) \bigr\| > 3 \sqrt{\tr \BB_{\JG}}
	\biggr)
	& \leq &
	\ex^{- \GPUB \JG/2} .
\label{emlam2trBla}
\end{EQA}
%

%
Now we prove the minimax risk bound for \( \QP = \Id \).
We use that the minimax risk is always smaller than the Bayes one whatever a prior is taken.
It is only important to check that the applied prior is concentrated on the set \( \cc{F}(\GPUB) \).

Condition \eqref{trIFJm1FG222} can be rewritten as
\begin{EQA}
	\tr ( \IF_{\JG}^{-1})
	& \leq &
	\CONST \| \IF_{\JG}^{-1} \| \, \tr ( \IF_{\JG}^{-2}) .
\label{trFlam1CFJm1m2}
\end{EQA}
Define a prior \( \fv \sim \ND(0,\GP_{1}^{-2}) \) for
\begin{EQA}
	\GP_{1}^{2}
	&=&
	\| \IF_{\JG}^{-1} \| \, \IF_{\JG}^{2} .
\label{G12IFGm1J2}
\end{EQA}
As \( \fv \cond \Yv \sim \ND(\tilde{\fv}_{\GP_{1}},\IF_{\GP_{1}}^{-1}) \),
it follows by the lower bound \eqref{Pxiv2dimAvp12m} of Theorem~\ref{TexpbLGA} 
with \( \dimA_{1} = \tr(\IF_{\GP_{1}}^{-1}) \) and 
\( \vp_{1}^{2} = \tr(\IF_{\GP_{1}}^{-1})^{2} 
\leq \| \IF_{\GP_{1}}^{-1} \| \tr(\IF_{\GP_{1}}^{-1}) \)
\begin{EQA}
	\P \Bigl( \| \fv - \tilde{\fv}_{\GP_{1}} \|^{2} < \dimA_{1} - 2 \vp_{1} \sqrt{\xx} \, \cond \Yv 
	\Bigr)
	& \leq &
	\ex^{-\xx}.
\label{EEtrIFm1QQT}
\end{EQA}
In view of \eqref{trIFJm1FG222} the value \( \dimA_{1} \) can be bounded from below by
\begin{EQA}
	\dimA_{1}
	=
	\tr\bigl( \IF_{\GP_{1}}^{-1} \bigr)
	&=&
	\tr\bigl\{ \bigl( \IF + \| \IF_{\JG}^{-1} \| \, \IF_{\JG}^{2} \bigr)^{-1} \bigr\}
	\geq 
	\frac{1}{2} \tr\bigl\{ \IF_{\JG}^{-2} / \| \IF_{\JG}^{-1} \| \bigr\}
	\geq 
	\frac{1}{2\CONST} \tr (\IF_{\JG}^{-1}) .
\label{CtrIFJm2Jm2FGm1}
\end{EQA}
Similarly
\begin{EQA}
	\| \IF_{\GP_{1}}^{-1} \|
	&=&
	\bigl\| \bigl( \IF + \| \IF_{\JG}^{-1} \| \, \IF_{\JG}^{2} \bigr)^{-1} \bigr\|
	\leq 
	\bigl\| \bigl( \| \IF_{\JG}^{-1} \| \, \IF_{\JG}^{2} \bigr)^{-1} \bigr\|
	=
	\| \IF_{\JG}^{-1} \| .
\label{GG1m1lGGG}
\end{EQA}
With \( \xx = \GPUB \JG/(9 \CONST) \), we derive
\begin{EQA}
	\dimA_{1} - 2 \vp_{1} \sqrt{\xx}
	& \geq &
	\sqrt{\dimA_{1}} \biggl( \sqrt{\dimA_{1}} - 2 \sqrt{\JG \| \IF_{\JG}^{-1} \|/(9\CONST)} \biggr)
	\geq 
	\CONST_{1} \tr (\IF_{\JG}^{-1}) 
\label{p1C1trFGm1m2v1}
\end{EQA}
and hence
\begin{EQA}
	\P \Bigl( \| \fv - \tilde{\fv}_{\GP_{1}} \|^{2} < \CONST_{1} \tr (\IF_{\JG}^{-1}) 
	\Bigr)
	& \leq &
	\ex^{- \CONST_{2} \GPUB \JG},
\label{EEtrIFm1QQT9}
\end{EQA}
where \( \CONST_{1} \) and \( \CONST_{2} \) depend on \( \CONST \) only.
It remains to show that this prior concentrates on the set of ``smooth'' functions with 
\( \| \GP \fv \| \leq 2 \sqrt{\GPUB} \).
Indeed, \( \fv = \GP_{1}^{-1} \gammav \) for \( \gammav \) standard normal, and
one can apply the deviation bound for Gaussian quadratic forms from Theorem~\ref{TexpbLGA} with
\( \BB = \BB_{1} = \bigl( \GP \GP_{1}^{-1} \bigr)^{\T} \GP \GP_{1}^{-1} \):
\begin{EQA}
	\P\bigl( \| \GP \fv \| > \zq(\BB_{1},\xx) \bigr)
	&=&
	\P\bigl( \bigl\| \GP \GP_{1}^{-1} \gammav \bigr\| > \zq(\BB_{1},\xx) \bigr)
	\leq 
	\ex^{-\xx} .
\label{exPm1Gm2Gf2}
\end{EQA}
Further, it holds in view of \( \JG \GP^{2} \leq \IF_{\JG} \) and \eqref{JFJm1trIGm1}
\begin{EQA}
	\tr \bigl( \GP_{1}^{-1} \GP^{2} \GP_{1}^{-1} \bigr)
	& \leq &
	\tr \bigl( \JG^{-1} \IF_{\JG} \IF_{\JG}^{-2} / \| \IF_{\JG}^{-1} \| \bigr)
	=
	\JG^{-1} \tr(\IF_{\JG}^{-1}) / \| \IF_{\JG}^{-1} \|
	=
	\GPUB
\label{JtrGFm1m12m2JG}
\end{EQA}
and similarly
\begin{EQA}
	\| \GP_{1}^{-1} \GP^{2} \GP_{1}^{-1} \|
	& \leq & 
	\JG^{-1} \| \IF_{\JG}^{-1} \| / \| \IF_{\JG}^{-1} \|
	=
	\JG^{-1} .
\label{exPm1Gm2Gf2GG1}
\end{EQA}
We derive with \( \xx = \GPUB \JG/2 \)
\begin{EQA}
	\zq(\BB_{1},\xx)
	& \leq &
	\sqrt{\tr \BB_{1}} + \sqrt{2 \xx \| \BB_{1} \|}
	\leq 
	2 \sqrt{\GPUB} .
\label{zqB1x2smutrB1}
\end{EQA}

\subsection{Proof of Theorem~\ref{TcommutAG} and Theorem~\ref{TnoncommutAG}}
For simplicity assume that there exists \( J \) with \( \sigma^{-2} \ksj_{J}^{2} = \JG \gp_{J}^{2} \).
Obviously 
\( \IF_{\JG} = \diag\bigl\{ (\sigma^{-2} \ksj_{j}^{2} + \JG \gp_{j}^{2}) \bigr\} \) and
\( \| \IF_{\JG}^{-1} \| = 1/ \bigl( 2 \JG \gp_{J}^{2} \bigr) \).
As \( \ksj_{j}^{-2} \) increase and \( \gp_{j} \) grow polynomially with \( j \), we bound
\begin{EQA}
	\tr(\IF_{\JG}^{-1})
	&=&
	\sum_{j} \frac{1}{\sigma^{-2} \ksj_{j}^{2} + \JG \gp_{j}^{2}}
	\leq 
	\sum_{j < J} \frac{1}{\sigma^{-2} \ksj_{j}^{2}}
	+ \sum_{j \geq J} \frac{1}{\JG \gp_{j}^{2}} 
	\leq 
	\frac{J}{\sigma^{-2} \ksj_{J}^{2}}
	+ \frac{\CONST \, J}{\JG \gp_{J}^{2}} 
	=
	\frac{(\CONST + 1) J}{\JG \gp_{J}^{2}} \, .
\label{trFlam1sj1s2a2}
\end{EQA}
Similarly 
\begin{EQA}
	\tr(\IF_{\JG}^{-1})
	& \geq &
	\sum_{j \geq J} \frac{1}{\JG \gp_{j}^{2}} 
	=
	\frac{\CONST^{-1} J}{\JG \gp_{J}^{2}} \, .
\label{trFlam1sj1s2a21}
\end{EQA}
Therefore, the relation \eqref{JFJm1trIGm1} corresponds to \( \GPUB\JG \asymp J \)
yielding the minimax risk of order
\( \tr (\IF_{\JG}^{-1}) \asymp \GPUB \, \gp_{J}^{-2} \).

For the non-commutative case, we proceed separately for the blocks of \( \IF_{\JG} \).
It holds in the similar way
\begin{EQA}
	\tr(\IF_{\JG,J^{c}}^{-1})
	& \asymp &
	\sum_{j \geq J} \frac{1}{\JG \gp_{j}^{2}} 
	\asymp
	\frac{J}{\JG \gp_{J}^{2}}
	\asymp
	\frac{\GPUB}{\gp_{J}^{2}}
	\, ,
\label{trFlam1sj1s2a21}
\end{EQA}
and by \eqref{kjm2AjTAjm1}
\begin{EQA}
	\tr (\IF_{\JG,J}^{-1})
	& \leq &
	J \bigl\| \bigl( \sigma^{-2} \KS^{\T} \KS \bigr)^{-1} \bigr\|
	\asymp
	\frac{J}{\JG \gp_{J}^{2}} 
	\asymp
	\frac{\GPUB}{\gp_{J}^{2}} \, .
\label{gjm2JATAtrFJ}
\end{EQA}
This yields the result by \eqref{Cm1FJJm1cc1}.

\appendix
\section{Tools}
Below we present some technical results and useful external references.

\subsection{Concavity and tail bounds}
\label{SGaussintegr}

Let \( f(\xv) \) be a function on \( \R^{\dimp} \).
Previous results describe the local behavior of \( f(\xv + \uv) \) for \( \uv \in \UV \)
under local smoothness conditions. 
Now we derive some upper bounds on \( f(\xv + \uv) \) for \( \uv \) large using that \( f \) is concave.
More precisely, we fix \( \xv \) and \( \uv \) and bound the values
\(f(\xv + t \uv) - f(\xv) - t f'(\xv,\uv) \) for \( \uv \in \UV \) and large \( t \).

\begin{lemma}
\label{LffpfppEL}
Suppose \eqref{1mffmxum34} with \( \delta_{m} \leq 1\) for \( m=3,4 \).
Let \( \xv + \UV \subset \Xs \). 
Let the function \(f(\xv + t \uv)\) be concave in \( t \). 
Then it holds for any \( \uv \in \UV\) and for \( t > 1 \)
\begin{EQA}
    f(\xv + t \uv) - f(\xv) - \bigl\langle \nabla f(\xv), \uv \bigr\rangle \, t
    & \leq &
    \Bigl( t - \frac{1}{2} \Bigr) 
    \Bigl\{ \langle \nabla^{2} f(\xv) \uv, \uv \rangle - 3 \delta_{3} \Bigr\} .
\label{tm12g1g026}
\end{EQA}
\end{lemma}

\begin{proof}
The  Taylor expansion of the third order for \( g(t) = f(\xv + t\uv) \) at \( t=0 \) yields
\begin{EQA}
    \biggl| g(1) - g(0) - g'(0) - \frac{1}{2} g''(0) \biggr|
    & \leq &
    \delta_{3} .
\label{d3Ad3g10}
\end{EQA}
Similarly one obtains
\begin{EQA}
    g'(1) - g'(0)
    &=&
    g'(1) - g'(0) - g''(0) + g''(0)
    \leq 
    g''(0) + 3 \delta_{3} \, .
\label{gpp063ut}
\end{EQA}
Concavity of \(g(\cdot)\) implies
\begin{EQA}
    g(t) - g(1)
    & \leq &
    (t-1) g'(1) .
\label{gtg1tm1gp1}
\end{EQA}
We summarize that
\begin{EQA}
    g(t) - g(0) - t g'(0)
    & = &
    g(t) - g(1) - (t-1) g'(1) + (t-1) \bigl\{ g'(1) - g'(0) \bigr\} + g(1) - g(0) - g'(0)
    \\
    & \leq &
    (t-1) \bigl\{ g''(0) + 3 \delta_{3} \bigr\}
    + \frac{1}{2} g''(0) + \delta_{3}
    \\
    & \leq &
    (t-1/2) \bigl\{ g''(0) + 3 \delta_{3} \bigr\} .
\label{t12gpp6d3uvt}
\end{EQA}
This implies the assertion in view of \(g''(0) = \langle \nabla^{2} f(\xv) \uv, \uv \rangle\).
\end{proof}

Now we specify the result of Lemma~\ref{LffpfppEL} for the elliptic set 
\( \UV(\rups)\) defined by the condition \(- \langle \nabla^{2} f(\xv) \uv, \uv \rangle \leq \rups^{2}\).
We write \( \delta_{3}(\rups)\) in place of 
\( \delta_{3}(\Xs,\UV(\rups))\).
We aim at bounding from above 
the value \(f(\xv + \uv) - f(\xv)\) for
\( \uv\) with \(- \langle \nabla^{2} f(\xv) \uv, \uv \rangle = \rr^{2} > \rups^{2}\).

\begin{lemma}
\label{CELLDtuvLt}
Consider \( \xv \in \Xs\) and 
\( \UV = \UV(\rups) = \bigl\{ \uv \colon - \langle \nabla^{2} f(\xv) \uv, \uv \rangle \leq \rups^{2} \bigr\}\).
Let \(f(\xv + \uv)\) be concave in \( \uv\).
Then for any \( \uv\) with \(- \langle \nabla^{2} f(\xv) \uv, \uv \rangle = \rr^{2} > \rups^{2}\), it holds
\begin{EQA}
  f(\xv + \uv) - f(\xv) - \bigl\langle \nabla f(\xv), \uv \bigr\rangle
  & \leq &
  - (\rr \rups - \rups^{2}/2) 
  \bigl\{ 1 - 3 \rups^{-2} \delta_{3}(\rups) \bigr\} .
\label{LGtuLGtr221622}
\end{EQA}
\end{lemma}

\begin{proof}
Define \(t = \rr/\rups\) and \( \uvd = \uv \rups/\rr\), so that 
\( - \langle \nabla^{2} f(\xv) \uv, \uv \rangle = \rups^{2} \) and \( \uvd \in \UV(\rups)\).
Then it holds by \eqref{tm12g1g026}
\begin{EQA} 
    f(\xv + \uv) - f(\xv) - \bigl\langle \nabla f(\xv), \uv \bigr\rangle
    &=&
    f(\xv + t \uvd) - f(\xv) - \bigl\langle \nabla f(\xv), \uv \bigr\rangle
    \\
    & \leq &
    - (\rr/\rups - 1/2) \bigl\{ \rups^{2} - 3 \delta_{3}(\rups) \bigr\}
    \\
    &=&
    - (\rr \rups - \rups^{2}/2) 
    \bigl\{ 1 - 3 \rups^{-2} \delta_{3}(\rups) \bigr\} 
\label{rr0mr0216r2}
\end{EQA}
and the result follows.
\end{proof}

The result is meaningful if \(3 \rups^{-2} \delta_{3}(\rups) < 1\).
Then with \( \CONSTru = 1 - 3 \rups^{-2} \delta_{3}(\rups)\), we obtain
for any \( \uv\) with \(- \langle \nabla^{2} f(\xv) \uv, \uv \rangle = \rr^{2} > \rups^{2}\)
\begin{EQA}
    f(\xv + \uv) - f(\xv) - \bigl\langle \nabla f(\xv), \uv \bigr\rangle
    & \leq &
    - \CONSTru (\rr \rups - \rups^{2}/2) .
\label{LGtuLGuuT22r22}
\end{EQA}

\subsection{Gaussian integrals}
Let \( \TAU \) be a linear operator in \( \R^{\dimp}\), \( \dimp \leq \infty \), 
with \( \| \TAU \|_{\oper} \leq 1\).
By \( \TAU^{\T} \) we denote the adjoint operator for \( \TAU \).
Given positive \( \rups\) and \( \CONSTru\), consider the following ratio 
\begin{EQA}[c]
    \frac{\int_{\| \TAU \uv \| > \rups}
    \exp\bigl(
      - \CONSTru \| \TAU \uv \| + \frac{1}{2} \CONSTru \rups^{2}
      + \frac{1}{2} \| \TAU \uv \|^{2} 
      - \frac{1}{2} \| \uv \|^{2} 
    \bigr) d\uv}
    {\int_{\| \TAU \uv \| \leq \rups}
    \exp\bigl( - \frac{1}{2}\| \uv \|^{2} \bigr) d\uv} \, .
\label{AAuvrupsc}
\end{EQA}
Obviously, one can rewrite this value as ratio of two expectations
\begin{EQA}[c]
	\frac{\E \Bigl\{ \exp\bigl(	- \CONSTru \rups \| \TAU \gaussv \| 
			+ \frac{1}{2} \CONSTru \rups^{2}
  			+ \frac{1}{2} \| \TAU \gaussv \|^{2} \bigr) 
			\Ind\bigl( \| \TAU \gaussv \| > \rups \bigr) \Bigr\}}
		 {\P\bigl( \| \TAU \gaussv \| \leq \rups \bigr) } \, ,
\label{EindAAgelerr}
\end{EQA}
where \( \gaussv \sim \ND(0,\Id_{\dimp}) \).
Note that without the linear term \(- \CONSTru \| \TAU \gaussv \|\) in the exponent, 
the expectation in the numerator can be infinite. 
We aim at describing \( \rups\) and \( \CONSTru\)-values which ensure that 
the probability in denominator is close to one while the expectation in the numerator is small.

\begin{lemma}
\label{TGaussintext}
Let \( \TAU\) be a linear operator in \( \R^{\dimp}\) with 
\( \| \TAU \|_{\oper} \leq 1\).
Define \( \dimAA = \tr(\TAU^{\T} \TAU)\).
For any \( \CONSTru, \rups\) with 
\(1/2 < \CONSTru \leq 1\) and 
\( \CONSTru \rups = 2 \sqrt{\dimAA} + \sqrt{\xx} \) for \( \xx > 0 \) 
\begin{EQA}
  	\E \Bigl\{ \exp\Bigl(
    		- \CONSTru \rups \| \TAU \gaussv \| + \frac{\CONSTru \rups^{2}}{2}
    		+ \frac{1}{2} \| \TAU \gaussv \|^{2} 
  	\Bigr) \Ind\bigl( \| \TAU \gaussv \| > \rups \bigr) \Bigr\}
  	& \leq &
	\CONST \ex^{ - (\dimAA + \xx)/2}
\label{EexCruAg12Ag2}
\end{EQA}
and 
\begin{EQA}
    \P\bigl( \| \TAU \gaussv \| \leq \rups \bigr)
    & \geq &
    1 - \exp \Bigl\{ - \frac{1}{2} (\rups - \sqrt{\dimAA})^{2} \Bigr\}
    \geq 
    1 - \ex^{- (\dimAA + \xx)/2}.
\label{EeAg12Ag2}
\end{EQA}
\end{lemma}

\begin{remark}
The result applies even if the full dimension \( \dimp\) is infinite and \( \gaussv\) is a Gaussian element in a Hilbert space, provided that \( \dimAA = \tr(\TAU^{\T} \TAU)\) is finite, that is, \( \TAU^{\T} \TAU\) is a trace operator.
\end{remark}

\begin{proof}
Define
\begin{EQA}
	\PAAr(\rr)
	& \eqdef &
	\P\bigl( \| \TAU \gaussv \| \geq \rr \bigr),
	\\
	f(\rr)
	& \eqdef &
	\exp \Bigl( 
		- \CONSTru \rups \rr + \frac{\CONSTru \rups^{2}}{2} + \frac{\rr^{2}}{2} 
	\Bigr) .
\label{PhirPAgrfrC0}
\end{EQA}
Then
\begin{EQA}
	&& \nquad
	\E \Bigl\{ \exp\Bigl(
        - \CONSTru \rups \| \TAU \gaussv \| + \frac{\CONSTru \rups^{2}}{2}
        + \frac{1}{2} \| \TAU \gaussv \|^{2} 
      \Bigr) \Ind\bigl( \| \TAU \gaussv \| > \rups \bigr) \Bigr\}
    \\
    &=&
    - \int_{\rups}^{\infty} f(\rr) d\PAAr(\rr)
    =
    f(\rups) \PAAr(\rups)
    + \int_{\rups}^{\infty} f'(\rr) \PAAr(\rr) \, d\rr \, .
\label{irifpPhrdrf}
\end{EQA}
Now we use that
\( \PAAr\bigl( \sqrt{\dimAA} + \sqrt{2\xx} \bigr) \leq \ex^{-\xx} \)
for any \( \xx > 0 \).
This can be rewritten as
\begin{EQA}
	\PAAr(\rr)
	& \leq &
	\exp \Bigl\{ - \frac{1}{2} \bigl( \rr - \sqrt{\dimAA} \bigr)^{2} \Bigr\}
\label{Pr12rsA2mr}
\end{EQA}
for \( \rr > \sqrt{\dimAA} \).
In particular, in view of \( \rups \geq 2 \sqrt{\dimAA} + \sqrt{\xx} \)
\begin{EQA}
	f(\rups) \PAAr(\rups)
	& \leq &
	\PAAr(\rups)
	\leq 
	\exp \Bigl\{ - \frac{1}{2} \bigl( \rups - \sqrt{\dimAA} \bigr)^{2} \Bigr\}
	\leq 
	\exp\bigl\{ - \frac{1}{2} \bigl( \dimAA + \xx \bigr) \bigr\} .
\label{?}
\end{EQA}
Now we use that \( f'(\rr) = (\rr - \CONSTru \rups) f(\rr) \) and 
\begin{EQA}
	\int_{\rups}^{\infty} f'(\rr) \PAAr(\rr) \, d\rr
	&=&
	\int_{\rups}^{\infty} (\rr - \CONSTru \rups) f(\rr) \PAAr(\rr) \, d\rr
	\\
	& \leq &
	\int_{\rups}^{\infty} (\rr - \CONSTru \rups) 
		\exp\Bigl\{ - \frac{1}{2} \bigl( \rr - \sqrt{\dimAA} \bigr)^{2} 
			- \CONSTru \rups \rr + \frac{\CONSTru \rups^{2}}{2} + \frac{\rr^{2}}{2}
		\Bigr\} \, d\rr
	\\
	&=&
	\int_{\rups}^{\infty} (\rr - \CONSTru \rups) 
		\exp\Bigl\{ 
			- \bigl( \CONSTru \rups - \sqrt{\dimAA} \bigr) \rr 
			+ \frac{\CONSTru \rups^{2}}{2} - \frac{\dimAA}{2} 
		\Bigr\} \, d\rr
	\\
	&=&
	\int_{0}^{\infty} (x + \rups - \CONSTru \rups) 
		\exp\Bigl\{ 
			- \bigl( \CONSTru \rups - \sqrt{\dimAA} \bigr) (x + \rups) 
			+ \frac{\CONSTru \rups^{2}}{2} - \frac{\dimAA}{2} 
		\Bigr\} \, dx .
\label{i0iex0ixexd}
\end{EQA}
The use of \( \int_{0}^{\infty} \ex^{-x} dx = \int_{0}^{\infty} x \ex^{-x} dx = 1 \)
yields
\begin{EQA}
	\int_{\rups}^{\infty} f'(\rr) \PAAr(\rr) \, d\rr
	& \leq &
	\left( 
		\frac{\rups - \CONSTru \rups}{\CONSTru \rups - \sqrt{\dimAA}} 
		+ \frac{1}{(\CONSTru \rups - \sqrt{\dimAA})^{2}} 
	\right)
	\exp\Bigl\{ \rups \sqrt{\dimAA} 
			- \frac{\CONSTru \rups^{2}}{2} - \frac{\dimAA}{2} \Bigr\} .
\label{irifprPr222}
\end{EQA}
It remains to check that for \( \CONSTru \in (1/2,1) \)
and \( \CONSTru \rups = 2\sqrt{\dimAA} + \sqrt{\xx} \)
\begin{EQA}
	- \rups \sqrt{\dimAA} 
			+ \frac{\CONSTru \rups^{2}}{2} + \frac{\dimAA}{2}
	& \geq &
	\frac{\xx + \dimAA}{2} \, .
\label{xp2r0sQfr2}
\end{EQA}
The result follows.
\end{proof}

\def\Xset{\cc{X}}
\subsection{Taylor expansions}
\label{STaylor}
Here we collect some useful bounds for various Taylor-type expansions for
a smooth function. 
Let \(f\) be a four time differentiable function on \(\R^{\dimp}\).
Here \( \dimp \leq \infty \).
By \(f^{(m)}(\xv,\uv)\) we denote the \(m\)th directional derivative at \(\xv\):
\begin{EQA}
    f^{(m)} (\xv,\uv)
    & \eqdef &
    \frac{d^{m}}{dt^{m}} f(\xv + t \uv) \bigg|_{t=0} \, .
\label{fmxudmdtmfxtu}
\end{EQA}
In particular,
\(f'(\xv,\uv) = \bigl\langle \nabla f(\xv), \uv \bigr\rangle\) and 
\(f''(\xv,\uv) = \bigl\langle \nabla^{2} f(\xv) \, \uv, \uv \bigr\rangle\).
Below we assume that some open set \(\Xset \subseteq \R^{\dimp}\) 
is fixed, and, in addition, for each \( \xv \in \Xset \), 
and a centrally symmetric convex set \(\UV(\xv)\) are fixed and
\begin{EQA}
    \frac{1}{m!} \bigl| f^{(m)}(\xv, \uv) \bigr|
    =
    \delta_{m}(\xv,\uv)
    & \leq &
    \delta_{m} \, ,
    \quad
    \xv \in \Xset, \uv \in \UV ,
    \quad
    m=3,4
\label{1mffmxum34}
\end{EQA}
for some constants \(\delta_{m} \) depending on \( \Xset \) and \( \UV \).
All bounds will be given in terms of \(\delta_{3}\) and \(\delta_{4}\).
The construction can be extended by making \( \UV \) dependent on \( \xv \in \Xset \)
at cost of more complicated notation. 

\begin{lemma}
\label{Lextftumu}
Suppose \eqref{1mffmxum34} with \(\delta_{m} \leq 1\) for \(m=3,4\).
Then for any point \(\xv \in \Xset\) 
\begin{EQA}
    && \nquad
    \left| \frac{1}{2} \left( \ex^{f(\xv + \uv) - f(\xv) - f'(\xv,\uv)} + \ex^{f(\xv - \uv) - f(\xv) + f'(\xv,\uv)} \right) 
    - \ex^{f''(\xv,\uv)/2} \right|
    \\
    & \leq &
    \ex^{f''(\xv,\uv)/2} \,\bigl( 4 \delta_{3}^{2} + 4 \delta_{4} \bigr) .
\label{12efxufxmufpp2}
\end{EQA}
Furthermore, 
\begin{EQA}
    \left| \ex^{f(\xv + \uv) - f(\xv) - f'(\xv,\uv)} - \ex^{f''(\xv,\uv)/2} \right|
    & \leq &
    \delta_{3} \, \ex^{f''(\xv,\uv)/2} \, .
\label{12efxufxmufpp23}
\end{EQA}    
\end{lemma}

\begin{proof}
Taylor expansions of the forth order imply 
\begin{EQA}
    f(\xv + \uv) - f(\xv) - f'(\xv, \uv) - \frac{1}{2} f''(\xv, \uv) - \frac{1}{6} f^{(3)}(\xv,\uv)  
    & = &
    \rho_{1} \, ,
    \qquad
    |\rho_{1}|
    \leq 
    \delta_{4} \, ,
    \\
    f(\xv - \uv) - f(\xv) + f'(\xv, \uv) - \frac{1}{2} f''(\xv, \uv) + \frac{1}{6} f^{(3)}(\xv,\uv) 
    & = &
    \rho_{2} \, ,
    \qquad
    |\rho_{2}|
    \leq 
    \delta_{4} \, .
\label{fu12un2fud2}
\end{EQA}
Further, define \(\kappa = f^{(3)}(\xv,\uv)/6\), so that 
\(|\kappa| \leq \delta_{3} \leq 1\).
Then
\begin{EQA}
    && \nquad
    \ex^{f(\xv + \uv) - f(\xv) - f'(\xv,\uv)} 
    + \ex^{f(\xv - \uv) - f(\xv) + f'(\xv,\uv)} - 2 \ex^{f''(\xv,\uv)/2}
    \\
    &=&
    \ex^{f''(\xv,\uv)/2} \left( \ex^{\kappa + \rho_{1}} + \ex^{-\kappa + \rho_{2}} - 2 \right) .
\label{12f1fm1fpp2m2}
\end{EQA}
The function 
\begin{EQA}
      g(s)
      & \eqdef &
      \frac{1}{2} \exp\bigl( s \, \kappa + \rho_{1} \bigr) 
      + \frac{1}{2} \exp\bigl( - s \, \kappa + \rho_{2} \bigr) - 1 
\label{gt12etKd1d2}
\end{EQA}
fulfills
\begin{EQA}
    |g(0)|
    &=&
    \Bigl| \frac{1}{2} \ex^{\rho_{1}} + \frac{1}{2} \ex^{\rho_{2}} - 1 \bigr|
    \leq 
    |\rho_{1}| + |\rho_{2}| , 
    \\
    |g'(0)|
    &=&
    \frac{1}{2}
    \bigl| \kappa \bigl( \ex^{\rho_{1}} - \ex^{\rho_{2}} \bigr) \bigr|
    \leq 
    |\rho_{1}| + |\rho_{2}| 
\label{g012ed12D12K}
\end{EQA}
and for any \(s \in [0,1]\) by simple algebra due to 
\(|\kappa| \leq 1\) and \(|\rho_{m}| \leq 1\) for
\(m=1,2\)
\begin{EQA}
    |g''(s)|
    &=&
    \frac{1}{2} \Bigl| 
    \kappa^{2} \Bigl\{ \exp\bigl( s \, \kappa + \rho_{1} \bigr) 
    + \exp\bigl( - s \, \kappa + \rho_{2} \bigr) \Bigr\} 
    \Bigr|
    \\
    & \leq &
    \frac{|\kappa|^{2} \ex}{2} \, \bigl( \ex^{ |\kappa|} + \ex^{-|\kappa|} \bigr)
    < 
    8 | \kappa |^{2} ,
\label{gpptK2gt12eD}
\end{EQA}
and thus
\begin{EQA}
    \bigl| g(1) \bigr|
    & \leq &
    \sup_{s \in [0,1]} \bigl| g(0) + g'(0) + \frac{1}{2} g''(s) \bigr|
    \leq 
    4 | \kappa |^{2} + 2 |\rho_{1}| + 2 |\rho_{2}| ,
\label{g1supt01g0gp0}
\end{EQA}
and \eqref{12efxufxmufpp2} follows. 
The bound \eqref{12efxufxmufpp23} can be obtained in a similar way 
using the Taylor expansion of the third order.
\end{proof}

Now we study the modulus of continuity for the gradient \(\nabla f(\xv)\) and the Hessian
\(\nabla^{2} f(\xv)\).

\begin{lemma}
\label{LHessTay}
Suppose \eqref{1mffmxum34} with \(\delta_{3} \leq 1\).
Let \(\xv \in \Xset\) and \(\uv \in \UV\) be such that
\(\xv + \uv \in \Xset\). 
Then, for any \(\wv \in \UV\)
\begin{EQA}
\label{wTnfxun2fu}
    \Bigl| \bigl\langle  \wv, \nabla f(\xv+\uv) - \nabla f(\xv) - \nabla^{2}f(\xv) \uv \bigr \rangle 
    \Bigr|
    & \leq &
    \CONST \delta_{3} \, ,
    \\
    \Bigl| \bigl\langle \wv, \bigl\{ \nabla^{2} f(\xv+\uv) - \nabla^{2}f(\xv) \bigr\} \wv \bigr\rangle \Bigr|
    & \leq &
    \CONST \delta_{3} \, .
\label{wTn2fxun2fw}
\end{EQA}
\end{lemma}

\begin{proof}
Let us fix any \(\xvd \in \Xset\) and \( \wvd \in \UV \) and define the function 
\begin{EQA}
	g(t)
	& \eqdef &
	f(\xvd + t \wvd) + f(\xvd - t \wvd) - 2 f(\xvd) - t^{2} f''(\xvd,\wvd) .
\label{fxdtwmtv222}
\end{EQA}
The Taylor expansion of the third order yields  
\begin{EQA}
    \bigl| g(1) \bigr|
    =
    \Bigl| f(\xvd+\wvd) + f(\xvd-\wvd) - 2 f(\xvd) - f''(\xvd,\wvd) \Bigr| 
    & \leq &
    2 \delta_{3}(\xvd,\wvd) \, .
\label{fxw2fxmw2fp2}
\end{EQA}
We apply this bound for \(\xvd = \xv\) and \(\xvd = \xv + \uv\) and take the difference 
between them. 
This implies
\begin{EQA}
    &&
    \bigl| f''(\xv,\wvd) - f''(\xv + \uv,\wvd) \bigr|
    \leq 
    \bigl| f(\xv+\wvd) + f(\xv-\wvd) - 2 f(\xv) 
    \\
    &&
    \qquad
    - \, f(\xv + \uv + \wvd) - f(\xv + \uv - \wvd) + 2 f(\xv + \uv) \bigr|
    + 2 \delta_{3}(\xv,\wvd) + 2 \delta_{3}(\xv+\uv,\wvd) \, . 
    \qquad \qquad
\label{12fxw224}
\end{EQA}
For given \( \xv, \uv, \wv \), and \(\xb = \xv + \uv/2\), define 
\begin{EQA}
    g(t)
    & \eqdef &
    f\bigl( \xb + t(\uv+\wv) \bigr) - f\bigl( \xb - t(\uv+\wv) \bigr)
    \\
    &&
    + \, f\bigl( \xb + t(\uv-\wv) \bigr) - f\bigl( \xb - t(\uv-\wv) \bigr)
    - 2 f\bigl( \xb + t\uv \bigr) + 2 f\bigl( \xb - t\uv \bigr) .
\label{gtfxtuw22fxtu2}
\end{EQA}
It is straightforward to see that \(g(0) = g'(0) = g''(0) = 0\).
Moreover, in view of \(\uv \in \UV\) and
\((\uv \pm \wv)/2 \in \UV\), it holds \( \delta_{3}(\xb,\uv/2) = \delta_{3}(\xb,\uv)/8 \) and
for any \(|t| \leq 1/2\)
\begin{EQA}
    \frac{1}{6} \bigl| g^{(3)}(t) \bigr|
    & \leq &
    \frac{5 \delta_{3}}{2} \, .
\label{124g46d4}
\end{EQA}
By Taylor expansion of the third order we derive
\begin{EQA}
    \bigl| g(1/2) \bigr|
    & \leq &
    \sup_{t \in [0,1]} \frac{1}{6} \bigl| g^{(3)}(t) \bigr|
    \leq 
    \frac{5 \delta_{3}}{2} \, . 
\label{g1t01124g4t6}
\end{EQA}
Note that \(g(1/2)\) is exactly the expression in the right hand-side of \eqref{12fxw224} with \( \wvd = \wv/2 \).
The use of \( \delta_{3}(\xvd,\wvd) = \delta_{3}(\xvd,\wv)/8 \)
together with \eqref{12fxw224} yields \eqref{wTn2fxun2fw} with \(\CONST = 3\).
\end{proof}

Now we specify the result to the case of an elliptic set \(\UV\) of the form
\begin{EQA}
    \UV 
    &=& 
    \bigl\{ \uv \colon \| \QP \uv \| \leq \rr \bigr\} 
\label{UVdefQPur}
\end{EQA}
for a positive invertible operator \(\QP\) and \(\rr > 0\).

\begin{lemma}
\label{LellUVD2w}
Let \(\UV\) be given by \eqref{UVdefQPur} with \(\QP > 0\), and let \(\xv \in \Xset\) and 
\(\uv \in \UV\) be such that \(\xv + \uv \in \Xset\).
Then
\begin{EQA}
\label{Qm1n2fxn2txurm1}
    \bigl\| 
    \QP^{-1} \bigl\{ \nabla f(\xv + \uv) - \nabla f(\xv) - \nabla^{2} f(\xv) \uv \bigr\} 
    \bigr\|
    & \leq &
    \CONST \rr^{-1} \delta_{3} \, ,
    \\
    \bigl\| \QP^{-1} \bigl\{ \nabla^{2} f(\xv) - \nabla^{2} f(\xv + \uv) \bigr\} \QP^{-1} \bigr\|
    & \leq &
    \CONST \rr^{-2} \delta_{3} \, .
\label{Qm1n2fxn2txurm2}
\end{EQA}
\end{lemma}

\begin{proof}
For any \(\wv \in \UV\), it holds by Lemma~\ref{LHessTay}
\begin{EQA}
    \Bigl| 
    	\bigl\langle \wv, \bigl\{ \nabla^{2} f(\xv+\uv) - \nabla^{2} f(\xv) \bigr\} \wv \bigr\rangle 
    \Bigr|
    & = &
    \Bigl| \bigl\langle \QP \wv, \QP^{-1} \bigl\{ \nabla^{2} f(\xv+\uv) - \nabla^{2} f(\xv) \bigr\} 
    \QP^{-1} (\QP \wv) \bigr\rangle \Bigr|
    \leq 
    \CONST \delta_{3} \, .
\label{wTn2fxun2fw3}
\end{EQA}
As this bound holds for all \(\wv \in \UV\) with \(\| \QP \wv \| \leq \rr\), the result follows.
\end{proof}

The result of Lemma~\ref{Lextftumu} can be extended to the integral of \(\ex^{f(\xv + \uv)}\)
over \(\uv \in \UV\).

\begin{lemma}
\label{Lintfxupp2}
Let \( \UV \) be a subset in \( \R^{\dimp} \).
Suppose \eqref{1mffmxum34} with \(\delta_{m} \leq 1\) for \(m=3,4\).
Then for any point \(\xv \in \Xset\) and any centrally symmetric set \( A \subset \UV \)
\begin{EQA}
\label{4d324d4efppm}      
    \biggl| 
    	\int_{A} \ex^{f(\xv+\uv) - f(\xv) - f'(\xv,\uv)} \, d\uv 
    	- \int_{A} \ex^{f''(\xv,\uv)/2} \, d\uv \biggr|
    & \leq &
    \err \int_{A} \ex^{f''(\xv,\uv)/2} \, d\uv \, 
    \qquad
\end{EQA}
with \( \err = 4 \delta_{3}^{2} + 4 \delta_{4} \). 
If \( A \) is not centrally symmetric then
\begin{EQA}
\label{d324d4efppm3}      
    \biggl| \int_{A} \ex^{f(\xv+\uv) - f(\xv) - f'(\xv,\uv)} \, d\uv - \int_{A} \ex^{f''(\xv,\uv)/2} \, d\uv \biggr|
    & \leq &
    \delta_{3} \int_{A} \ex^{f''(\xv,\uv)/2} \, d\uv \, .
    \qquad
\end{EQA}
\end{lemma}

\begin{proof}
By symmetricity of \(\UV\), it holds 
\begin{EQA}
    \int_{A} \ex^{f(\xv+\uv) - f(\xv) - f'(\xv,\uv)} \, d\uv
    &=&
    \frac{1}{2} \int_{A}\left( \ex^{f(\xv+\uv) - f(\xv) - f'(\xv,\uv)} 
    + \ex^{f(\xv-\uv) - f(\xv) + f'(\xv,\uv)} \right) d\uv ,
\label{12efuAefmuA}
\end{EQA}
and the first result is proved by \eqref{12efxufxmufpp2}.
The final bound for any \( A \) follows from \eqref{12efxufxmufpp23}.
\end{proof}

\subsection{Deviation bounds for Gaussian quadratic forms}
\label{SdevboundGauss}
The next result explains the concentration effect of 
\( \langle \BB \gaussv, \gaussv \rangle \)
for a standard Gaussian vector \(\gaussv\) and a symmetric trace operator \(\BB\) in \( \R^{\dimp} \),
\( \dimp \leq \infty \).
We use a version from \cite{laurentmassart2000}.

\begin{theorem}
\label{TexpbLGA}
\label{Lxiv2LD}
\label{Cuvepsuv0}
Let \(\gaussv\) be a standard normal Gaussian element in \( \R^{\dimp} \) and \(\BB\) be symmetric non-negative trace operator in \( \R^{\dimp}\).
Then with \(\dimA = \tr(\BB)\), \(\vA^{2} = \tr(\BB^{2})\), and 
\(\supA = \| \BB \|_{\oper}\), it holds for each \(\xx \geq 0\)
\begin{EQA}
\label{Pxiv2dimAvp12}
	&&\P\Bigl( \langle \BB \gaussv, \gaussv \rangle > \zq^{2}(\BB,\xx) \Bigr)
	\leq
	\ex^{-\xx} ,
	\\
	&& \zq(\BB,\xx)
	\eqdef
	\sqrt{\dimA + 2 \vA \xx^{1/2} + 2 \supA \xx} \,\, .
\label{zqdefGQF}
\end{EQA}
In particular, it implies 
\begin{EQA}
	\P\bigl( \| \BB^{1/2} \gaussv \| > \dimA^{1/2} + (2 \supA \xx)^{1/2} \bigr)
	& \leq &
	\ex^{-\xx} .
\label{Pxiv2dimAxx12}
\end{EQA}
Also
\begin{EQA}
	\P\bigl( \langle \BB \gaussv, \gaussv \rangle < \dimA - 2 \vA \xx^{1/2} \bigr)
	& \leq &
	\ex^{-\xx} .
\label{Pxiv2dimAvp12m}
\end{EQA}
If \(\BB\) is symmetric but non necessarily positive then
\begin{EQA}
	\P\bigl( \bigl| \langle \BB \gaussv, \gaussv \rangle - \dimA \bigr| > 2 \vA \xx^{1/2} + 2 \supA \xx \bigr)
	& \leq &
	2 \ex^{-\xx} .
\label{PxivTBBdimA2vp}
\end{EQA}
\end{theorem}

As a special case, we present a bound for the chi-squared distribution 
corresponding to \(\BB = \Id_{\dimp}\), \( \dimp < \infty \).
Then \(\tr (\BB) = \dimp\), \(\tr(\BB^{2}) = \dimp\) and \(\supA(\BB) = 1\).

\begin{corollary}
\label{Cchi2p}
Let \(\gaussv\) be a standard normal vector in \(\R^{\dimp}\).
Then for any \( \xx > 0 \)
\begin{EQA}[ccl]
\label{Pxi2pm2px}
	\P\bigl( \| \gaussv \|^{2} \geq \dimp + 2 \sqrt{\dimp \xx} + 2 \xx \bigr)
	& \leq &
	\ex^{-\xx},
	\\
	\P\bigl( \| \gaussv \| \,\,  \geq \sqrt{\dimp} + \sqrt{2 \xx} \bigr)
	& \leq &
	\ex^{-\xx} ,
\label{Pxi2pm2px12}
	\\
	\P\bigl( \| \gaussv \|^{2} \leq \dimp - 2 \sqrt{\dimp \xx} \bigr)
	& \leq &
	\ex^{-\xx}	.
\label{Pxi2pm2px22}
\end{EQA}
\end{corollary}

\subsection{Deviation bounds for non-Gaussian quadratic forms}
\label{Sprobabquad}
This section collects some probability bounds for non-Gaussian quadratic forms.
The presented results can be viewed as a slight improvement of the bounds from \cite{SP2011}.
The proofs are very similar to ones from \cite{SP2011} and are omitted by the space reasons. 

\label{SdevboundnonGauss}

Let a random vector \(\xiv \in \R^{\dimp}\) has some exponential moments.
More exactly, suppose  for some fixed \(\gm > 0\) that 
\begin{EQA}[c]
    \log \E \exp\bigl( \langle \gammav, \xiv \rangle \bigr)
    \le
    \| \gammav \|^{2}/2,
    \qquad
    \gammav \in \R^{\dimp}, \, \| \gammav \| \le \gm .
\label{expgamgm}
\end{EQA}
First we present a bound for the norm \( \| \xiv \| \) assuming \( \dimp \lesssim \gm^{2} \).
For ease of presentation, assume below that \(\gm\) is sufficiently large, namely, 
\(0.3 \gm \ge \sqrt{\dimp}\).
In typical examples of an i.i.d. sample, \(\gm \asymp \sqrt{n}\).
Define
\begin{EQA}
	\xxc
	& \eqdef &
	\gm^{2}/4,
	\\
	\zqc^{2}
	& \eqdef &
	\dimp + \sqrt{\dimp \gm^{2}} + \gm^{2}/2 
	=
	\gm^{2} \bigl( 1/2 + \sqrt{\dimp/\gm^{2}} + \dimp/\gm^{2} \bigr),
	\\
	\gmc
	& \eqdef &
	\frac{\gm \, \bigl( 1/2 + \sqrt{\dimp/\gm^{2}} + \dimp/\gm^{2} \bigr)^{1/2}}
		 {1 + \sqrt{\dimp/\gm^{2}}} .
\label{xxcgm24mucyyc2}
\end{EQA}
Note that with \(\alp = \sqrt{\dimp / \gm^{2}} \leq 0.3\), one has
\begin{EQA}
	\zqc^{2}
	&=&
	\gm^{2} \bigl( 1/2 + \alp + \alp^{2} \bigr),
	\quad
	\gmc
	=
	\gm \,\, \frac{\bigl( 1/2 + \alp + \alp^{2} \bigr)^{1/2}}{1 + \alp} \,\, ,
\label{zqcgmcalp2}
\end{EQA}
so that \(\zqc^{2} / \gm^{2} \in [1/2,1]\) and \(\gmc^{2} / \gm^{2} \in [1/2,1]\).

\begin{theorem}
\label{LLbrevelocro}   
Let \eqref{expgamgm} hold and 
\(0.3 \gm \ge \sqrt{\dimp}\).
Then for each \(\xx > 0\)
\begin{EQA}
    \P\bigl( \| \xiv \| \ge \zq(\dimp,\xx) \bigr)
    & \le &
    2 \ex^{-\xx} + 8.4 \ex^{-\xxc } \Ind(\xx < \xxc) ,
\label{PxivbzzBBro}
\end{EQA}    
where \(\zq(\dimp,\xx)\) is defined by
\begin{EQA}
\label{PzzxxpBro}
    \zq(\dimp,\xx)
    & \eqdef &
    \begin{cases}
      \bigl( \dimp + 2 \sqrt{\dimp \xx} + 2 \xx\bigr)^{1/2}, &  \xx \le \xxc  , \\
      \zqc + 2 \gmc^{-1} (\xx - \xxc)   , & \xx > \xxc .
    \end{cases}
\label{zzxxppdBlro}
\end{EQA}    
\end{theorem}

Depending on the value \(\xx\), we have two types of tail behavior of the 
quadratic form \(\| \xiv \|^{2}\). 
For \(\xx \le \xxc = \gm^{2}/4\), we have the same deviation bounds as in the Gaussian case
with the extra-factor two in the deviation probability.
Remind that one can use a simplified expression 
\(\bigl( \dimp + 2 \sqrt{\dimp \xx} + 2 \xx\bigr)^{1/2} \leq \sqrt{\dimp} + \sqrt{2 \xx}\).
For \(\xx > \xxc\), we switch to the special regime driven by the exponential moment
condition \eqref{expgamgm}.
Usually \(\gm^{2}\) is a large number (of order \(n\) in the i.i.d. setup) and 
the second term in \eqref{PxivbzzBBro} can be simply ignored. 
The result applies with \( \gm = \infty \) yielding the simplified bound 
\begin{EQA}
    \P\bigl( \| \xiv \| \ge \zq(\dimp,\xx) \bigr)
    & \le &
    2 \ex^{-\xx} ,
    \\
    \zq(\dimp,\xx) 
    &=& 
    \sqrt{\dimA + 2 \dimA \xx^{1/2} + 2 \xx }
    \leq 
    \sqrt{\dimA} + \sqrt{2 \xx} \, .
\label{PxivbzzBBroinf}
\end{EQA}    

\medskip

Next we present a bound for a quadratic form \( \langle \BB \xiv, \xiv \rangle \), where 
\(\xiv\) satisfies \eqref{expgamgm} and \(\BB\) is a given symmetric non-negative 
operator in \( \R^{\dimp} \).
Here we relax \( \dimp < \infty \) to \( \tr \BB < \infty \).
Define 
\begin{EQA}[c]
    \dimB
    \eqdef
    \tr \bigl( \BB \bigr) ,
    \qquad 
    \vpB^{2}
    \eqdef
    \tr(\BB^{2}) ,
    \qquad
    \lambdaB \eqdef \lambda_{\max}\bigl( \BB \bigr). 
\label{BBrddB}
\end{EQA}   
For ease of presentation, 
suppose that \(0.3 \gm \ge \sqrt{\dimA}\) so that 
\(\alp = \sqrt{\dimA / \gmb^{2}} \leq 0.3\).
The other case only changes the constants in the inequalities. 
Define also
\begin{EQA}
	\xxc
	& \eqdef &
	\gm^{2}/4,
	\\
	\zqc^{2}
	& \eqdef &
	\dimA + \vp \gm + \supA \gmb^{2}/2 ,
	\\
	\gmc
	& \eqdef &
	\frac{ \sqrt{\dimA/\supA + \gm \vp / \supA + \gm^{2}/2}}{1 + \vp / (\supA \gmb)} .
\label{xxcgm24mucyyc2B}
\end{EQA}

\begin{theorem}
\label{LLbrevelocroB}   
Let \eqref{expgamgm} hold and 
\(0.3 \gm \ge \sqrt{\dimA/\supA}\).
Then for each \(\xx > 0\)
\begin{EQA}
    \P\bigl( \langle \BB \xiv, \xiv \rangle \ge \zq^{2}(\BB,\xx) \bigr)
    & \le &
    2 \ex^{-\xx} + 8.4 \ex^{-\xxc} \Ind(\xx < \xxc) ,
\label{PxivbzzBBroB}
\end{EQA}    
where \(\zq(\BB,\xx)\) is defined by
\begin{EQA}
\label{PzzxxpBroB}
    \zq(\BB,\xx)
    & \eqdef &
    \begin{cases}
      \sqrt{ \dimA + 2 \vp \xx^{1/2} + 2 \supA \xx }, &  \xx \le \xxc, \\
      \zqc + 2 \lambdaB (\xx - \xxc)/\gmc , & \xx > \xxc.
    \end{cases}
\label{zzxxppdBlroB}
\end{EQA}    
If \eqref{expgamgm} hold with \( \gm = \infty \), then
\begin{EQA}
    \P\bigl( \langle \BB \xiv, \xiv \rangle \ge \zq^{2}(\BB,\xx) \bigr)
    & \le &
    2 \ex^{-\xx}  ,
    \quad
    \zq(\BB,\xx)
    =
    \sqrt{ \dimA + 2 \vp \xx^{1/2} + 2 \supA \xx }
    \leq 
    \sqrt{\dimA} + \sqrt{2 \supA \xx} \, .
    \qquad
\label{PxivbzzBBroBinf}
\end{EQA}    

\end{theorem}

Similarly to the case \(\BB = \Id_{\dimp}\), the upper quantile 
\(\zq(\BB,\xx) = \sqrt{ \dimA + 2 \vp \xx^{1/2} + 2 \supA \xx }\) can be upper bounded 
by \(\sqrt{\dimA} + \sqrt{2 \supA \xx}\):
\begin{EQA}
\label{PzzxxpBroBu}
    \zq(\BB,\xx)
    & \leq &
    \begin{cases}
        \sqrt{\dimA} + \sqrt{2 \supA \xx}, &  \xx \le \xxc, \\
        \zqc + 2 \lambdaB (\xx - \xxc)/\gmc , & \xx > \xxc.
    \end{cases}
\end{EQA}    

\subsection{Gaussian comparison}
\label{SGausscomp}
Let \(\HM\) be a Hilbert space and 
\(\Sigma_{\xiv}\) be a covariance operator  of an arbitrary Gaussian random element in \(\HM\). 
By \(\{\lambda_{k\xiv}\}_{k \geq 1}\) we denote the set of its eigenvalues arranged in the non-increasing order, i.e. 
\(\lambda_{1\xiv} \geq \lambda_{2\xiv} \geq \ldots \), and let 
\(\lambdav_{\xiv} \eqdef \diag(\lambda_{j \xiv})_{j=1}^{\infty}\). 
Note that  \(\sum_{j=1}^{\infty} \lambda_{j \xiv} < \infty\).  
Introduce the following quantities
\begin{EQA}[c]
    \Frobg_{k\xiv}^{2} 
    \eqdef 
    \sum_{j=k}^{\infty} \lambda_{j \xiv}^{2}, \quad k = 1,2,
\label{Lambda def}
\end{EQA}

\begin{theorem}[\cite{GNSUl2017}]
\label{Tgaussiancomparison3}
Let \(\xiv\) and \(\etav\) be Gaussian elements in \(\HM\) with zero mean and covariance operators 
\(\Sigma_{\xiv}\) and \(\Sigma_{\etav}\) respectively.
Then for any \(\av \in \HM\)
\begin{EQA}[rcl]
    && \nquad
    \sup_{x > 0} \left|\P( \| \xiv - \av \| \leq x) - \P( \| \etav  \| \leq x) \right| 
    \\ 
    &&
    \lesssim  
    \bigg(\frac{1}{(\Frobg_{1\xiv}\Frobg_{2\xiv})^{1/2}} + \frac{1}{(\Frobg_{1\etav}\Frobg_{2\etav})^{1/2}}\bigg) 
    \bigg( \| \lambdav_{\xiv} - \lambdav_{\etav}  \|_{1} + \| \av\|^{2}\bigg).
\label{expl_gauss22}
\end{EQA}
Moreover, assume that
\begin{EQA}[c]
    3 \| \Sigma_{\xiv}\|^{2} \le \| \Sigma_{\xiv}\|_{\Fr}^{2} 
    \quad \text{ and }  \quad 
    3 \| \Sigma_{\etav}\|^{2} \le \| \Sigma_{\etav}\|_{\Fr}^{2} \, .
\label{2plusdelta}
\end{EQA}
Then for any \(\av \in \HM\) 
\begin{EQA}
    \sup_{x > 0} \left|\P(\| \xiv - \av \| \leq x) - \P(\| \etav  \| \leq x) \right| 
    & \lesssim &
    \biggl(
    \frac{1}{\| \Sigma_{\xiv}\|_{\Fr}} 
    + \frac{1}{\| \Sigma_{\etav}\|_{\Fr}}
    \biggr) 
    \biggl( \| \lambdav_{\xiv} - \lambdav_{\etav} \|_{1} + \| \av\|^{2}\biggr).
\label{expl_gauss 2}
\end{EQA}
\end{theorem}

\bibliography{exp_ts,listpubm-with-url}

\end{document}